\documentclass[a4paper,12pt]{article}
\usepackage{amsfonts,amsmath,amssymb,amsthm}
\usepackage{txfonts}
\usepackage[utf8]{inputenc}
\usepackage{amsmath}
\usepackage{amssymb}
\usepackage{mathtools}
\usepackage{bbm}
\usepackage{float}

\RequirePackage[colorlinks,citecolor=red,urlcolor=blue, linkcolor=blue]{hyperref}
\usepackage[margin=0.8in]{geometry}

\usepackage[dvips]{graphics}
\usepackage{epsfig,rotating}
\usepackage{tabularx}
\usepackage{bm}
\usepackage{bbm}
\usepackage{array}
\newcolumntype{P}[1]{>{\raggedright\arraybackslash}p{#1}}

\usepackage{mdframed}
\usepackage{lipsum}
\usepackage{enumitem}

\usepackage[authoryear,square]{natbib}

\usepackage[nottoc,numbib]{tocbibind} 
\usepackage[toc,page]{appendix}

\newtheorem{thm}{Theorem}
\newtheorem{lemma}{Lemma}
\newtheorem{remark}{Remark}
\newtheorem{prop}{Proposition}

\newtheorem{cond}{Condition}
\newtheorem{class}{Class}

\DeclareMathOperator{\Var}{Var}

\newcommand{\beaa}{\begin{eqnarray*}}
\newcommand{\eeaa}{\end{eqnarray*}}
\newcommand{\bea}{\begin{eqnarray}}
\newcommand{\eea}{\end{eqnarray}}

\newcommand{\la}{\left\{}
\newcommand{\ra}{\right\}}
\newcommand{\lb}{\left(}
\newcommand{\rb}{\right)}
\newcommand{\mb}{\mathbb}
\newcommand{\ve}{\varepsilon}
\newcommand{\mc}{\mathcal}
\begin{document}
\title{Detecting non-uniform patterns on high-dimensional hyperspheres}
\author{Tiefeng Jiang \\ \small{School of Data Science, Chinese University of Hong Kong, Shenzhen} \\ \small{\href{mailto:jiang040@cuhk.edu.cn}{jiang040@cuhk.edu.cn}}
   \and Tuan Pham \\ \small{Department of Statistics and Data Science, University of Texas, Austin} \\  \small{\href{mailto:tuan.pham@utexas.edu}{tuan.pham@utexas.edu} }}
\date{}

\maketitle

\begin{abstract}

We propose a new probabilistic characterization of the uniform distribution on the hypersphere in terms of the distribution of pairwise inner products, extending the ideas of \citep{cuesta2009projection,cuesta2007sharp} in a data-driven manner. This characterization naturally leads to an Ingster-type distance for quantifying deviations from uniformity, whose asymptotic behavior can be analyzed systematically via Edgeworth-type expansions. Perhaps surprisingly, we show that this distance captures the minimax rates for testing uniformity simultaneously across several high-dimensional parametric models, even in the models where densities with respect to the uniform law do not exist.

We then introduce a simple test for spherical uniformity based on this distance and study its detection rates and consistency against various classes of alternatives, both local and non-local. The proposed test is universally consistent in fixed dimensions, minimax-optimal over a variety of high-dimensional parametric models, and consistent against non-local high-dimensional alternatives. This is different from previously studied high-dimensional Sobolev tests and extreme-value-based tests, which are rate-suboptimal or inconsistent against one or more classes of alternatives. We also establish the local asymptotic distribution of the proposed test under the considered classes of alternatives, along with new information lower bounds.

\end{abstract}




\tableofcontents

\section{Introduction}

Testing whether a sample from an unknown distribution is uniformly distributed over a domain is a classical problem in statistical theory. In the discrete case, this problem has been extensively studied by statisticians, computer scientists, and probabilists; see \cite{bhattacharya2024sparse,balakrishnan2018hypothesis} and references therein. For the continuous case, one of the most common and intriguing settings is the unit hypersphere, not only because of its rich mathematical structure but also due to its importance in statistical analysis on non-Euclidean spaces. Below, we briefly formulate the problem and review relevant literature.

 Consider the hypersphere $\mb{S}^{p-1}= \la x \in \mb{R}^{p}: \|x\|_2 =1 \ra$, where $\|.\|_2$ is the Euclidean distance. The observed data points are denoted by $\bm{X}_1, \bm{X}_2,...,\bm{X}_n$ with $\bm{X}_i \in \mb{S}^{p-1}$ for all $i=1,2,...,n$. We are mainly interested in the high-dimensional case where one assumes $p=p_n$ is a sequence diverging to infinity. Assume that the data $\bm{X}_i$'s are drawn independently from an unknown distribution $\mu$ supported on the hypersphere $\mb{S}^{p-1}$.  The uniform distribution on the hypersphere $\mb{S}^{p-1}$ is denoted by $\mbox{Unif}(\mb{S}^{p-1})$. The uniformity testing problem can be formulated as
\begin{align} \label{uniform-test}
    H_0: \mu= \mbox{Unif}(\mb{S}^{p-1}) \ \ \ \mbox{against}\ \ \ H_{1}:\mu \neq \mbox{Unif}(\mb{S}^{p-1}).
\end{align} 

Before going further, let us explain why the continuous case, and especially the spherical setting, is more interesting than the discrete one. The sphere carries a rich geometric structure that discrete sample spaces do not. While testing on discrete spaces is essentially a multinomial goodness-of-fit problem based on counts, testing on $\mathbb{S}^{p-1}$ involves geometry, rotation invariance, and distinctive high-dimensional effects such as the curse of dimensionality. As a result, departures from uniformity on the sphere can take many qualitatively different forms, such as concentration around a direction, axial structure, or support near a lower-dimensional subset. This makes the problem both mathematically richer and more closely connected to modern applications involving directional and high-dimensional data.

Moreover, unlike in the discrete setting, on high-dimensional spheres there is a distinction between {\it local} and {\it non-local} alternatives. One may view local alternatives as those arising from smoothly parametrized models containing the null, whereas non-local alternatives do not arise from such local perturbations and may remain separated from the null in a stronger sense. It is therefore of interest not only to detect local alternatives with high power, but also to develop tests that remain consistent against non-local ones. Perhaps surprisingly, in high dimensions, procedures that are effective against local alternatives may fail completely against non-local alternatives, and vice versa. More generally, it remains poorly understood how high dimensionality influences the minimax rates in the testing problem \eqref{uniform-test}.


In fixed and small dimensions, the uniformity testing problem has been investigated extensively in the last few decades. An incomplete list of early results concerning the case $p=2$ includes the Kuiper test (see, for example, \cite{Kuiper}), Watson test (see, for example, \cite{watson}) and Hodges--Ajne test (see, for example, \cite{Ajne}). 
In arbitrary but fixed dimensions, the class of Sobolev-based tests was introduced in \cite{Gine} and shown to be universally consistent against any absolutely continuous alternative with $L^2$-integrable densities. Notable developments of Sobolev tests include the data-driven procedures proposed by \cite{Bogdan} and \cite{Jupp}. The readers are referred to the survey papers \cite{survey-uni,pewsey2021recent} for recent progress on this problem and for a list of recent testing procedures. The consistency and optimality of various testing procedures in fixed dimensions have been well-studied in the literature, and can also be found in \cite{survey-uni}. Recent results in the fixed-dimensional settings include \cite{garcia2021cramer,garcia2023projection,fernandez2023new,boucher2025modified,boucher2025runs}.

In the era of big data, there has been an increasing interest in studying high-dimensional directional statistics, which assumes that the dimensions diverge to infinity. For example, in shape analysis and nonparametric statistics, a popular approach is to consider sign-based procedures, in which one projects the observations onto the hyperspheres and carries out statistical inference based on the projected data. This approach is robust in high dimensions since the concentration of measure phenomenon implies that the majority of information from the data is captured by the directions rather than the magnitudes of the observations. We briefly review the high-dimensional directional statistics literature below.

 In \cite{Dryden05}, the author investigates the asymptotic properties of high-dimensional spherical distributions and their applications to brain shape modeling. Specifically, the study involved statistical modeling of a sample of \(n = 74\) MRI images of adult brains. After normalization, each brain image was represented as a unit vector with dimension \(p = 62,501\). A natural question in this modeling task is whether some simple, well-known distributions (such as the uniform distribution) provide a good fit for the data. 
Clustering analysis on high-dimensional hyperspheres has been studied in \cite{Banerjee04,Banerjee03}. Potential applications of high-dimensional uniformity tests were illustrated in \cite{Juan2001}, in which the authors relate the multivariate outliers detection problem to the uniformity testing problem. Sign-based procedures in high dimensions have been considered in \cite{Zou14} in the context of sphericity testing and  in \cite{WPL15}, where the authors propose a high-dimensional nonparametric mean test. 

Another motivation for studying the high-dimensional analog of \eqref{uniform-test} arises from deep learning theory. In overparameterized neural networks, regularization is crucial for preventing overfitting and improving generalization. In \cite{xie2017diverse}, it was shown that optimizing one-hidden-layer neural networks with approximately uniformly distributed neurons can help avoid spurious local minima. Furthermore, empirical studies in \cite{lin2020regularizing,liu2018learning} have shown that regularization methods promoting uniformity among neurons effectively reduce the generalization error in deep networks. Such methods are fundamentally tied to the question of whether a random set of points on the unit hypersphere is approximately uniformly distributed. The overparameterized nature of deep networks makes it natural to study this question in the high-dimensional settings. Given the complexity of many deep networks, including heavy-tailed or strongly correlated structures (see \cite{mahoney2019traditional} for details), we focus on detecting non-uniformity in a nonparametric manner. This perspective shifts the attention away from traditional parametric modeling goals—such as optimality and asymptotic local power within a parametric class of distributions—towards prioritizing simplicity of implementation and universal consistency.

\subsection{Related literature}
Despite the vast literature on fixed-dimensional tests, much less is known about the uniformity testing problem in the high-dimensional setting with diverging dimensions. Most approaches in the literature are designed to detect only specific types of alternatives, and those that are universally consistent are difficult to analyze in high dimensions. To the best of the authors' knowledge, only a few high-dimensional tests have been investigated in the literature. We give a brief overview of these tests below.

\begin{enumerate}
    \item {\it Rayleigh test in \cite{Cutting-P-V} and \cite{Ley-P}}. This test can be formulated in terms of a U-statistic of the data points with the inner product kernel, i.e.
     \begin{align}
          R_n &:= \frac{\sqrt{2p}}{n} \sum_{1 \leq i<j \leq n} \bm{X}^{\top}_i \bm{X}_j. \label{Rayleigh}
     \end{align}
    \item {\it Bingham test in \cite{Cutting-P-V2, Zou14} and \cite{Ley-P}}. This test is also based on a U-statistic of the data points, but with a quadratic inner product kernel, i.e. 
    \begin{align}
            B_n &:= \frac{p}{n} \sum_{1 \leq i<j \leq n} \Big[ \lb \bm{X}^{\top}_i \bm{X}_j \rb^2 - \frac{1}{p}  \Big]. \label{Bingham}
    \end{align}

    \item {\it Packing test in \cite{Jiang13}.} This test is based on the smallest angle, i.e.
    \begin{align}
            P_n &:= p \cdot  \max_{1 \leq i<j \leq n} \lb \bm{X}^{\top}_i \bm{X}_j \rb^2 -4 \log n + \log \log n \label{packing}.
    \end{align}

    \item {\it High-dimensional Sobolev tests in \cite{ebner2025high}.} A family of high-dimensional Sobolev tests that generalize \eqref{Rayleigh} and \eqref{Bingham} are proposed. 
\end{enumerate}

The null asymptotic distributions of these test statistics, as well as certain non-null results, have been studied rigorously in recent years; see, for example, \cite{Cutting-P-V,Cutting-P-V2,Ley-P,Ley-P-2,Ley-P-V,ebner2025high}.   

It is known that the Rayleigh test $R_n$ is optimal (in the sense of Le Cam) against the class of Fisher--von Mises--Langevin (FvML) distributions (which correspond to $f(x)=x$ in \eqref{semiparametric model} below), and that it is blind under the class of Watson distributions (which correspond to $f(x)=x^2$ in \eqref{semiparametric model} below), regardless of the signal strength; see \cite{Cutting-P-V,Cutting-P-V2} for details. On the other hand, the Bingham test \eqref{Bingham} achieves a sub-optimal rate for testing uniformity within the class of FvML distributions, but is optimal for testing uniformity within the class of Watson distributions. Regarding the packing test $P_n$, it is known that, under the null hypothesis and the mild assumption $p \gg (\log n)^2$, $P_n$ converges in distribution to the Gumbel distribution with CDF $\exp \lb -(8\pi)^{-1/2} e^{-x/2} \rb$; see \cite{Jiang13} and also \cite{Jiang12}. The detection threshold and consistency properties of the Packing test remain open in many models, although simulation evidence suggests that it is often sub-optimal.

Each of these three tests has its own advantages and disadvantages. 
However, a common limitation is that they are each optimal only for a specific class of (parametric) alternatives: 
they perform well against certain models but may be essentially powerless outside those classes. We note that the results in \cite{ebner2025high} do not provide high-dimensional non-null theory beyond the class of FvML distributions.

\subsection{Our contribution}

The primary objective of this article is to approach \eqref{uniform-test} from a probabilistic and geometric perspective, rather than relying on the likelihood-based framework commonly used in statistics. 
We introduce a novel \emph{pseudometric} to quantify deviations from uniformity (see \eqref{metric} below). 
Unlike most classical distances between probability measures, this distance takes into account geometric deviations from uniformity; see Section~\ref{when to use} for further discussion. The key message is that, although this distance is geometrically defined, it is able to capture the minimax rates across various high-dimensional parametric models for the problem \eqref{uniform-test} simultaneously, while also remaining consistent against non-local alternatives. Consequently, the resulting test enjoys these high-dimensional optimality and consistency properties, and remains universally consistent in fixed dimensions. To the best of our knowledge, no previous procedure has been shown to enjoy all of these properties simultaneously. Our contributions can be summarized as follows.
\begin{itemize}
    \item We propose a new distance $d$ (see \eqref{metric} for a precise definition) to quantify deviations from uniformity. 
    This distance does not require the alternatives to be absolutely continuous with respect to the uniform distribution and is therefore well suited for analyzing singular alternatives.
    A natural test statistic associated with this distance (see $T_n$ in \eqref{def-T}) is introduced and shown to converge in distribution to the supremum of a Brownian bridge under the null (Theorem \ref{null-dist}). 

    \item We prove an information lower bound of order $1/n$ for testing with respect to $d$ (Theorem \ref{information bound}), and we show that the proposed test achieves this lower bound (Theorem \ref{consistency}). 
    Both results are established without imposing any restriction on the relationship between the sample size $n$ and the dimension $p$, and the setting is completely nonparametric. We also propose a general framework based on Edgeworth-type expansions for obtaining the asymptotic behavior of $d$ under parametric models (Section \ref{sec:minimax}). 

    \item {\it \textbf{ Local alternatives}.} We investigate how the distance $d$ behaves when restricted to high-dimensional parametric models, and how it relates to the minimax rates in those settings. 
    In particular, we study three concrete models: the FvML model, the Watson model, and a low-dimensional uniform distribution model \eqref{low rank uniform}. The low-dimensional model is new and has not been studied in the literature before. We show that the distance $d$ captures precisely the minimax rates in these three models (Propositions \ref{Kolmogorov expansion fvml}, \ref{watson d asymptotic} and \ref{Kolmogorov low-rank}). This means that the threshold at which $d \asymp 1/n$ matches the minimax rate for testing uniformity in the corresponding parametric model. This phenomenon seems to hold for other models as well; see the discussion in Section \ref{when to use}.

    \item {\it \textbf{Local limiting distribution}.} For each of the models considered, we show that the proposed test achieves the minimax rate in the corresponding model and converges weakly to shifted Brownian bridges (Propositions \ref{local power FvML}, \ref{power low-rank} and Theorem \ref{watson local distribution}), where the shifts are smooth functions that vanish at the end points.

    \item {\it \textbf{Non-local alternatives}.} We investigate two models of non-local alternatives: one involving projections of heavy-tailed distributions and one involving deterministic sets of points whose moments match those of the uniform distribution; see Section \ref{sec: non-local} for more details. We show that the existing tests fail to be consistent under these non-local alternatives, while the proposed test remains consistent. 
    
\item {\it \textbf{Technical contribution}.} We obtain an information-theoretic lower bound for testing against low-dimensional uniform distributions (Theorem \ref{low-rank information bound}). This result is new and of independent interest. We also improve the information lower bound for the Watson model in the regime $p \gg n$ and show that the detection boundary undergoes a phase transition: in the previously studied regime $p=o(n)$, the detection boundary is of order $\Theta\!\left(p^{3/2}/\sqrt{n}\right)$ \cite{Cutting-P-V2}, whereas in the high-dimensional regime $n=o(p)$, it is of order $p/2-\Theta(\sqrt{pn})$; see Theorem \ref{watson minimax}. Various technical tools are also developed to study the limiting distribution under local alternatives in the Watson model, where Le Cam's third lemma does not apply directly because a LAN expansion is not available.

\end{itemize}
We would also like to point out that there are other approaches in the literature that are not based on likelihood inference, such as the family of Sobolev tests originally proposed in \cite{Gine} and the projection-based tests introduced in \cite{cuesta2009projection}, with further developments in recent works \cite{garcia2021cramer,garcia2023projection,fernandez2023new,ebner2025high}. Each of these approaches relies on a different characterization of the uniform distribution: Sobolev tests are based on the eigenfunctions of the Laplacian, while projection-based tests are based on the one-dimensional distributions obtained by projecting the data onto all possible directions. However, the detection rates of these tests remain largely unexplored in high dimensions, even under very simple models (see Table \ref{tab:detection-rates} for some known results and comparisons). Even in fixed dimensions, many of these procedures are not universally consistent, since they are based on Sobolev tests, and it is known that universal consistency of Sobolev-based tests requires the associated weight sequence to be infinite. The recent work \cite{ebner2025high} studies such tests in high dimensions, but only under the null and FvML alternatives.


The rest of the paper is organized as follows. The proposed pseudometric, the proposed test and a general consistency result are presented in Section \ref{main-results}. The general lower bound and asymptotic results in high-dimensional parametric models are provided in Section \ref{sec:asymptotic d}. Local limiting distribution and minimax rates for high-dimensional parametric models are presented in Section \ref{sec lower-bound}. The results for non-local alternatives are given in Section \ref{sec: non-local}. Comparison with projection-based tests and simulation studies can be found in Section \ref{when to use} and \ref{sec:simulation}, respectively.
Section \ref{discussion} contains the conclusions and some remarks.

\section{Measuring uniformity deviation and testing procedure} \label{main-results}

\subsection{Notation and preliminaries} \label{subsec-notation}
Throughout the paper, we consider the hypersphere $\mb{S}^{p-1}= \la x \in \mb{R}^{p}: \|x\|_2 =1 \ra$, where $\|.\|_2$ is the Euclidean distance. Throughout, unless stated otherwise, we implicitly assume that the dimension $p=p_n$ diverges to infinity (as $n \to \infty$). The observed data points are denoted by $\bm{X}_1, \bm{X}_2,...,\bm{X}_n$ with $\bm{X}_i \in \mb{S}^{p-1}$ for all $i=1,2,...,n$.  We assume that $\bm{X}_i$'s are drawn independently from an unknown distribution $\mu$ supported on the hypersphere $\mb{S}^{p-1}$.  The uniform distribution on the hypersphere $\mb{S}^{p-1}$ is denoted by $\mbox{Unif}(\mb{S}^{p-1})$. 

For a pair of data points $(\bm{X}_i, \bm{X}_j)$, we denote by $\bm{X}^{\top}_i \bm{X}_j$ the inner product formed by $\bm{X}_i$ and $\bm{X}_j$. Under $H_0$, the distribution of $\bm{X}_1 ^\top \bm{X}_2$ is known to have the law (see, for example, Lemma 11 and 12 in \cite{Jiang13})
\begin{align}
    \mb{P} \lb a \leq \bm{X}_1 ^\top \bm{X}_2 \leq b  \rb &=  \frac{1}{\sqrt{\pi}} \frac{\Gamma\left(\frac{p}{2}\right)}{\Gamma\left(\frac{p-1}{2}\right)} \cdot \int_{a}^{b} \left(1-\rho^2\right)^{\frac{p-3}{2}} d\rho \label{density-inner}.
\end{align}
In the formula (\ref{density-inner}) above, $a$ and $b$ are taken to be in $(-1,1)$. The CDF and density of a standard normal distribution $N(0,1)$ will be denoted by $\Phi(t)$ and $\phi(t)$, respectively. 
Throughout the paper, we will use $\mu_0$ or $\mb{P}_0$ to denote the uniform distribution. The notation $\mb{P}_{0n}$ is used to indicate the $n$-fold product measures of the uniform distribution. 

\subsection{A characterization of the uniform distribution}
Let us start with an important observation regarding the random inner product of two i.i.d.\ points on the hypersphere. 

\begin{prop} \label{identifiable}
  Let \(\nu\) be any Borel probability measure on \(\mathbb{S}^{p-1}\), and let \(\mu_0\) denote the uniform distribution on \(\mathbb{S}^{p-1}\). Suppose \(\bm{X}_1, \bm{Y}_1\) are drawn independently from \(\nu\), and \(\bm{X}, \bm{Y}\) are drawn independently from \(\mu_0\). If  
\begin{align} \label{characterization}
    \bm{X}_1^\top \bm{Y}_1 \stackrel{d}{=} \bm{X}^\top \bm{Y},
\end{align}
then \(\nu \equiv \mu_0\).
\end{prop}


Proposition \ref{identifiable} establishes the identifiability of the uniform distribution in terms of the inner product, asserting that the distribution of the inner product uniquely characterizes \(\mathrm{Unif}(\mathbb{S}^{p-1})\). To the best of our knowledge, this characterization is new and may be of independent interest. Notably, we do not impose any regularity assumptions on the measure \(\nu\) in the statement of Proposition \ref{identifiable}, and the result holds even if \(\nu\) is singular. For instance, the characterization applies to probability measures supported on lower-dimensional sets. Proposition \ref{identifiable} allows us to construct tests that work against any type of alternative, distinguishing it from other omnibus tests in the literature, which typically require alternatives to have \(L^2\)-integrable densities with respect to \(\mbox{Unif}(\mathbb{S}^{p-1})\).

A related but distinct characterization of the uniform distribution in terms of projections was proposed in \cite{cuesta2009projection}; see Section \ref{comparison} for further details and comparisons with the class of projection-based tests, which are built upon this characterization. Some recent extensions of this characterization to other types of distributions can be found in \cite{fraiman2023cramer,fraiman2024application,fraiman2023quantitative}. Specifically, the characterization in \cite{cuesta2009projection} is as follows. Let \(\bm{X}_1, \bm{X}_2\) be random variables on \(\mathbb{S}^{p-1}\) for some \(p \geq 1\), and let \(U \sim \mbox{Unif}(\mathbb{S}^{p-1})\). Then, under some regularity conditions,
\begin{align} \label{identifiable-2}
    \bm{X}_1^\top \bm{U} \stackrel{d}{=} \bm{X}_2^\top \bm{U} \Leftrightarrow \bm{X}_1  \stackrel{d}{=} \bm{X}_2. 
\end{align}
Broadly speaking, the characterization \eqref{identifiable-2} relies on projections onto independent, uniformly distributed directions. To construct testing procedures using \eqref{identifiable-2}, one often needs to sample the direction $\bm{U}$ repeatedly. In contrast, the left-hand side of \eqref{characterization} does not involve any uniformly distributed direction and is completely data-driven. Consequently, our method requires neither integrating over all directions $\bm{U}$, as was done in \cite{escanciano2006consistent,garcia2023projection,fernandez2023new}, nor sampling $\bm{U}$ repeatedly, as was done in \cite{cuesta2009projection}.

The key challenge in proving Proposition \ref{identifiable} is that the arguments used to establish \eqref{identifiable-2} are no longer applicable. Specifically, the proof of the characterization \eqref{identifiable-2} in \cite{cuesta2007sharp,cuesta2009projection} relies on a sharp version of the Cramér--Wold device, which cannot be applied to \eqref{characterization} because \(\bm{Y}_1\) and \(\bm{Y}\) follow different laws. The proof of Proposition \ref{identifiable}, instead, relies on a subtle application of the Lebesgue differentiation theorem.
 \subsection{Testing procedure}
Proposition \ref{identifiable} suggests the following test statistic:

\begin{align} \label{def-T}
    T_n :&= \sup_{t \in [-1,1]} \Big|  \frac{2}{n(n-1)} \sum_{1 \leq i<j \leq n} \mathbf{1}_{\la \bm{X}^{\top}_i \bm{X}_j  \leq t \ra} - m(t) \Big|,
\end{align}
where
\begin{align} \label{cdf-uni}
    m(t):= \mb{P} \lb \bm{X}^{\top}_1 \bm{X}_2 \leq t \Big| H_0  \rb= \mb{P} \lb 1-2U \leq t \rb,
\end{align}
and $U \sim \mbox{Beta}\lb \frac{p-1}{2}, \frac{p-1}{2}\rb$.

We choose the $L^\infty$ distance rather than an $L^2$ distance in \eqref{def-T} primarily for convenience in establishing lower bounds. In particular, the analysis of the distance $d$ in \eqref{metric} becomes substantially more tractable under the $L^\infty$ formulation than under an $L^2$ formulation, since the available empirical process tools are better suited to the $L^\infty$ distance.

Recall that the Brownian bridge $\la B_t; 0 \leq t \leq 1 \ra$ has the same distribution as $\la W_t - tW_1; 0 \leq t \leq 1 \ra$, where $\la W_t \ra$ is a standard Brownian motion. Our next result establishes the asymptotic distribution of $T_n$ under $H_0$, namely 

 \begin{thm} \label{null-dist}
 Let $\la B_t \ra_{0 \leq t \leq 1}$ be a standard Brownian bridge on $[0,1]$. If $\min \la n, p \ra \to \infty$, then
    $$\sqrt{\frac{n(n-1)}{2}} T_n \xrightarrow{d} \max_{t \in [0,1]} |B_t|$$
    where $T_n$ is defined in (\ref{def-T}).
 \end{thm}

The asymptotic distribution in Theorem \ref{null-dist} is the same as that of the classical nonparametric Kolmogorov--Smirnov test. The exact expression for the maximum of a Brownian bridge is known to be 
\begin{align} \label{bridge-asymptotic}
    \mb{P} \lb \max_{t \in [0,1]} |B_t| > x \rb = 2\sum_{k=1}^{\infty} (-1)^{k+1} \exp \lb -2k^2x^2 \rb;
\end{align}
see, for example, \cite{brownian-bridge}. 

At first glance, one can see that $T_n$ is the supremum of a degenerate U-process. Moreover, the underlying distribution of the U-process is also allowed to change with the dimension. To the best of our knowledge, this situation has not been explored in the literature, although limit theory for non-degenerate U-processes is well studied. To establish the null distribution of $T_n$, we rely on a special property of the uniform distribution on the sphere: under $H_0$, the normalized inner products $\sqrt{p} \bm{X}^{\top}_i \bm{X}_j$ are pairwise independent and asymptotically normal. Moreover, as discussed above, the dependence between them becomes weaker as the dimension increases. Thus, one should expect the asymptotic distribution of $T_n$ to be the same as that of the classical Kolmogorov--Smirnov test in the i.i.d.\ setting. Remarkably, Theorem \ref{null-dist} presents a stark departure from classical results concerning degenerate U-processes with a {\it fixed distribution}. Typically, the asymptotic distributions of such statistics lack closed-form expressions; see \eqref{fix-p} below. However, the divergence of $p$ gives a convenient asymptotic distribution, as demonstrated in Theorem \ref{null-dist}.

Thanks to Theorem \ref{null-dist} and the expression (\ref{bridge-asymptotic}), one can easily calculate the critical value $c_{\alpha}$ for the $\alpha$-level test. The test rejects the null hypothesis if $T_n > \sqrt{2}c_{\alpha} / \sqrt{n(n-1)}$, where $c_{\alpha}$ is chosen such that 
$$\mb{P} \lb \max_{t \in [0,1]} |B_t| \leq c_{\alpha} \rb = 1- \alpha.$$
The quantile $c_{\alpha}$ of the Kolmogorov--Smirnov distribution is well understood in the literature and can be calculated precisely via (\ref{bridge-asymptotic}); for example, $c_{0.95} = 1.36$. This is also the critical value we choose in the simulation. From now on, we will use $\phi_n$ to indicate the $\alpha$-level test based on $T_n$, which is

\begin{align}
    q_n(\alpha)  := \frac{\sqrt{2}c_{\alpha}}{\sqrt{n(n-1)}} ~~~~~\text{and}~~~~~
    \phi_n :=  \mathbf{1}_{\la T_n \geq q_n(\alpha) \ra}, \label{test-phi}
\end{align}
where $T_n$ is defined in (\ref{def-T}).

From Theorem \ref{null-dist}, one can see that the test $\phi_n$ is doubly robust: there is no restriction on how $p$ diverges to infinity. Among the three known high-dimensional tests, only the Rayleigh test $R_n$ in \eqref{Rayleigh} and the Bingham test $B_n$ in \eqref{Bingham} satisfy this property. The packing test $P_n$ requires a mild regularity condition $p/(\log n)^2 \to \infty$ and thus is not doubly robust. It is worth noting that the test $\phi_n$ is also valid in the fixed-$p$ scenario, albeit without the asymptotic distribution provided in Theorem \ref{null-dist}. Indeed, when $p$ is fixed, it is known that
\begin{align} \label{fix-p}
	      \sqrt{\frac{n(n-1)}{2}} T_n \xrightarrow{d} \max_{t \in [0,1]} |Q_t|
\end{align}
where $T_n$ is defined in (\ref{def-T}) and $\la Q_t \ra_{0 \leq t \leq 1}$ is a stochastic process whose marginal distributions are equal to linear combinations of independent chi-squared distributions; see Theorem 7 in \cite{nolan1988functional}. The asymptotic distribution in \eqref{fix-p} does not have a tractable expression, and Monte Carlo simulation is needed to approximate the critical value of the test.

\subsection{Model-free consistency}
For any two probability measures $\mu$ and $\nu$ supported on $\mb{S}^{p-1}$, define the pseudometric
\begin{align} \label{metric}
	d(\mu,\nu) := \sup_{t \in [-1,1]} \Big| \mb{P}_{\mu} \lb \bm{X}_1^\top \bm{Y}_1 \leq t \rb - \mb{P}_{\nu} \lb \bm{X}^\top \bm{Y} \leq t \rb \Big|
\end{align}
where $\bm{X}_1$, $\bm{Y}_1$ are drawn independently from $\mu$, and $\bm{X}$, $\bm{Y}$ are drawn independently from $\nu$.

The distance $d$ in \eqref{metric} is only a pseudometric, in the sense that $d(\mu,\nu)=0$ does not imply $\mu \equiv \nu$. However, in light of Proposition \ref{identifiable}, we can see that $d(\mu,\mbox{Unif}\lb \mb{S}^{p-1} \rb)=0$ yields $\mu \equiv \mbox{Unif}\lb \mb{S}^{p-1} \rb$. Thus, the pseudometric $d$ can be used as a quantitative measure of deviation from the null. 

The reader may wonder why the pseudometric \eqref{metric} is more convenient than existing projection-based characterizations. For example, one may naturally consider
\begin{align*}
    d_1(\mu,\nu) &=  \sup_{\bm{u} \in \mb{S}^{p-1}} \sup_{t \in [-1,1]} \left| \mb{P}_{\mu} \lb \bm{u}^\top \bm{X} \leq t \rb - \mb{P}_{\nu} \lb \bm{u}^\top \bm{X} \leq t \rb \right|; \\
    d_2 \lb \mu, \nu \rb &= \int_{-1}^1 \int_{\mb{S}^{p-1}} \left[  \mb{P}_{\mu} \lb \bm{u}^\top \bm{X} \leq t \rb - \mb{P}_{\nu} \lb \bm{u}^\top \bm{X} \leq t \rb \right]^2 d\bm{u}\, dt.
\end{align*}
In fact, certain variants of $d_2$ have been used in \cite{garcia2023projection} to construct tests for \eqref{uniform-test} in the fixed-dimensional setting. The difficulty with $d_1$ is that its empirical counterpart involves taking the supremum over the class of functions
\[
\la \mathbf{1}_{\la \bm{u}^\top \bm{x} \leq t \ra} \ra_{\bm{u} \in \mb{S}^{p-1},\, t \in [-1,1]},
\]
whose VC dimension is of order $\Theta(p)$.

More concretely, let
\[
F_{\mu,\bm u}(t):=\mb P_\mu(\bm u^\top \bm X\le t),
\qquad
\hat F_{n,\bm u}(t):=\frac1n\sum_{i=1}^n \mathbf 1_{\{\bm u^\top \bm X_i\le t\}}.
\]
The natural empirical counterpart of $d_1$ is
\[
\hat d_1(\mu,\nu)
:=
\sup_{\bm u\in \mb S^{p-1}}\sup_{t\in[-1,1]}
\bigl|\hat F_{n,\bm u}(t)-F_{\nu,\bm u}(t)\bigr|,
\]
and a direct plug-in bound gives
\[
\bigl|\hat d_1(\mu,\nu)-d_1(\mu,\nu)\bigr|
\le
\sup_{\bm u\in \mb S^{p-1}}\sup_{t\in[-1,1]}
\bigl|\hat F_{n,\bm u}(t)-F_{\mu,\bm u}(t)\bigr|.
\]
Now the class
\[
\mathcal H_p
=
\Bigl\{
\mathbf 1_{\{\bm u^\top \bm x\le t\}}
:\bm u\in \mb S^{p-1},\ t\in[-1,1]
\Bigr\}
\]
has VC dimension of order $\Theta(p)$, so empirical process conccentration inequalities imply that
\[
\sup_{\bm u,t}
\bigl|\hat F_{n,\bm u}(t)-F_{\mu,\bm u}(t)\bigr|
=
O_{\mb P}\!\left(\sqrt{\frac{p}{n}}\right).
\]
Consequently, the concentration of $\hat d_1$ around $d_1$ is controlled by an empirical process whose complexity grows linearly with $p$. This is undesirable, since $\sqrt{p/n}$ diverges when $p \gg n$. Likewise, $d_2$ is defined through a double integration, which makes it difficult to analyze asymptotically under general high-dimensional alternatives, as well as to quantify the concentration between $d_2$ and its empirical counterpart. Moreover, many common high-dimensional spherical distributions exhibit concentration around an unknown location parameter, so the deviation from uniformity is often localized rather than spread out over the sphere. In such cases, a max-type distance such as $d_1$ is intuitively more sensitive than an integrated distance such as $d_2$.

In contrast, the empirical counterpart of \eqref{metric}, namely $T_n$ in \eqref{def-T}, involves a supremum over a function class of VC dimension $2$, which makes both its asymptotic behavior and its concentration properties much easier to analyze. In fact, in Section \ref{sec:minimax} we show that simple Edgeworth-type expansions can be used to study the behavior of $d$ under various high-dimensional parametric alternatives, and in Theorem \ref{consistency} we show that the consistency of the resulting test \eqref{test-phi} is determined entirely by the behavior of $d$.


Based on this pseudometric \eqref{metric}, we define a consistency criterion, namely the {\it $1/n$-separation condition}. The precise definition is as follows.

{
\begin{cond}[Separation condition]  Given a sequence $\la  (n,p_n) ; {n \geq 1} \ra$, let $\mu_n$ be a sequence of probability measures on $\mb{S}^{p_n-1}$. We say that the sequence $\la \mu_n \ra_{n \geq 1}$ satisfies the {\it $1/n$-separation condition} if
    \begin{align} \label{rootn}
        n \cdot d \lb \mu_n, \mbox{Unif}\lb \mb{S}^{p_n-1} \rb \rb \to \infty
    \end{align}
    where $d$ is the pseudometric defined in \eqref{metric}.
\end{cond}}

The separation condition (\ref{rootn}) measures departure from the null hypothesis in terms of the pseudometric $d$, which is the Kolmogorov distance between the random inner products drawn under $H_0$ and $H_1$. Interestingly, the rate in \eqref{rootn} is of order $n^{-1}$, which differs from the usual $n^{-1/2}$ rate. This is due to the degenerate nature of $H_0$ and the form of $d$.  
Interestingly, condition (\ref{rootn}) is nonparametric in nature and requires neither parametric assumptions nor regularity conditions: the sequence of alternatives $\mu_n$ may or may not have densities with respect to the uniform measure $\mu_0$, and it is not restricted to any parametric class of distributions containing $\mbox{Unif}(\mb{S}^{p-1})$.

We do not require $p_n$ to converge to infinity in \eqref{rootn}, and the fixed-$p$ setting is also covered by \eqref{rootn}. The assumption that $p$ is diverging is only required to control the size of the test via the asymptotic result in Theorem \ref{null-dist}. In the fixed-$p$ scenario, one can use Monte Carlo simulation to approximate the critical value of the test, as stated in \eqref{fix-p}. If we keep $p$ fixed and consider a fixed alternative, then by Proposition \ref{identifiable}, \eqref{rootn} always holds. Therefore, in the fixed-dimensional case, \eqref{rootn} is the same as the universal consistency property of Sobolev tests. Furthermore, condition \eqref{rootn} remains valid in high-dimensional settings, making it a natural analogue of universal consistency that operates in both fixed- and high-dimensional cases. Next, given the separation condition (\ref{rootn}), it can be shown that the test $\phi_n$ is consistent.

{
\begin{thm} \label{consistency}
    Let $\mu_n$ be a sequence of probability measures on $\mb{S}^{p_n-1}$ that satisfies the separation condition \eqref{rootn}. Then
    \[
    \lim_{n\to\infty}\mb{P}_{\mu_n}(\phi_n =1) = 1.
    \]
    Here $\phi_n=\phi_n(\bm{X}_1,\bm{X}_2,\dots,\bm{X}_n)$ is the test defined in \eqref{test-phi}.
\end{thm}
}

{ The rate $1/n$ in condition \eqref{rootn} is sharp: we prove a matching lower bound in Theorem \ref{information bound} below, which also does not impose any restriction on the relationship between $p$ and $n$. The proof of Theorem \ref{consistency} is short but somewhat delicate. For a general sequence of alternatives $\mu_n$, the alternatives may lie in a non-local regime in which $d(\mu_n,\mu_0)$ is bounded away from zero, while the first few moments of the inner products under $\mu_n$ and $\mu_0$ might be close. We avoid the need to analyze this phenomenon directly by identifying the dominant contribution through the Hoeffding decomposition.

A notable feature of Theorem \ref{consistency} is that it yields a provably consistent test in both fixed- and high-dimensional settings, without any restriction on the way the dimension diverges. By contrast, even for the classical Sobolev tests \cite{Gine,Jupp}, universal consistency in fixed dimensions requires an infinite sequence of weights. This leads to statistics that are analogous to \eqref{Rayleigh} or \eqref{Bingham}, but involve infinite-order polynomial expansions rather than the linear or quadratic kernels appearing in \eqref{Rayleigh} and \eqref{Bingham}.

Interestingly, when one restricts attention to a high-dimensional parametric class, the threshold at which the distance $d$ scales like $1/n$ often coincides with the minimax rate in that model. We do not have a result of this type for arbitrary sequences of alternatives. However, we will investigate three examples of this phenomenon in detail below.}


\section{Asymptotic results for the distance $d$} \label{sec:asymptotic d}

Recall the distance $d$ from \eqref{metric}. In this section, we highlight an interesting phenomenon that it captures precisely the minimax rates of various parametric models and develop a general framework for understanding the asymptotic properties of $d$ when restricted to a parametric model containing the null.

\subsection{High-dimensional parametric models}

Let us start by introducing the three parametric models that will be considered in the rest of the paper. Let $f: \mb{R} \to \mb{R}^+$ be a smooth function. Define the family of densities 
\begin{align} \label{semiparametric model}
\bm{x} \mapsto c_{f,\kappa} \cdot \exp \Big[ \kappa f \lb  \bm{x}^\top \bm{\mu}    \rb
\Big] d\mu_0 \lb \bm{x} \rb
\end{align}
Here $\kappa>0$ is the concentration parameter and $\bm{\mu} \in \mb{S}^{p-1}$ is the location parameter. Most of the common distributions in directional statistics belong to this class. Two common choices are
\begin{itemize}
    \item {\it  FvML distributions.} This corresponds to the case  $f(x)=x$ \cite{Cutting-P-V}.
       \item {\it  Watson distributions.} This corresponds to the case $f(x)=x^2$ \cite{Cutting-P-V2}.
\end{itemize}
These two classes of distributions are among the most common in directional analysis; see \cite{M-Jupp} for more details. The minimax rates for testing in these models have been studied in \cite{Cutting-P-V,Cutting-P-V2}. It is known that under the FvML model, the minimax rate is $\kappa \sim p^{3/4}/\sqrt{n}$ and is achieved by the Rayleigh test. It is also shown there that the Rayleigh test $R_n$ is optimal in the Le Cam sense under this model. Regarding the class of Watson distributions, the minimax rate is $\kappa \sim p^{3/2}/\sqrt{n}$ and is achieved by the Bingham test $B_n$ in the regime $p=o(n)$.

In addition to the two models above, we will also study another model involving the uniform distribution supported on the intersection of a lower-dimensional hyperplane and the hypersphere. This model is an example of a class of distributions that are singular with respect to the uniform law. Consider the set of $k$-dimensional hyperplanes in $\mathbb{R}^p$. We denote this set by $G(k,p)$, which is known to form the Grassmannian manifold; see \cite{chikuse2003statistics} for a comprehensive overview. Consider the class of distributions
\begin{align} \label{low rank uniform}
    \la \mu: \exists\  (k,\mathcal{H}) \ \text{s.t.} \ \mathcal{H}\in G(k,p)  \text{ and } 
\mu \sim \mathrm{Unif}\!\left( \mathcal{H} \cap \mathbb{S}^{p-1} \right)  \ra
\end{align}
for some $k \in \{1, \dots, p-1\}$.

Obviously, the case $k=p$ corresponds to the null hypothesis $H_0$. The uniformity testing problem in this model is essentially about detecting whether the uniform distribution is supported on a low-dimensional hyperplane. We are interested in the hard regime where $k$ is close to $p$, and we will thus assume that $k=k_n$ such that $\min \la k,p,n \ra \to \infty$ and $k/p \to 1$. The information lower bound for testing uniformity in this model is studied in Theorem \ref{low-rank information bound}, where the minimax rate is shown to be $k=p\lb 1 - \Theta(1/n) \rb$.

\subsection{An information lower bound with respect to $d$}
Define the set of test functions based on a sample of size $n$ as
\begin{align} \label{test class}
\mathcal{T}_n:= \la \phi=\phi \lb \bm{X}_1, \bm{X}_2,\dots,\bm{X}_n \rb: \lb \mb{S}^{p-1} \rb^n \to \la 0,1 \ra \ra.
\end{align}
Our first result is a general lower bound for testing uniformity under the distance $d$. By using Le Cam's mixture argument, we can show that
\begin{thm} \label{information bound}
Suppose $\min \la p,n \ra \to \infty$. For $\ve$ small enough, we have 
\begin{align} \label{minimax}
    \liminf_{n \to \infty} \la \inf_{\phi \in \mathcal{T}_{n}} \la \mb{P}_{\mu_0} \lb \phi=1 \rb + \sup_{d(\nu,\mu_0) \geq \frac{\ve}{n}}  \mb{P}_\nu \lb \phi =0 \rb \ra \ra \geq 1/4.
\end{align}
\end{thm}
In the fixed-dimensional settings, Theorem~\ref{information bound} is straightforward to prove since one can directly apply Le Cam’s two-point argument to a perturbation of size $\Theta(1/\sqrt{n})$ of the uniform distribution. The non-trivial aspect of Theorem~\ref{information bound} lies in establishing the result in high-dimensional settings, and doing so {\it without} imposing any growth condition on $p$ and $n$. Theorem \ref{information bound} shows that the minimax lower bound for testing uniformity with respect to $d$ is of order $1/n$, and over the whole class of probability distributions over the hypersphere with no extra assumption. Theorem \ref{consistency} above shows that the proposed test in \eqref{test-phi} achieves this rate.

The worst-case construction in the proof of Theorem~\ref{information bound} is based on the Fisher–von Mises–Langevin (FvML) distributions. This choice is motivated by simulation results showing that the proposed test in \eqref{test-phi} exhibits power very close to that of the Rayleigh test, which is the optimal invariant test (in the Le Cam's sense) within this model \cite{Cutting-P-V}.

\subsection{Minimaxity in high-dimensional parametric models} \label{sec:minimax}
 We now show that the distance $d$ scales like $\Theta(1/n)$ at the minimax rates under the three models described above. 

\begin{prop} \label{Kolmogorov expansion fvml}
    Consider the FvML model in \eqref{semiparametric model} with $f(x)\equiv x$. Let $\kappa =\tau p_n^{3/4}/\sqrt{n}$ with $\tau \in (0,\infty)$. Suppose $\bm{X}, \bm{Y}$ are two i.i.d. random points on $\mb{S}^{p-1}$. Then, for any $u \in \mb{R}$, we have 
    \begin{align*}
    n \left[ \mb{P}_{\mu_n} \lb \sqrt{p}\bm{X}^\top \bm{Y} \leq u  \rb - \mb{P}_{\mu_0} \lb \sqrt{p}\bm{X}^\top \bm{Y} \leq u  \rb \right] \to \frac{-\tau^2}{\sqrt{2\pi}} \exp \lb -u^2/2 \rb
    \end{align*}
   as $\min \la n,p \ra \to \infty$, where $\mu_n$ is an FvML distribution on $\mb{S}^{p_n-1}$ with concentration parameter $\kappa$. Moreover, the convergence above is uniform in $u$, which also means 
   $$n \cdot d \lb \mu_n, \mu_0 \rb \to \tau^2/\sqrt{2 \pi}$$ 
   where $d$ is the distance in \eqref{metric}.
\end{prop}
Proposition \ref{Kolmogorov expansion fvml} and Theorem \ref{consistency} show that the proposed test \eqref{test-phi} achieves the minimax rate \(p^{3/4}/\sqrt{n}\) in this model. We also note that \cite{ebner2025high} shows that {\it any} Sobolev test based on a Gegenbauer-polynomial kernel of degree \(k>1\) does not achieve this minimax rate. In particular, any linear combination of the tests \eqref{Rayleigh} and \eqref{Bingham} fails to attain this rate. Although \cite{ebner2025high} proposes using a suitable infinite-order kernel to recover the minimax rate, the resulting test statistic does not admit a simple closed-form expression. By contrast, the test \eqref{test-phi} has a simple closed-form expression and still achieves the minimax rate.

Regarding the Watson model, we have
\begin{prop} \label{watson d asymptotic}
    Suppose $\kappa < p/2$, $p/n^{2/3} \to \infty$. Denote by $\mu_n$ the Watson distribution with concentration parameter $\kappa$. Assume that
    \begin{align} \label{watson scaling limit}
        \frac{n \kappa^2}{p(p/2- \kappa)^2} \to \tau \in (0,\infty).
    \end{align}
    Suppose $\bm{X}, \bm{Y}$ are two i.i.d. random points on $\mb{S}^{p-1}$. Then, for any $u \in \mb{R}$, we have 
    \begin{align*}
    n \left[ \mb{P}_{\mu_n} \lb \sqrt{p}\bm{X}^\top \bm{Y} \leq u  \rb - \mb{P}_{\mu_0} \lb \sqrt{p}\bm{X}^\top \bm{Y} \leq u  \rb \right] \to \frac{-\tau}{2\sqrt{2\pi}} \cdot u\exp \lb -u^2/2 \rb
    \end{align*}
  In particular, 
  $$n\cdot d \lb \mu_n, \mu_0 \rb \to \frac{\tau}{2\sqrt{2e \pi}}$$
  where $d$ is the distance in \eqref{metric}.
\end{prop}
We did not specify the location parameter of the FvML and Watson distributions in the statements of Propositions \ref{Kolmogorov expansion fvml} and \ref{watson d asymptotic}, since the distribution of the inner product is independent of the choice of location. The condition $p/n^{2/3} \to \infty$ in Proposition \ref{watson d asymptotic} can be further relaxed by expanding the Taylor expansion to higher orders. In the proof of Proposition \ref{watson d asymptotic}, we do this only up to second order for the sake of simplicity. 

It may not be immediate that \eqref{watson scaling limit} is the same as the minimax rate $\kappa \sim p^{3/2}/\sqrt{n}$ in \cite{Cutting-P-V2}. To relate \eqref{watson scaling limit} to the information lower bound in \cite{Cutting-P-V2}, note that
\begin{align} \label{relation}
        \frac{n \kappa^2}{p(p/2- \kappa)^2} = \tau_n \iff \kappa = \frac{p^{3/2} \sqrt{\tau_n}}{2 \lb \sqrt{n} + \sqrt{\tau_n p} \rb}.
\end{align}
The last expression is exactly the same as the minimax lower bound in \cite{Cutting-P-V2} for the regime $p=o(n)$, namely $\Theta \lb p^{3/2}/n^{1/2} \rb$. Therefore, in the regime $n^{2/3} \ll p \ll n$, the rate in \eqref{watson scaling limit} coincides with that in \cite{Cutting-P-V2}. However, it is of strictly smaller order when $p/n$ diverges. In such case, the minimax rate scales like $p/2-\Theta \lb  \sqrt{pn}\rb$.  A matching lower bound is given in Theorem \ref{watson minimax} below. 

Finally, regarding the low-dimensional model, we can prove that
\begin{prop} \label{Kolmogorov low-rank}
Suppose $k \leq p$, $p/n \to \infty$, and $\lb 1 - k/p \rb n \to \tau \in (0,\infty)$. Suppose $\bm{X}, \bm{Y}$ are two i.i.d. random points on $\mb{S}^{p-1}$. Denote by $\mu_n$ the low-dimensional uniform distribution in \eqref{low rank uniform} for some $k$. Then, we have 
\[
n \cdot \left[ \mb{P}_{\mu_n} \lb \sqrt{p} \bm{X}^\top \bm{Y} \leq u \rb - \mb{P}_{\mu_0} \lb \sqrt{p} \bm{X}^\top \bm{Y} \leq u \rb  \right] \to \frac{-\tau}{2} \cdot \frac{u \cdot \exp \lb -u^2/2 \rb}{\sqrt{2 \pi}}
\]
as $n \to \infty$, uniformly in $u \in \mb{R}$. In particular, 
$$n \cdot d \lb \mu_n, \mu_0 \rb \to \frac{\tau}{2\sqrt{2e \pi}}$$
where $d$ is the distance in \eqref{metric}.
\end{prop}
Propositions \ref{Kolmogorov expansion fvml}, \ref{watson d asymptotic}, and \ref{Kolmogorov low-rank} reveal an interesting phenomenon: the local alternatives at the minimax rate with respect to $d$ (which is indeed of order $1/n$ by Theorem \ref{information bound} below) coincide exactly with the local alternatives in the usual parametric sense. More precisely, for each of the three models above, there exists a local sequence of alternatives $\mu_n$ such that
\[
n\, d(\mu_n,\mu_0)\to c \in (0,\infty)
\]
if and only if the model parameter lies at the usual minimax scale. Concretely, this critical regime is given by
\[
\kappa_n \asymp \frac{p_n^{3/4}}{\sqrt{n}}
\qquad \text{for the FvML model},
\]
\[
\frac{n\kappa_n^2}{p_n(p_n/2-\kappa_n)^2} \asymp 1
\qquad \text{for the Watson model},
\]
and
\[
n\left(1-\frac{k_n}{p}\right)\asymp 1
\qquad \text{for the low-dimensional model}.
\]
Equivalently, in all three cases,
\[
d(\mu_n,\mu_0)\asymp \frac{1}{n}
\]
at the minimax rates.

We conjecture that this phenomenon holds for a much broader range of models, although we do not currently have a proof in such generality. Nevertheless, none of the existing high-dimensional tests achieves the minimax rates simultaneously across all of these models. For example, the Rayleigh test $R_n$ defined in \eqref{Rayleigh} is blind to the Watson distributions at the minimax threshold \cite{Cutting-P-V2}. The Bingham test $B_n$ defined in \eqref{Bingham} achieves a sub-optimal rate under the class of FvML distributions, and the Packing test $P_n$ defined in \eqref{packing} is blind to the FvML distributions at the minimax threshold \cite{jiang2025asymptotic}.

The analysis of the distance $d$ in the three models above follows a common scheme, despite the differences between the underlying alternatives. In each case, one starts from the random inner product $\sqrt{p}\,\bm{X}^{\top}\bm{Y}$, where $\bm{X},\bm{Y}$ are i.i.d.\ from the alternative distribution, and rewrites it as a perturbation of the null inner product. More precisely, we decompose
\[
\sqrt{p}\,\bm{X}^{\top}\bm{Y}
=
A+(1+B)\Xi,
\]
where $\Xi$ has the same distribution, or the same leading-order asymptotic behavior, as the normalized inner product under the null, while $A$ and $B$ are (small) correction terms that encode the effect of the model parameter and are independent from $\Xi$.  We then expand
\[
u \mapsto \mb{P}\lb A+(1+B)\Xi \leq u \rb
\]
around the null distribution by a Taylor or Edgeworth-type expansion, often after conditioning on $A, B$ (which are independent from $\Xi$). This reduces the study of $d$ to a comparison of a few low-order moments of $A$ and $B$ under the null and the alternative. The leading term is determined by the first order at which this moment matching fails, while the remaining terms are controlled uniformly through derivative bounds in the expansion.

Regarding the proof of Proposition \ref{Kolmogorov expansion fvml}, we use a different and somewhat tricky argument in order to preserve the clean condition $\min \la n,p \ra \to \infty$: the FvML model admits a LAN expansion \cite{Cutting-P-V}, which can be combined with Le Cam's third lemma to establish the result. To the best of our knowledge, LAN expansions of this type are not available for any model other than the FvML model.

\section{Lower bounds and power under local alternatives} \label{sec lower-bound}

\subsection{Local limiting distribution under the FvML alternatives} \label{sec fvml}


In this subsection, we will investigate the local power and consistency of the proposed test under the class of FvML distributions. It is known that within the class of FvML distributions, the threshold $\kappa \sim p^{3/4}/\sqrt{n}$ is the minimax rate: when $\kappa$ is below this threshold, no test can be consistent. Moreover, when $\kappa$ is above this threshold, the Rayleigh test is consistent and is also optimal in the sense of Le Cam. 

Let $\Phi^{-1}(t):[0,1] \to \mb{R}$ be the quantile function of the standard normal distribution, and let $\phi$ be the standard Gaussian density, $\phi(x)=(1/\sqrt{2\pi}) \exp \lb -x^2/2 \rb$. Our main result regarding the FvML alternatives is
\begin{prop} \label{local power FvML}
 Let $\kappa= \tau_n p^{3/4}/\sqrt{n}$.   Then, if $\tau_n \to \tau \in (0,\infty)$, then 
     \[
     \sqrt{\frac{n(n-1)}{2}}T_n \stackrel{d}{\to} \sup_{t \in [0,1]} \left| B_t - \frac{\tau^2}{\sqrt{2}} \phi \lb \Phi^{-1}(t) \rb  \right|.
     \]
     under the sequence of FvML alternatives with concentration parameter $\kappa_n$, where $\la B_t; 0 \leq t \leq 1 \ra$ is the Brownian bridge. 
 \end{prop}
It follows directly from Proposition \ref{local power FvML} that the asymptotic power of $T_n$ under the class of FvML distributions is given by
     \[
     \mb{P} \lb \sup_{t \in [0,1] } \left| B_t - \frac{\tau^2}{\sqrt{2}} \phi \lb \Phi^{-1}(t) \rb  \right| \geq c_{\alpha}   \rb.
     \]
From the display above, we can see that $T_n$ is consistent at the contiguity rate $p^{3/4}/\sqrt{n}$. By Proposition \ref{Kolmogorov expansion fvml}, we get
\[
n \times  d \lb \mbox{FvML} \lb \tau_n p^{3/4}/\sqrt{n}  \rb, \mu_0 \rb \to \frac{\tau^2}{\sqrt{2 \pi}}
\]
whenever $\tau_n \to \tau$. The display above indicates that the local alternatives at the minimax threshold for $d$ are the same as those of the parametric FvML model. In other words, the distance $d$ captures precisely the minimax rate of testing uniformity in the FvML model.

The asymptotic power above does not have a closed-form expression, but we observe in simulation that its power is slightly lower than that of the Rayleigh test, which is expected due to the LAN expansion in \cite{Cutting-P-V}. Note that in this regime, the Packing test $P_n$ and the Bingham test $B_n$ are both blind. We further know from \cite{Cutting-P-V2} that the detection threshold for the Bingham test in this model is $p^{3/4}/n^{1/4}$, which is strictly sub-optimal.

\subsection{Local limiting distribution under the low-dimensional uniform distributions} \label{sec low-rank}
 We are interested in testing uniformity against the class of low-dimensional uniform distributions:
\begin{align} \label{low-rank}
H_1:\quad \exists\  (k,\mathcal{H}) \ \text{s.t.} \ \mathcal{H}\in G(k,p)  \text{ and } 
\mu \sim \mathrm{Unif}\!\left( \mathcal{H} \cap \mathbb{S}^{p-1} \right)
\end{align}
for some $k \in \{1, \dots, p-1\}$.

\begin{prop} \label{power low-rank}
    Suppose $p/n \to \infty$ and $\lb 1 - k/p \rb n \to \tau \in (0,\infty)$. Let $\la B_t \ra$ be the Brownian bridge. Then, 
    \[
    \sqrt{\frac{n(n-1)}{2}}T_n \stackrel{d}{\to} \sup_{t \in[0,1]} \left| B_t - \frac{\tau}{2\sqrt{2}} \cdot \Phi^{-1}(t) \phi \lb \Phi^{-1}(t)  \rb  \right|
    \]
    where $T_n$ is defined as in \eqref{def-T}.
\end{prop}
The proof of Proposition \ref{power low-rank} follows directly from Proposition \ref{Kolmogorov low-rank} and Theorem \ref{null-dist}. Proposition \ref{power low-rank} shows that the proposed test has non-trivial power at the threshold $k=(1 - \tau/n)p$, with asymptotic power given by
\[
\mb{P} \lb \sup_{t \in[0,1]} \left| B_t - \frac{\tau}{2\sqrt{2}} \cdot \Phi^{-1}(t) \phi \lb \Phi^{-1}(t)  \rb  \right| \geq c_{\alpha}  \rb.
\]
It is natural to ask whether the rate $k=(1 - \Omega(1/n))p$ is optimal. The answer is yes, which is the claim of the theorem below.
\begin{thm} \label{low-rank information bound}
    Suppose $p/n \to \infty $ and put $\delta(k)= (1-k/p)n$. Then, for some $\ve>0$ sufficiently small, we have
    \[
      \liminf_{n \to \infty} \la \inf_{\phi \in \mathcal{T}_{n}} \la \mb{P}_{\mu_0} \lb \phi=1 \rb + \sup_{\mu_k \in H_1:\delta(k) \geq \ve} \la \mb{P}_{\mu_k} \lb \phi =0 \rb \ra \ra \ra \geq 1/4
    \]    
    where $\mathcal{T}_n$ is the class of tests based on the data as defined in \eqref{test class}, and the supremum is taken over all alternatives $\mu_k$ of the form \eqref{low-rank} such that $\delta(k) \geq \ve$.
\end{thm}
Theorem \ref{low-rank information bound} claims that in the high-dimensional regime, as long as $(1-k/p)n$ is small enough, no test based on a sample of size $n$ can be consistent. This suggests that the proposed test is rate-optimal in this model. To the best of our knowledge, this information lower bound is new and has not been studied before. The most technical part of the proof is to analyze the likelihood ratio against a random distribution over the Grassmannian $G(k,p)$.

By Proposition \ref{Kolmogorov low-rank}, we have 
\begin{align*}
    n \times d \lb \mbox{Unif} \lb H_k \cap \mb{S}^{p-1} \rb, \mu_0 \rb \to \frac{\tau}{2} \cdot \sup_{u \in \mb{R}} |u \phi(u)|
\end{align*}
where $\phi$ is the standard Gaussian density. Therefore, the local alternatives at the minimax threshold for $d$ are the same as those of the low-dimensional model. In other words, the distance $d$ captures precisely the minimax rate of testing uniformity in the low-dimensional model.

Let us now compare the local power of the four tests $R_n, B_n, P_n, T_n$. The nice feature of this low-dimensional model is that the asymptotic power and detection thresholds of each of the four tests can be computed precisely.

\begin{itemize}
    \item Recall the Rayleigh test $R_n$ from \eqref{Rayleigh}. It is easy to check that
$
\sqrt{p/k} \cdot R_n \stackrel{d}{\to} N \lb 0, 1 \rb
$
as $n \to \infty$. Thus, $R_n$ is not consistent, even when $k/p \to 0$. Its maximum power will not exceed $1/2$. However, a two-sided version of $R_n$, which rejects if $|R_n|$ is large, is consistent in the regime $k/p \to 0$.

\item For the Bingham test $B_n$ in \eqref{Bingham}, we have
$$
\frac{k}{p} \left[  B_n - \frac{p(n-1)}{2} \lb \frac{1}{k} - \frac{1}{p} \rb \right]  \stackrel{d}{\to} N \lb 0, 1 \rb
$$
as $n \to \infty$, under $H_1$. Therefore, in the regime $n(1-k/p) \to \tau \in (0,\infty)$, the asymptotic power of the Bingham test is
\[
1-\Phi \lb z_{\alpha} - \tau/2 \rb
\]
where $z_\alpha$ is the $(1-\alpha)$-quantile of the standard normal distribution. This shows that the Bingham test achieves the optimal rate suggested by Theorem \ref{low-rank information bound}, with local power given above.

\item Finally, regarding the Packing test $P_n$ in \eqref{packing}, we have 
\[
\frac{p}{k} \cdot P_n - \lb 1 - \frac{k}{p} \rb \lb 4 \log n - \log \log n \rb \to G
\]
where $G$ is a standard Gumbel law. Thus, the test $P_n$ is consistent if and only if $(\log n) (1 - k/p) \to \infty$. This detection threshold is strictly sub-optimal, but it is still better than that of the Rayleigh test.
\end{itemize}
From the above, we can see that only the proposed test \eqref{test-phi} and the Bingham test $B_n$ achieve the optimal rate. Although the power function of the proposed test does not have a closed-form expression, we find in simulation studies that the local power of the Bingham test is slightly higher than that of the proposed test.

\subsection{Local limiting distribution under the class of Watson distributions}

In this subsection, we investigate the detection rate and asymptotic distribution, along a sequence of local alternatives, of the proposed test under the class of Watson distributions, which corresponds to $f(x)=x^2$ in \eqref{semiparametric model}. The analysis of the Watson distributions is much more delicate than that of the FvML case, and the only high-dimensional result currently known in the literature is \cite{Cutting-P-V2}, which is limited to the regime $p=o(n)$. We will highlight some interesting phenomena regarding the minimax rates in this model. Let us first prove a matching lower bound corresponding to \eqref{watson scaling limit}.

\begin{thm} \label{watson minimax}
Let
\[
H^W_{1,n}:=
\left\{
\mu_n:\ \text{$\mu_n$ is a Watson distribution on $\mathbb S^{p-1}$ with concentration $\kappa<\frac{p}{2}$}
\right\}.
\]
 Define
\begin{align} \label{delta 1}
\delta_1(n)
:=
\frac{n}{p}\left(\frac{\kappa}{p/2-\kappa}\right)^2.
\end{align}
Then there exists \(\varepsilon>0\) sufficiently small such that
\[
\liminf_{n \to \infty}
\inf_{\varphi \in \mathcal{T}_{n}}
\left\{
\mathbb{P}_{\mu_0}(\varphi=1)
+
\sup_{\mu_n \in H^W_{1,n}:\,\delta_1(n) \ge \varepsilon}
\mathbb{P}_{\mu_n}(\varphi =0)
\right\}
\ge \frac14
\]
where $\mathcal{T}_n$ is the class of tests based on a sample of size $n$.
\end{thm}

In light of the relation \eqref{relation}, \eqref{delta 1} can be rewritten as
\begin{align*}
    \kappa = \frac{p^{3/2} \sqrt{\delta_1(n)}}{2 \lb \sqrt{n} + \sqrt{\delta_1(n) p} \rb}.
\end{align*}
If one keeps $\delta_1(n)=\Theta(1)$, then the minimax rate for the Watson model undergoes a phase transition between low and high dimensions: when $p=o(n)$, it scales like $\kappa=\Theta(p^{3/2}/\sqrt{n})$, as in \cite{Cutting-P-V2}, whereas when $n=o(p)$, it scales like $\kappa=p/2-\Theta(\sqrt{pn})$. Our argument is based on the truncated second-moment method similar to those in \cite{george2022robust} for high-dimensional linear model. Of course, the non-trivial part is to figure out how to truncate the likelihood ratio.

In terms of the local asymptotic distribution, we can prove that, similarly to the FvML case and the low-dimensional case, the test statistic $T_n$ in \eqref{def-T} converges in distribution to a shifted Brownian bridge.  

\begin{thm} \label{watson local distribution}
    Suppose that
    $$\min \la  \frac{p - 2\kappa}{p^{1/2}n^{1/4}}; \frac{p}{n^{2/3}}; \frac{n^{1/4}}{\log(p)} \ra \to \infty$$ and
    \[
    \frac{n}{p}\left(\frac{\kappa}{p/2-\kappa}\right)^2 \to \tau \in (0,\infty).
    \]
    Then, along the sequence of Watson distributions with concentration $\kappa$,
    \[
    \sqrt{\frac{n(n-1)}{2}}T_n \stackrel{d}{\to} \sup_{t \in[0,1]} \left| B_t - \frac{\tau}{2\sqrt{2}} \cdot \Phi^{-1}(t) \phi \lb \Phi^{-1}(t)  \rb  \right|
    \]
    where $T_n$ is defined as in \eqref{def-T} and $\la B_t\ra_{0 \leq t \leq 1}$ is the standard Brownian bridge.
\end{thm}
The proof of Theorem \ref{watson local distribution} is long and highly technical. The main difficulty in this setting is that we do not have any LAN expansion, as was used in the proof of Proposition \ref{local power FvML}. To the best of our knowledge, a high-dimensional LAN expansion (under unspecified location) is not yet available for the Watson model. We therefore require a number of different tools to obtain the limiting distribution in this local regime. Combining Theorems \ref{watson minimax} and \ref{watson local distribution}, we see that the proposed test achieves the minimax rates adaptively in both low and high-dimensional regime: it can detect local alternatives of order $\Theta \lb p^{3/2}/\sqrt{n} \rb$ when $p \ll n$ and local alternatives of order $p/2-\Theta \lb \sqrt{pn} \rb$ when $p \gg n$. Note that it is not known in \cite{Cutting-P-V2} whether the Bingham test achieves the minimax rate in the high-dimensional regime $p \gg n$.

To prove Theorem \ref{watson local distribution}, we decompose $T_n$ into four terms: one term behaves like $T_n$ under the null, two term are small random perturbations, and one term is a smooth deterministic shift. The assumption $p/n^{2/3} \to \infty$ in Theorem \ref{watson local distribution} is required to control the deterministic shift in $T_n$ along this sequence of alternatives, while the condition $(p-2\kappa) \gg p^{1/2}n^{1/4}$ is the weakest condition we can obtain in order to control the random shift via the Brascamp--Lieb inequality. The condition $\log p = o(n^{1/4})$ is needed for the empirical process maximal inequality to apply. It stems from a well-known difficulty in the empirical process literature: most existing maximal inequalities involve an upper bound of the form
$$|\log \lb \text{1/variance profile} \rb|/n^{1/2},$$
which blows up if the variance decays exponentially fast. This issue has been partially resolved for empirical processes over (weak) VC-major classes in \cite{baraud2016bounding}, but to the best of our knowledge, it has not been addressed for U-processes. Nevertheless, the condition $\log (p) \ll n^{1/4}$ is mild and valid for most high-dimensional problems.


\section{Non-local alternatives} \label{sec: non-local}
We would like a test not only to detect local alternatives around the minimax thresholds, but also to remain consistent against non-local alternatives, namely alternatives that do not arise from a smooth parametrization of a model containing the null. In this section, we consider two such examples: one arising from projections of heavy-tailed distributions, and another arising from deterministic sets of points on the hypersphere that satisfy a suitable moment-matching property with respect to polynomials on the sphere.

\subsection{Projections of heavy-tailed distributions}

\begin{class}[$\alpha$-spherical distributions] \label{heavy-tailed}
Let $\alpha \in (0,2)$. We say that a probability measure $\mu_{\alpha,p}$ on $\mathbb{S}^{p-1}$ is $\alpha$-spherical if it is the law of
\[
\frac{\mathbf{X}}{\|\mathbf{X}\|},
\]
where $\mathbf{X}=(X_1,\dots,X_p)^\top$ has i.i.d.\ symmetric coordinates that are regularly varying with index $\alpha$.
\end{class}

This class was introduced in \cite{heiny2022limiting,dornemann2025limiting} in the study of a heavy-tailed analogue of the Marchenko--Pastur law for sample correlation matrices. Although $\mu_{\alpha,p}$ does not represent a local alternative to the uniform distribution, its geometry is close to that of $\mu_0$ in the sense that the sample points are still nearly orthogonal in high dimensions. The difference is that, under $\mu_{\alpha,p}$, a small number of points are either unusually close to each other or nearly aligned along straight lines through the origin. As a consequence, tests based on a single polynomial of the inner products, such as the Rayleigh and Bingham tests, fail to be consistent.

By \cite{cohen2020heavy} (see the discussion after Theorem~4.1), it follows that
\[
p^{1/\alpha} \, \mathbf{X}^\top \mathbf{Y} \;\xrightarrow{d}\; Z_{\alpha}
\]
for some non-degenerate random variable $Z_{\alpha}$ that can be written as the ratio of independent stable random variables. Since $\alpha \in (0,2)$, we have $1/\alpha - 1/2 > 0$, and hence $p^{1/\alpha - 1/2}/2 \to \infty$ as $p \to \infty$. Therefore,
\begin{align*}
    d \bigl( \mu_{\alpha,p}, \mu_0 \bigr) 
    &= \sup_{t  \in [-1,1]} 
       \left|  
       \mathbb{P}_{\mu_{\alpha,p}}\!\left( \mathbf{X}^\top \mathbf{Y} \le t \right) 
       - \mathbb{P}_{\mu_0}\!\left( \mathbf{X}^\top \mathbf{Y} \le t \right) 
       \right| \\
    &\ge 
    \left|  
    \mathbb{P}_{\mu_{\alpha,p}}\!\left( 
        p^{1/\alpha} \mathbf{X}^\top \mathbf{Y} \le \frac{p^{1/\alpha - 1/2}}{2} 
    \right) 
    - \mathbb{P}_{\mu_0}\!\left( 
        \sqrt{p}\,\mathbf{X}^\top \mathbf{Y} \le \frac{1}{2} 
    \right)  
    \right|.
\end{align*}
As $p \to \infty$, the first probability converges to $1$, while the second converges to $\Phi(1/2)$. Therefore, for all sufficiently large $p$,
\[
n \cdot d \bigl( \mu_{\alpha,p}, \mu_0 \bigr) 
    \;\ge\; n\cdot \frac{1 - \Phi(1/2)}{2} \to \infty.
\]
Thus, the proposed test \eqref{test-phi} is also consistent in this model. In contrast, it is known that both the Rayleigh test and the Bingham test are not consistent in this model  while the Packing test is consistent \cite{jiang2025asymptotic}.

\subsection{Alternatives generated from the spherical $2$-designs}

In this subsection, we construct a class of alternatives such that the Rayleigh test, the Bingham test, and the Packing test are all blind to it. The construction is based on a spherical \(2\)-design on \(\mb{S}^{p-1}\) \cite{bondarenko2013optimal,bannai2009survey} with \(p+1\) points. Recall that a set of points \(\la \bm{x}_1,\dots,\bm{x}_N \ra \subset \mb{S}^{p-1}\) is called a {\it spherical \(t\)-design} if, for every multivariate polynomial \(f\) of degree no greater than \(t\),
\[
\mb{E}_{\mu_N} f = \mb{E}_{\mu_0} f
\]
where \(\mu_N\) denotes the empirical distribution over the set \(\la \bm{x}_1,\dots,\bm{x}_N \ra\) and $\mu_0$ is the uniform distribution on $\mb{S}^{p-1}$.

Constructing deterministic sets that are spherical designs is a classical and difficult problem in numerical analysis and combinatorics \cite{bondarenko2013optimal,bannai2009survey}. We will only need the existence of a spherical $2$-design on $\mb{S}^{p-1}$ of size $p+1$ such that every pair of distinct points has the same inner product, equal to $-1/p$, which is easy to construct. We first claim that for any $p \geq 3$, there exist $p+1$ points $\bm{v}_1,\dots,\bm{v}_{p+1}$ on $\mb{S}^{p-1}$ such that
\[
\bm{v}_r^\top \bm{v}_s=
\begin{cases}
1, & r=s,\\[1mm]
-\dfrac1p, & r\neq s.
\end{cases}
\]
To see this, define
\[
\mathbf 1:=(1,\dots,1)^\top\in\mathbb R^{p+1},
\qquad
H:=\Bigl\{ \bm{x}\in\mathbb R^{p+1}:\sum_{k=1}^{p+1}x_k=0\Bigr\}.
\]
Then $H$ is a $p$-dimensional subspace of $\mathbb R^{p+1}$. For $i=1,\dots,p+1$, define
\[
\bm{u}_i:=\bm{e}_i-\frac1{p+1}\mathbf 1\in H,
\qquad
\tilde{\bm{s}}_i:=\sqrt{\frac{p+1}{p}}\,\bm{u}_i.
\]
It is easy to check that 
\[
\tilde{\bm{s}}_i^\top \tilde{\bm{s}}_j=
\begin{cases}
1, & i=j,\\[1mm]
-\dfrac1p, & i\neq j.
\end{cases}
\]
Since $H$ is $p$-dimensional, there is an orthogonal transformation that maps $H$ to $\la \bm{x} \in \mb{R}^{p+1}: x_{p+1}=0 \ra$. One can then identify this space with $\mb{R}^p$, and since orthogonal transformations preserve both distance and inner products, the images of $\tilde{\bm{s}}_1,\dots,\tilde{\bm{s}}_{p+1}$ can be viewed as a set of $p+1$ points on $\mb{S}^{p-1}$ satisfying the desired properties. 

We are now ready to construct the class of alternatives. Choose $\ve_p:=1/(4p)$ and define the spherical caps
\[
C_{r,p}:=\left\{ \bm{x}\in\mathbb S^{p-1}: \arccos(\bm{x}^\top \bm{v}_r)\le \ve_p\right\},
\qquad r=1,\dots,p+1.
\]
Let \(\nu_{r,p}\) be the uniform distribution on \(C_{r,p}\), and define
\[
\mu_n
:=
\frac{1}{p+1}\sum_{r=1}^{p+1}\nu_{r,p}.
\]
In other words, $\mu_n$ is a mixture of uniform distributions over spherical caps of small width. The use of spherical caps is meant only to ensure that the resulting distribution is continuous, since otherwise samples from $\mu_n$ might have repeated outcomes. Under this class of alternatives, we can show that

\begin{thm} \label{inconsistency}
    Suppose $p/n^2 \to \infty$ and $\mu_n$ is as above with cap width $\ve_p=1/(4p)$. Then the following hold:
    \begin{itemize}
        \item $\liminf_{n \to \infty} d(\mu_n, \mu_0) \geq 1/10$, where $d$ is the distance in \eqref{metric}.
        \item $R_n$ converges to $0$ in probability, where $R_n$ is the Rayleigh test in \eqref{Rayleigh}.
        \item $B_n$ diverges to $-\infty$ in probability, where $B_n$ is the Bingham test in \eqref{Bingham}.
        \item $P_n$ diverges to $-\infty$ in probability, where $P_n$ is the Packing test in \eqref{packing}.
    \end{itemize}
\end{thm}
Theorem \ref{inconsistency} states that all three existing tests are inconsistent against this class of alternatives, whereas the proposed test remains consistent. The intuition behind this class of alternatives is that its moments nearly match the first two moments of the inner products under the uniform distribution, while the observations drawn from it remain nearly orthogonal. The difference between this class and Class \ref{heavy-tailed} is that there is no small subset of points that stay close together or are nearly aligned. This construction might be extendable to higher orders, but we do not know whether there exists a corresponding spherical design that would make the argument work.

\section{Discussion of the distance $d$ and the proposed test} \label{when to use}

In this section, we discuss several advantages of the distance $d$ and the proposed test $T_n$. 
First, the distance $d$ is a measure of ``symmetry'' and differs from classical metrics between probability measures, such as total variation, Hellinger, or chi-squared distance. These distances are not tailored to the orthogonally invariant structure of the problem and, in particular, they do not reflect geometric features such as concentration along lower-dimensional subspaces.

To illustrate the difference from such metrics, consider the class of low-dimensional uniform distributions introduced in Section~\ref{sec low-rank}. 
In terms of total variation distance, we always have
\[
d_{\rm TV}\!\left( \mathrm{Unif}\!\left(\mathbb{S}^{p-1}\right),\, 
\mathrm{Unif}\!\left(H \cap \mathbb{S}^{p-1}\right) \right) = 1
\]
for all subspaces $H$ with dimension less than or equal to $p-1$.

Thus, density-based distances such as total variation are not appropriate for alternatives that are singular with respect to the uniform distribution. 
In contrast, the distance $d$ is sensitive to geometric deviations: it detects that low-dimensional uniform distributions have large “empty regions’’ compared to the standard uniform distribution, and it also captures the concentration patterns of FvML distributions or Watson distributions at the optimal rates.

Regarding testing procedures for the problem \eqref{uniform-test}, a good high-dimensional test should ideally satisfy the following properties:
\begin{enumerate}
    \item In fixed dimensions, it should be universally consistent. In other words, we would like a test that generalizes naturally from fixed-dimensional to high-dimensional settings.
    \item For common high-dimensional parametric models, it should achieve the optimal detection rates.
    \item It should also be able to detect non-local alternatives in high-dimensional settings.
\end{enumerate}

To the best of our knowledge, no previously studied high-dimensional test has been shown to satisfy all of these properties simultaneously. Most of the tests studied in the literature are Sobolev-based tests \cite{Gine,Ley-P,Cutting-P-V,Cutting-P-V2,ebner2025high}. 
For Sobolev tests, however, universal consistency requires an infinite sequence of weights. In particular, tests with only finitely many nonzero coefficients are simpler to compute, whereas tests with nonzero coefficients across all spherical harmonic degrees are consistent against all alternatives, but are generally much less tractable and typically do not admit a simple closed-form expression. This was one of the main motivations behind the data-driven Sobolev procedures of \cite{Bogdan,Jupp}, which automatically truncate the series. By contrast, our procedure admits a simple closed-form test statistic.

Extending Sobolev tests to the high-dimensional setting is also difficult, since it requires a detailed understanding of the spherical Laplacian when the dimension diverges. Recent work \cite{ebner2025high} studies the high-dimensional null distribution of Sobolev tests on hyperspheres and derives their non-null behavior under FvML alternatives, but its analysis does not include Watson-type alternatives. In contrast, the present paper shows that our proposed test is rate-optimal or consistent across all classes of alternatives considered here, including the FvML model, the Watson model, the low-dimensional uniform distributions model, projections of heavy-tailed distributions, and spherical $2$-design alternatives. In this sense, our method provides a more unified notion of optimality than some existing Sobolev-based approaches.

Another interesting feature of the distance $d$, for which we do not yet have a fully general theory, is that when one restricts attention to many parametric, high-dimensional classes of distributions, the threshold at which $n d \lb \mu_n, \mu_0 \rb$ converges to a non-zero limit often coincides with the minimax rate for testing uniformity within that family. One can show this for some other models, such as the spiked covariance distributions. The task of obtaining an asymptotic expansion for the distance $d$ along a parametric model can often be carried out using an Edgeworth-type expansion, similar to that of Proposition \ref{watson d asymptotic}.

The behavior of the four tests can be summarized in Tables \ref{tab:detection-rates} and \ref{tab:nonlocal-detection}. One can see that the proposed test \eqref{test-phi} is the only one that remains consistent or rate-optimal across all the models considered here. All the other tests are sub-optimal or inconsistent for at least one model. We note, however, that their behavior has not been fully analyzed in all models and regimes. For example, the Packing test $P_n$ in \eqref{packing} is known to be blind at the minimax threshold in the FvML model, but its precise detection threshold in that model remains unknown. Likewise, for the Watson model, it is also unknown what detection threshold the Packing test can achieve.


\begin{table}[t]
\centering
\footnotesize
\setlength{\tabcolsep}{4pt}
\renewcommand{\arraystretch}{1.15}
\caption{Asymptotic detection boundaries under FvML, Watson, and low-dimensional alternatives.}
\label{tab:detection-rates}
\begin{tabularx}{\textwidth}{|c|X|X|X|}
\hline
\textsf{Test / model} 
& \textsf{FvML} 
& \textsf{Watson}
& \textsf{low-dimensional} \\
\hline
$R_n$ \eqref{Rayleigh}
& $\dfrac{p^{3/4}}{\sqrt{n}}$, \textbf{optimal}; \cite{Cutting-P-V}
& sub-optimal (blind); \cite{Cutting-P-V2}
& $k=o(p)$, sub-optimal \\ 
\hline
$B_n$ \eqref{Bingham}
& $\dfrac{p^{3/4}}{n^{1/4}}$, sub-optimal; \cite{Cutting-P-V}
& $\dfrac{p^{3/2}}{\sqrt{n}}$, \textbf{optimal}; \cite{Cutting-P-V2}
& $k=\bigl(1-\Omega(1/n)\bigr)p$, \textbf{optimal} \\ 
\hline
$P_n$ \eqref{packing}
& blind at $\dfrac{p^{3/4}}{\sqrt{n}}$, sub-optimal; \cite{jiang2025asymptotic}
& unknown
& $k=\bigl(1-\Omega(1/\log n)\bigr)p$, sub-optimal \\ 
\hline
proposed \eqref{test-phi}
& $\dfrac{p^{3/4}}{\sqrt{n}}$, \textbf{optimal}; Proposition~\ref{local power FvML}
& $\dfrac{p^{3/2}}{\sqrt{n}}$ in the regime $p=o(n)$; more generally \eqref{watson scaling limit}; \textbf{optimal}
& $k=\bigl(1-\Omega(1/n)\bigr)p$, \textbf{optimal}; Proposition~\ref{power low-rank} \\ 
\hline
\end{tabularx}
\end{table}

\begin{table}[t]
\centering
\footnotesize
\setlength{\tabcolsep}{4pt}
\renewcommand{\arraystretch}{1.15}
\caption{Consistency under non-local alternatives and in fixed dimensions.}
\label{tab:nonlocal-detection}
\begin{tabularx}{\textwidth}{|c|X|X|X|}
\hline
\textsf{Test / model} 
& \textsf{Projections of heavy-tailed distributions} 
& \textsf{spherical $2$-design alternatives}
& \textsf{universal consistency (fixed $p$)} \\
\hline
$R_n$ \eqref{Rayleigh}
& inconsistent; \cite{jiang2025asymptotic}
& inconsistent; Theorem~\ref{inconsistency}
& no \\ 
\hline
$B_n$ \eqref{Bingham}
& inconsistent; \cite{jiang2025asymptotic}
& inconsistent; Theorem~\ref{inconsistency}
& no \\ 
\hline
$P_n$ \eqref{packing}
& \textbf{consistent}; \cite{jiang2025asymptotic}
& inconsistent; Theorem~\ref{inconsistency}
& no  \\ 
\hline
proposed \eqref{test-phi}
& \textbf{consistent}
& \textbf{consistent}; Theorem~\ref{inconsistency}
& \textbf{yes} \\ 
\hline
\end{tabularx}
\end{table}

\section{Simulation studies and comparison with other tests} \label{sec:simulation}

\subsection{Comparison with projection-based tests} \label{comparison}

At a high level, our newly developed procedure follows the philosophy of projection-based tests, initially developed in \cite{cuesta2009projection}. Their test relies on the characterization \eqref{identifiable-2}. 
This characterization can be shown using a variant of the Cramér-Wold device, although the same argument does not apply to Proposition \ref{identifiable}. Based on \eqref{identifiable-2}, \cite{cuesta2009projection} proposed a test that rejects for large values of  
\begin{align} \label{D_nU}
    D_{n,\bm{U}} := \sup_{x \in [-1,1]} \left| \frac{1}{n}\sum_{i=1}^{n} \mathbb{I}\left\{ \bm{X}_{i}^\top\bm{U} \leq x \right\} - \mathbb{P}\left( \bm{e}_{1}^\top\bm{U} \leq x \right) \right|,
\end{align}
where \(\bm{U} \sim \text{Uni}(\mathbb{S}^{p-1})\) is drawn independently from the data and $\bm{e}_1=(1,0,\dots,0)$. 

The test in \cite{cuesta2009projection} uses the same critical value as the Kolmogorov--Smirnov test. In practice, \(\bm{U}\) is drawn multiple times from $\mbox{Unif}\lb \mb{S}^{p_n-1} \rb$ and one gets a corresponding $p$-value for every such $\bm{U}$. The test rejects if the smallest \(p\)-value is below a threshold. More specifically, one picks a large number $k$ and draws $\bm{U}_1,\dots,\bm{U}_k$ independently from $\mbox{Unif}\lb \mb{S}^{p_n-1} \rb$. The test in \cite{cuesta2009projection} rejects at $\alpha$-level if 
\[
\min_{1 \leq i \leq k} \mb{P} \lb D_{n,\bm{U}_i} > K_{\alpha} \Big| \bm{X}_{1},\dots,\bm{X}_{n} \rb \leq c_{\alpha}
\]
where $K_{\alpha}$ is the critical value of the Kolmogorov--Smirnov test and $c_{\alpha}$ is the $(1-\alpha)$-quantile of the left-hand side. However, as the asymptotic theory for this test remains unresolved, computationally intensive Monte Carlo methods are often required to approximate $c_{\alpha}$.

Subsequent works, such as \cite{escanciano2006consistent,garcia2023projection,garcia2021cramer} (see also the references therein), addressed this issue by integrating over all possible directions \(\bm{U}\), resulting in test statistics of the form
\begin{align} \label{L2-test}
    \mathbb{E}_{\bm{U}} \left[ \int_{-1}^{1} \left( \frac{1}{n}\sum_{i=1}^{n} \mathbb{I}\left\{ \bm{X}_{i}^\top\bm{U} \leq x \right\} - \mathbb{P}\left( \bm{e}_{1}^\top\bm{U} \leq x \right) \right)^2 w(x)dx \right],
\end{align}
for some weight function \(w(x) \in L^2\lb [-1,1] \rb\), where the expectation above is taken with respect to $\bm{U}$. 

Test statistics like \eqref{L2-test} exhibit desirable properties, similar to the Anderson--Darling and Cramér--von Mises tests. However, the expectation with respect to \(\bm{U}\) in \eqref{L2-test} often lacks a closed-form expression, requiring Monte Carlo simulations for approximation. Additionally, their asymptotic distributions frequently involve weighted sums of chi-squared distributions, complicating the computation of tail probabilities. For example, the tests in \cite{garcia2021cramer,garcia2023projection} rely on Imhof's method to approximate the critical value.  Our proposed test offers two key advantages over projection-based tests:

\begin{enumerate}
    \item \textbf{Simplicity:} Unlike projection-based methods, our test avoids sampling random directions or integrating over all possible directions, which often requires complex procedures to approximate the critical values. Theorem \ref{null-dist} demonstrates that when the dimension is large, the tail probabilities of our test statistic are much simpler to approximate, eliminating the need for Monte Carlo simulation.

    \item \textbf{Adaptivity over different types of alternatives:}  As demonstrated by our analysis, the proposed test is the only one among the existing high-dimensional tests that is provably minimax-optimal across the parametric models containing the null while also being able to detect non-local alternatives. We are not aware of any theory for projection-based tests in high-dimensional settings. 
\end{enumerate}

\subsection{Simulation studies}

In this subsection, we conduct simulation studies to compare the power of $T_n$ with that of the three existing tests, $R_n$, $B_n$, and $P_n$, defined in \eqref{Rayleigh}, \eqref{Bingham}, and \eqref{packing}, respectively. In the experiments, we set $n=p=80$ for the FvML model and the low-dimensional uniform distributions model. For the Watson model, we use $p=600$ and $n=400$ to get better empirical result, because the convergence seems to be slow in this case.  The empirical powers of the four tests are reported in Figure \ref{fig-power} below. The red curve corresponds to the proposed test $T_n$, while the dashed curve represents the asymptotic power predicted by our theoretical results. We see that the empirical power agrees well with the theoretical analysis: the proposed test is the only one that exhibits strong power across the three high-dimensional parametric models considered. In particular, The Bingham test is blind to FvML local alternatives at the minimax scale. The Rayleigh test is blind to both the low-dimensional alternatives and the Watson alternatives. The Packing test appears to be blind to both the FvML and Watson alternatives.

\begin{figure}[H]
    \centering
    \includegraphics[width=9cm,height=8cm]{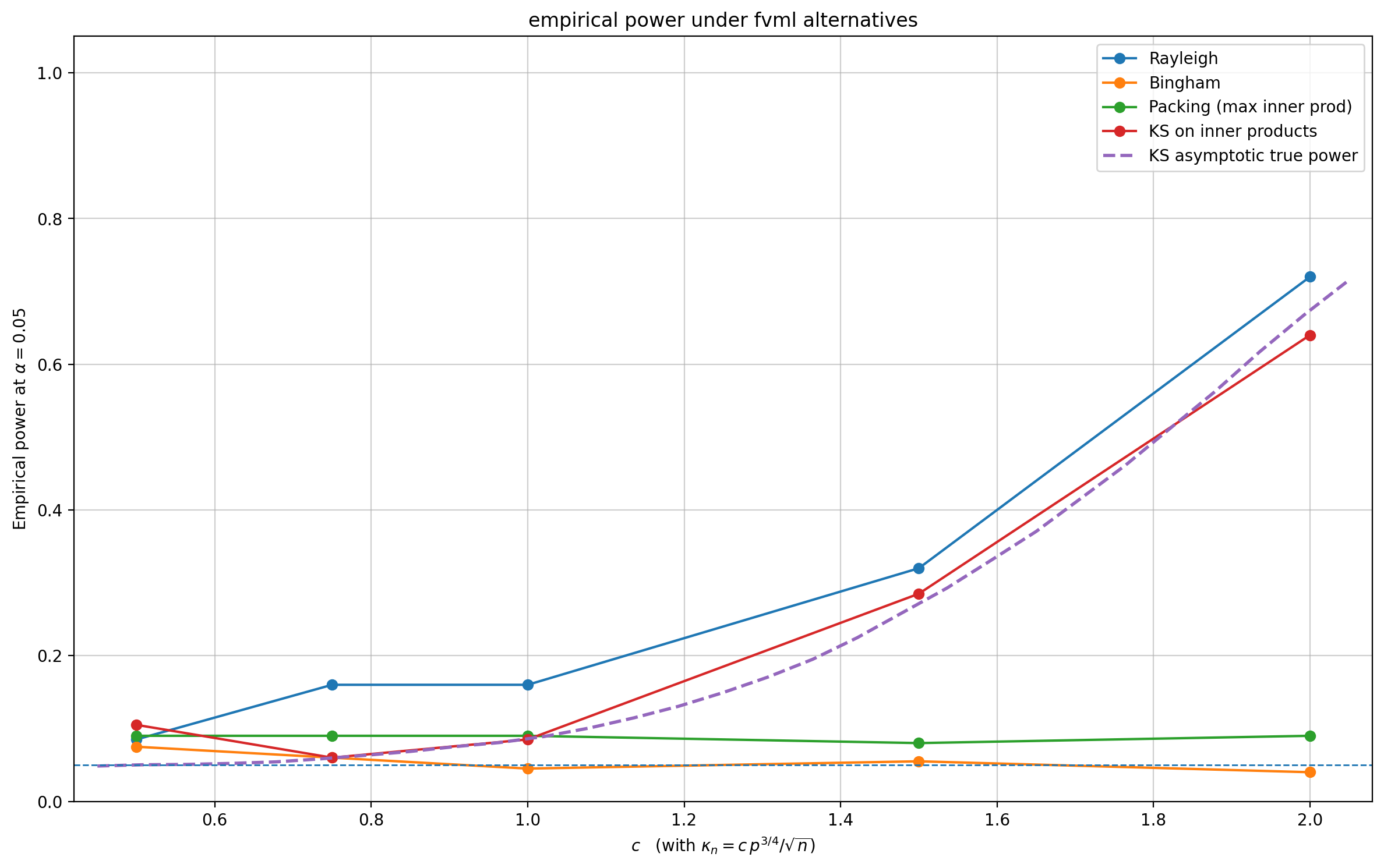}
    \includegraphics[width=9cm,height=8cm]{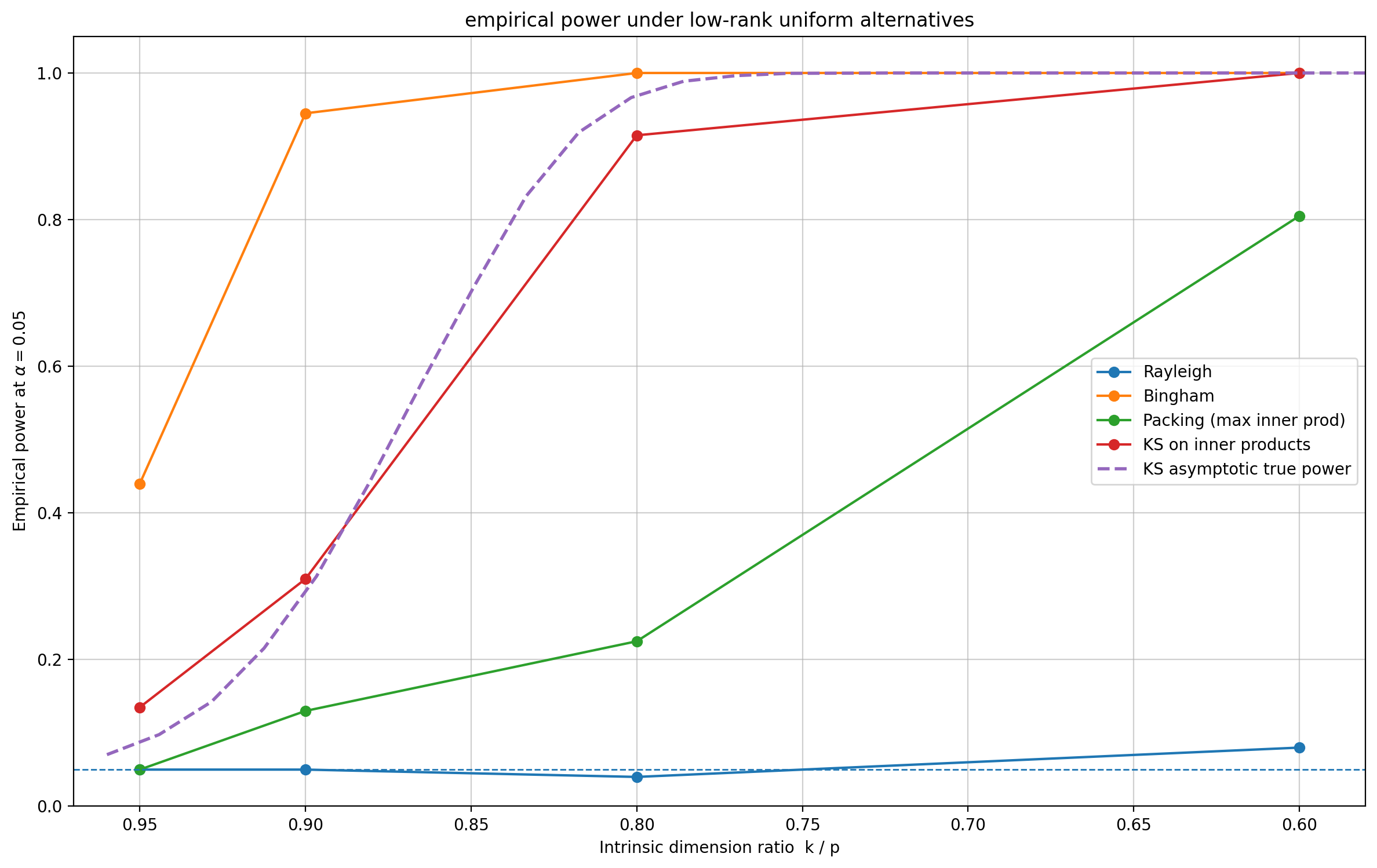}
    \caption{Empirical power under the FvML model and the low-dimensional model}
    \label{fig-power}
\end{figure}

\begin{figure}[H]
    \centering
    \includegraphics[width=9cm,height=8cm]{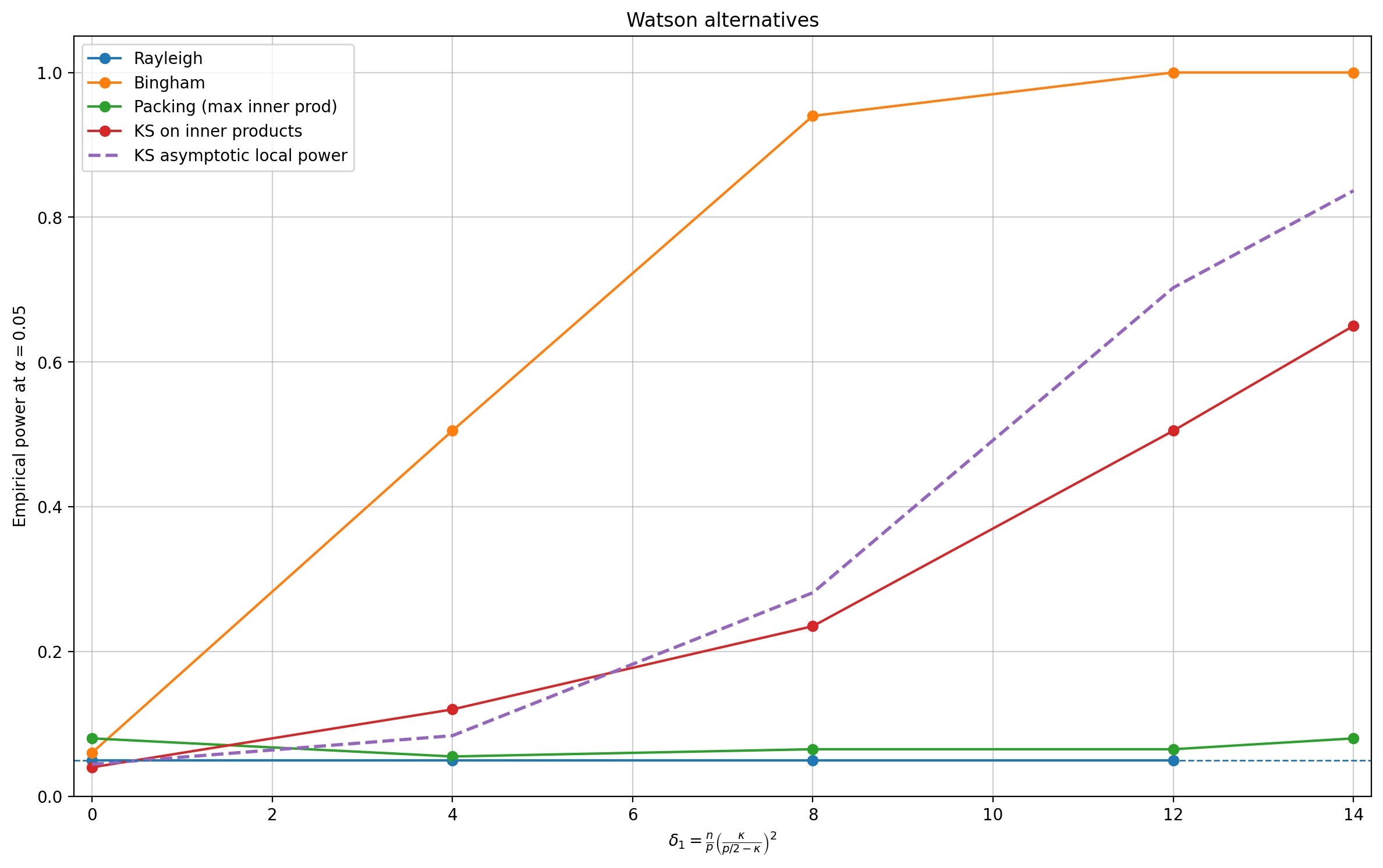}
    \caption{Empirical power under the Watson model}
    \label{fig-power2}
\end{figure}

\section{Conclusions and remarks} \label{discussion}

In this paper, we propose a novel distance to quantify deviations from uniformity, together with a test naturally associated with this distance. 
We show that the test enjoys very simple asymptotic properties in high dimensions and admits a nonparametric consistency theory. 
We establish optimal detection rates with respect to the proposed distance and show that the test attains these rates. 
Furthermore, we show that, when restricted to various parametric models, the proposed distance precisely captures the minimax rates for testing uniformity in all of these models; this is verified for the FvML model, the Watson model, and a low-dimensional uniform distribution model. 
As a consequence of our analysis, we also obtain new information-theoretic lower bounds that match the upper bounds. We conclude with a few remarks.

\begin{enumerate}
    \item It is of independent interest to extend Proposition~\ref{identifiable} to other types of spherical distributions. 
    We conjecture that the conclusion of Proposition~\ref{identifiable} remains valid, up to an orthogonal transformation, for any two Borel probability measures on the sphere under suitable regularity conditions.

    \item We believe that the proposed distance characterizes the minimax rates for other models as well. 
    For example, one can show this for certain spiked covariance models, although the analysis in those cases requires specific restrictions on the joint growth of $p$ and $n$. However, the analysis of the local power is expected to be more involved, since the corresponding contiguity properties are not yet well understood. We expect the framework used in the proof of Theorem \ref{watson local distribution}, as well as that of \cite{han2023contiguity}, to be useful in this direction. 

    \item One can also investigate a procedure based on an $L^2$-type distance instead, for which similar results are expected to hold. 
    We leave this as a direction for future work.
\end{enumerate}

\section*{Acknowledgement}
The second author gratefully acknowledge NSF grant 2217069.

\section{Proofs} \label{proofs}


\subsection{Proof of Theorem \ref{null-dist}}
 
Define
\begin{align} \label{S_nt}
S_n(t):= \sqrt{\frac{2}{n(n-1)}} \sum_{i<j} \left[ \mathbf{1}_{\la \sqrt{p} \bm{X}_i^\top \bm{X}_j \leq t \ra} - \mb{P}_{\mu_0} \lb  \sqrt{p} \bm{X}_1^\top \bm{X}_2 \leq t \rb \right].
\end{align}
We will show that the process $\la S_n(t); n \geq 1 \ra$ converges in distribution to $\la B_{\Phi(t)} \ra$ in the Skorohod space $D[a,b]$, for all $a<b \in \mb{R}$. Some basic properties of the topology on this space can be found in \cite{Dehling}.

Here $\Phi$ is the CDF of a standard normal distribution, and $\la B_{u}; 0 \leq u \leq 1 \ra$ is the Brownian bridge. We do not work directly with the space $D\lb \mb{R} \rb$ (see \cite{vogel2010weak} for more details) since the supremum functional is not almost surely continuous on this space; see Remark \ref{D(R) continuity} below.

{\it \underline{Step 1: Convergence in $D[a,b]$}.} Suppose $a<b$. To show convergence in $D[a,b]$, we need to check the following conditions.

\begin{cond}[Finite-dimensional convergence in distribution] \label{cond-fd}
For any grid $a \leq t_1 < t_2 < \dots <t_k \leq b$, one has $(S_n(t_1), S_n(t_2), \dots, S_n(t_k))$ converges in distribution to $\lb B_{\Phi(t_1)}, B_{\Phi(t_2)},\dots, B_{\Phi(t_k)} \rb$. 
\end{cond}

\begin{cond}[Tightness]
    For any $\ve >0$, we have
    \begin{align} \label{tightness}
       \lim_{\delta\to 0} \limsup_{n \to \infty} \mb{P} \lb \sup_{|t-s| \leq \delta} \Big| S_n(t)-S_n(s)  \Big| > \ve \rb =0.
    \end{align}
\end{cond}





 To check Condition \ref{cond-fd}, we will make use of the following result, whose proof is given in Section \ref{proof-clt Ustat}.

\begin{prop} \label{clt for Ustat}
    Let $h_n: \mb{R} \mapsto \mb{R}$ be a sequence of measurable functions such that $\mb{E} h_n(\bm{X}^{\top}_1 \bm{X}_2)=0$ and
\begin{align}
    \mbox{Var}\lb h_n(\bm{X}^{\top}_1 \bm{X}_2) \rb &\to \sigma^2 >0;  \label{conv-var} \\
        \frac{\mb{E} \lb h_n^4(\bm{X}^{\top}_1 \bm{X}_2) \rb}{n  \lb \mb{E}  h_n^2(\bm{X}^{\top}_1 \bm{X}_2) \rb ^2} &\to 0.\label{lindeberg}
\end{align}
Then,
\begin{align} \label{reduced-fd}
  \sqrt{\frac{2}{n(n-1)}} \sum_{1 \leq i<j \leq n} h_n \lb \bm{X}^{\top}_i \bm{X}_j \rb \xrightarrow{d} N(0,\sigma^2).
\end{align}
\end{prop}
By applying Proposition \ref{clt for Ustat} to kernels of the form 
$$h_n(x)=\mathbf{1}_{\la  \sqrt{p} \cdot x \leq t \ra} - \mb{P}_{\mu_0} \lb \sqrt{p} \bm{X}_1^\top \bm{X}_2 \leq t \rb,$$ 
we obtain the convergence of finite-dimensional distributions. Note that condition \eqref{lindeberg} is satisfied because $\sqrt{p}\bm{X}_1^\top \bm{X}_2$ is asymptotically standard normal. 

We now check the tightness condition \eqref{tightness}. Note that under uniformity,
\[
c_{n,\delta}:= \sup_{|t-s| \leq \delta} \mb{P} \lb s< \sqrt{p} \bm{X}_1^\top \bm{X}_2 <t \rb = \sup_{|t-s| \leq \delta} \mb{P} \lb s \leq  Z \leq t \rb + O \lb p^{-1/2} \rb
\]
as $p \to \infty$. 

By applying Lemma \ref{U-process concentration} to the class of functions $\la \mathbf{1}_{\la s\leq  \sqrt{p} \cdot  \bm{x}^\top \bm{y} \leq t \ra } \ra_{a\leq s<t \leq b}$ (which has $\mbox{VC}$ dimension $2$) and using the degeneracy of the kernels, we obtain
\begin{align*}
    \mb{P} \lb \sup_{|t-s| \leq \delta} \mid S_n(t) - S_n(s) \mid \geq \ve \rb \lesssim c_{n,\delta} \left[1 + \log \lb c_{n,\delta} \rb \right] + \frac{\log \lb c_{n,\delta}^{-1} \rb}{\sqrt{n}}.
\end{align*}
The proof is completed by first letting $n \to \infty$ and then letting $\delta \to 0$.

{\it \underline{Step 2: Continuous mapping and negligibility of the tail}.} By the continuous mapping theorem, for all $a>0$, we have 
\[
\sup_{t \in [-a,a]} |S_n(t)| \stackrel{d}{\to} \sup_{t \in [-a,a]} |B_{\Phi(t)}|.
\]
To deduce the result, it suffices to show that for all $\ve>0$
\begin{align} \label{tail-negligible}
\lim_{a \to \infty} \mb{P} \lb \sup_{|t|>a}|S_n(t)|>\ve  \rb = 0.
\end{align}
The above is equivalent to showing that
\[
\lim_{a \to \infty} \mb{P} \lb \sup_{t>a}|S_n(t)|>\ve  \rb = 0, \ \text{and} \  
\lim_{a \to \infty} \mb{P} \lb \sup_{t<-a}|S_n(t)|>\ve  \rb = 0.
\]
Since the proofs of these two limits are identical, we will only prove the former. We again apply Lemma \ref{U-process concentration} to the VC-type class of functions
\[
\la \mathbf{1}_{\la \sqrt{p} \cdot  \bm{x}^\top \bm{y} \leq t \ra}  \ra_{t \in (a,\infty)}
\]
to deduce that
\[
\mb{P} \lb \sup_{t>a}|S_n(t)|>\ve  \rb  \lesssim \tau_a \left[ 1 + \log \lb \tau_a \rb \right] + \frac{1}{\sqrt{n}}
\]
where the variance profile $\tau_a$ is defined as
\[
\tau_a:= \sup_{t>a} \la \mbox{Var} \lb \mathbf{1}_{\la \sqrt{p} \bm{X}_1^\top \bm{X}_2 \leq t \ra} \rb \ra \leq  1 - \mb{P} \lb  \sqrt{p}\cdot\bm{X}_1^\top \bm{X}_2 \leq a \rb.
\]
It is easy to check that $\tau_a$ converges to $0$ as $a \to \infty$. This finishes the proof. $\hfill$ $\square$

\begin{remark} \label{D(R) continuity}
    The reason we do not work directly with $D( \mb{R})$ is that the functional
\begin{align*}
\mathcal{S}: D(\mb{R}) &\to \mb{R}^+ \\
 f &\to \sup_{t \in \mb{R}} |f(t)|
\end{align*}
is not almost surely continuous at $\la B_{\Phi(u)} \ra_{u \in \mb{R}}$.

In fact, the topology on $D(\mb{R})$ is equivalent to the coarsest topology such that the projection map from $D(\mb{R})$ to $D[a,b]$ is continuous for all $a<b$ (see Section 3 in \cite{vogel2010weak} for more details). Since this topology only sees the behavior of the process on bounded intervals, modifying the process on a diverging sequence still yields convergence in $D(\mb{R})$, but the supremum functional can blow up. 
\end{remark}


\subsection{Proof of Theorem \ref{consistency}}

Consider a sequence of laws $\mu_n$ such that $nd \lb \mu_n, \mu_0 \rb \to \infty$. Put $h_t(\bm{x},\bm{y})=\mathbf{1}_{\la  \langle \bm{x}, \bm{y}  \rangle \leq t \ra}$ and define
\[
 g_{t,n} (\bm{x}):= \mb{E}_{\mu_n} \lb h_t \lb \bm{x}, \bm{y} \rb \mid \bm{x} \rb.
 \]
For $t \in [-1,1]$, rewrite $T_n(t)$ in terms of the Hoeffding projection as 
 \begin{align*}
     T_n(t)= T_{n,1}(t) +T_{n,2}(t) +  d_t
 \end{align*}
where
\begin{align*} 
T_{n,1}(t) :&=   \frac{2}{n}  \sum_{i=1}^n \Big[ g_{t,n} \lb \bm{X}_i \rb - \mb{E} g_{t,n} \lb \bm{X}_i \rb  \Big];    \\
T_{n,2}(t):&=  \frac{2}{n(n-1)}\sum_{1\leq i<j \leq n} \Big[  h_{t} \lb \bm{X}_i, \bm{X}_j \rb -  g_{t,n} \lb \bm{X}_i \rb - g_{t,n} \lb \bm{X}_j \rb + \mb{E}h_{t} \lb \bm{X}_i, \bm{X}_j \rb \Big]; \\
d_t:&= \mb{P}_{\mu_n} \lb \bm{X}_1^\top \bm{X}_2 \leq t \rb - \mb{P}_{\mu_0} \lb \bm{X}_1^\top \bm{X}_2 \leq t \rb.
\end{align*}
Define
\begin{align*}
V_n:&= \max_{t \in [-1,1]} \mbox{Var} \lb g_{t,n} \lb \bm{X}_1 \rb \rb; \\
t_n:&= \mbox{argmax}_{t \in [-1,1]} \mbox{Var} \lb g_{t,n} \lb \bm{X}_1 \rb \rb; \\
\alpha_n:&= \mbox{argmax}_{t\in [-1,1]} \left|  d_t    \right|;   \\
v:&= \limsup_{n \to \infty} ( nV_n ). 
\end{align*}
Roughly speaking, $t_n$ is the point at which the main contribution from the Hoeffding projection term $T_{n,1}$ is attained, and $\alpha_n$ is the point at which the main contribution from the deterministic perturbation $d_t$ comes. 

It suffices to consider the following two cases.

{\it \underline{Case 1:  $v< \infty$}.} In this case, we can estimate 
\begin{align*}
    n \cdot \sup_{t \in [-1,1]}|T_n(t)| &\geq n\cdot |T_n(\alpha_n)| - n\cdot \left| T_{n,1} \lb \alpha_n \rb + T_{n,2}(\alpha_n) \right| \\
    &= n\cdot d\lb \mu_n, \mu_0 \rb - n\cdot \left| T_{n,1} \lb \alpha_n \rb + T_{n,2}(\alpha_n) \right|.
\end{align*}
It is easy to check that 
\begin{align*}
\mbox{Var} \lb T_{n,1} \lb \alpha_n \rb  \rb \leq  \frac{V_n}{n} \ \ \text{and} \ \
\left| T_{n,2}(\alpha_n) \right| \leq \sup_{t \in [-1,1]} \left| T_{n,2}(t) \right| = O_{\mb{P}} \lb n^{-1} \rb
\end{align*}
where the second inequality follows from Lemma \ref{U-process concentration}. 

Consequently,
\begin{align*}
      n \cdot \sup_{t \in [-1,1]}|T_n(t)| &\geq       n |T_n(\alpha_n)| \\
      &\geq  nd_{\alpha_n} - n \left|  T_{n,1} \lb \alpha_n \rb \right| - n \left|  T_{n,2} \lb \alpha_n \rb \right|  \\
      &= n\cdot d\lb \mu_n, \mu_0 \rb - O_{\mb{P}} \lb \sqrt{nV_n} \rb - O_{\mb{P}}(1).
\end{align*}
Since $\sup_{n \geq 1} \lb nV_n \rb = v< \infty$ and $nd\lb \mu_n,\mu_0 \rb \to \infty$, the test rejects with probability tending to one.

{\it \underline{Case 2:  $v = \infty$}.}    In this case, by a subsequence argument, we may assume that $nV_n \to \infty$. Notice that
\begin{align*}
\mb{P} \lb \text{reject the null hypothesis}  \rb &= \mb{P} \lb   \sup_{t \in [-1,1]} |T_n(t)| \geq \frac{q_{\alpha}}{n}(1+o(1)) \rb \\
&\geq \mb{P} \lb \left| T_{n,1}(t_n) + d_{t_n}  \right| \geq  \frac{q_{\alpha}}{n}(1+o(1)) + \left| T_{n,2} \lb t_n \rb \right| \rb \\
&= \mb{P} \lb \left| \sqrt{\frac{n}{V_n}} \cdot T_{n,1}(t_n) + d_{t_n}\sqrt{\frac{n}{V_n}}  \right| \geq T^{*}_n \rb
\end{align*}
where
\[
T^{*}_n:=  \frac{q_\alpha(1+o(1))}{\sqrt{nV_n}} + \sqrt{\frac{n}{V_n}} \left| T_{n,2} \lb t_n \rb \right|,
\]
and the second line follows from the fact that
\[
 \sup_{t \in [-1,1]} |T_n(t)|  \geq \left| T_{n,1}(t_n) + d_{t_n}  \right| - \left| T_{n,2} \lb t_n \rb \right| .
\]
By Lemma \ref{U-process concentration} and the assumption that $nV_n \to \infty$, we deduce that $T_n^*=o_{\mb{P}}(1)$. Now, thanks to the Berry--Esseen bound for sums of i.i.d.\ random variables (see, for example, Theorem 3.7 in \cite{chen2010normal}) and the fact that $|g_{n,t_n}| \leq 1$, we have
\begin{align*}
\sup_{x \in \mathbb{R}}
\left|
\mathbb{P}\!\left(
\sqrt{\frac{n}{V_n}} \, T_{n,1}(t_n) \le x
\right)
-
\mathbb{P}\!\left( N(0,4) \le x \right)
\right| & \lesssim \frac{ \mb{E} \left| g_{n,t_n} \lb \bm{X}_i \rb  \right|^3}{\sqrt{n} \cdot V_n^{3/2}} \\
&\lesssim \frac{ \mb{E} \left| g_{n,t_n} \lb \bm{X}_i \rb  \right|^2}{ \sqrt{n} \cdot V_n^{3/2}}  \\
&= \frac{V_n}{\sqrt{n} \cdot V_n^{3/2}} = \frac{1}{\sqrt{n V_n}} \to 0.
\end{align*}
The proof in this case is completed by employing Lemma \ref{anti concentration} below with
\[
X_n= \sqrt{\frac{n}{V_n}} \, T_{n,1}(t_n); \quad Y_n= T_n^*; \quad a_n=  d_{t_n}\sqrt{\frac{n}{V_n}}.
\]
to get 
\[
\lim_{n\to \infty} \mb{P} \lb \left| \sqrt{\frac{n}{V_n}} \cdot T_{n,1}(t_n) + d_{t_n}\sqrt{\frac{n}{V_n}}  \right| \geq T^{*}_n \rb = 1.
\]
$\hfill$ $\square$ \medskip

\begin{lemma} \label{anti concentration}
Suppose $\{X_n\}$ is a sequence of random variables such that
\[
\sup_{t \in \mb{R}} \left|   \mb{P} \lb X_n \leq t \rb - \mb{P} \lb N(0,4) \leq t \rb \right| \to 0
\]
where $N(0,4)$ is a normal distribution with variance $4$. 

Let $\{a_n\}$ be any sequence of real numbers (not necessarily bounded), and let $\{Y_n\}$ be a sequence of random variables such that
$
Y_n \stackrel{\mb{P}}{\to} 0
$.
Then
\[
\lim_{n\to\infty}
\mathbb{P}\!\left( |X_n + a_n| \ge |Y_n| \right) = 1.
\]
\end{lemma}

\noindent \textbf{Proof of Lemma \ref{anti concentration}.}   
Fix $\ve>0$ and write 
\begin{align*}
     \mathbb{P}\!\Big( |X_n + a_n| \leq |Y_n| \Big)    &\leq      \mathbb{P}\!\Big( |X_n + a_n| \leq |Y_n|, |Y_n| \leq \ve \Big) + \mb{P} \Big( |Y_n| >\ve \Big)   \\
     &\leq \mb{P} \lb |X_n+a_n| \leq \ve \rb + \mb{P} \Big( |Y_n| >\ve \Big) \\
     &= \mb{P} \lb -\ve-a_n \leq N(0,4) \leq \ve-a_n \rb + \mb{P} \Big( |Y_n| >\ve \Big) + o(1) \\
     &\leq 2\ve \cdot \sup_{t \in \mb{R}} \la  \frac{1}{2\sqrt{2 \pi}} \exp \lb -t^2/8 \rb \ra + \mb{P} \Big( |Y_n| >\ve \Big) + o(1). \\
     &\leq  2\ve + \mb{P} \Big( |Y_n| >\ve \Big) + o(1).
\end{align*}
The proof is completed by taking $n \to \infty$ and then letting $\ve \to 0$. $\hfill$ $\square$

\subsection{Proof of Theorem \ref{information bound}}

Fix $\ve>0$ sufficiently small and define
\[
\kappa_n:= \frac{\ve p^{3/4}}{\sqrt{n}}; \quad \frac{d\mu_{n,\bm{\theta}}}{d \mu_0} (\bm{x}) \propto \exp \lb \kappa_n \lb \bm{x}^\top \bm{\theta} \rb \rb.
\]
In other words, $\mu_n$ is an FvML distribution with location $\theta$ and concentration parameter $\kappa_n$. Consider the least favorable distribution $\mu_n^*$ defined as
\begin{align} \label{FvML randomized location}
\frac{d\mu_n^{*}}{d \mu_0^{\otimes n}} \lb \bm{X}_1,\dots,\bm{X}_n \rb:= \mb{E}_{\bm{\theta} \sim \mu_0} \left[  \prod_{i=1}^n \frac{d\mu_{n,\bm{\theta}}}{d \mu_0} (\bm{X}_i)  \right].
\end{align}
By Proposition \ref{Kolmogorov expansion fvml}, for all $\bm{\theta} \in \mb{S}^{p-1}$, we have 
\[
d \lb \mu_{n,\bm{\theta}}, \mu_0 \rb \geq \frac{\ve^2}{10n}
\]
whenever $\min \la p,n \ra$ is sufficiently large and $\ve$ is small enough.  

Consequently, Le Cam's mixture argument yields
\begin{align*}
\liminf_{n \to \infty} \la  \inf_{\phi \in \mathcal{T}_{n}} \la \mb{P}_{\mu_0} \lb \phi=1 \rb + \sup_{d(\nu,\mu_0) \geq \frac{\ve}{n}}  \mb{P}_\nu \lb \phi =0 \rb \ra \ra &\geq  \liminf_{n \to \infty} \Big\{ 1 - d_{\rm TV} \lb \mu_0, \mu_n^{*}  \rb \Big\} \\
&\geq  \liminf_{n \to \infty} \Big\{ 1 - \sqrt{\mb{E} L_n^2 -1 } \Big\}. 
\end{align*}
where $L_n$ is the likelihood ratio defined in \eqref{fvml LLR}.

Thanks to Proposition \ref{2nd moment LLR}, we know that $\mb{E} L_n^2 -1 \leq e^{\ve^2}-1= O(\ve^2)$ for small $\ve>0$. Thus, we get \eqref{minimax} by choosing $\ve$ small enough. The proof is completed. $\hfill$ $\square$


 \subsection{Proof of Theorem \ref{low-rank information bound}}

Let us start with a useful result for calculating likelihood ratios between distributions that are invariant under group actions. For terminology related to group actions and maximal invariants, we refer the reader to Chapters~2 and~3 of the monograph \cite{eaton1989group}. 

For the reader's convenience, we briefly recall the relevant concepts. A group $G$ is said to act on a space $\mathcal{X}$ if there exists a mapping $G \times \mathcal{X} \to \mathcal{X}$ that is compatible with the group operation. A measurable mapping $T : \mathcal{X} \to \mathcal{Y}$ is called an invariant if
\[
T(x) = T(gx), \qquad \forall\, g \in G.
\]
An invariant $T$ is called a maximal invariant if, whenever $T(x) = T(y)$ for some $x,y \in \mathcal{X}$, there exists $g \in G$ such that $x = gy$.

\begin{lemma}\label{group action}
Let $\mathcal{X}$ be a Polish space, and suppose a compact group $G$ acts on $\mathcal{X}$ continuously. 
Let $\mathbb{P}$ and $\mathbb{Q}$ be two Borel probability measures on $\mathcal{X}$ that are invariant under the action of $G$. Let $T : \mathcal{X} \to \mathcal{Y}$ be a continuous maximal invariant for some Polish space $\mathcal{Y}$. 
Define the induced laws
\[
\mathbb{P}_T := \mathbb{P} \circ T^{-1},
\qquad
\mathbb{Q}_T := \mathbb{Q} \circ T^{-1}.
\]
Then $\mathbb{P} \ll \mathbb{Q}$ whenever $\mathbb{P}_T \ll \mathbb{Q}_T$. 
Moreover, when this holds and $X \sim \mathbb{Q}$, we have
\[
\frac{d\mathbb{P}}{d\mathbb{Q}}(X)
=
\frac{d\mathbb{P}_T}{d\mathbb{Q}_T}\!\bigl(T(X)\bigr)
\quad \mathbb{Q}\text{-almost surely}.
\]
\end{lemma}

The proof of Lemma~\ref{group action} can be found in Appendix~\ref{appendix group action}. 
We now construct the least favorable alternative. 
Let $\Pi_{k,p}$ denote the normalized left Haar measure on the Grassmannian $G(k,p)$ (so that it is a probability measure). 
Define
\begin{align*}
    \mathbb{P}_{0n}
    &:= 
    \underbrace{\mathrm{Unif}\!\left(\mathbb{S}^{p-1}\right) \otimes \cdots \otimes 
    \mathrm{Unif}\!\left(\mathbb{S}^{p-1}\right)}_{n\text{-times}}, \\
    \mathbb{P}_{1n}
    &:=
    \int_{G(k,p)}
    \underbrace{\mathrm{Unif}\!\left(H \cap \mathbb{S}^{p-1}\right) \otimes \cdots \otimes 
    \mathrm{Unif}\!\left(H \cap \mathbb{S}^{p-1}\right)}_{n\text{-times}}
    \, \Pi_{k,p}(dH).
\end{align*}
Roughly speaking, $\mathbb{P}_{0n}$ is the joint distribution of 
$\mathbf{X}_1,\ldots,\mathbf{X}_n$ under $H_0$, 
while $\mathbb{P}_{1n}$ is the law obtained by first sampling a 
$k$-dimensional subspace $H \sim \Pi_{k,p}$ and then sampling
\[
(\mathbf{X}_1,\ldots,\mathbf{X}_n)\mid H 
\stackrel{\text{i.i.d.}}{\sim} 
\mathrm{Unif}\!\left(H \cap \mathbb{S}^{p-1}\right).
\]
From now on, we let $\mathbf{X}$ denote the $p \times n$ data matrix whose columns are 
$\mathbf{X}_1,\ldots,\mathbf{X}_n$. 
We will apply Lemma~\ref{group action} to show that 
$d\mathbb{P}_{1n}/d\mathbb{P}_{0n}$ exists and to derive its explicit form. Note that, although one can also use the Blaschke--Petkantschin formula to compute this integral (see, for example, Chapter~7 of \cite{schneider2008stochastic}), the computation is quite lengthy. Define
 \begin{align*}
 \mc{X}:&= \lb \mb{S}^{p-1} \rb^n; \\
 \mc{Y}:&= \la \mc{C} \in \mbox{Sym}_n \lb \mb{R} \rb: \mc{C} \succ 0 \  \text{and} \ \mc{C}_{ii}=1, \ \forall 1\leq i \leq n \ra.
 \end{align*}
where $\mbox{Sym}_n \lb \mb{R} \rb$ is the set of all symmetric matrices of size $n$. 

By Lemma \ref{maximal invariant}, the map 
 \begin{align*}
 T: \mc{X} &\to \mc{Y} \\
 \bm{X} &\to \bm{X}^\top \bm{X}
 \end{align*}
is a maximal invariant under the action of the group of orthogonal matrices. 

Note that the matrix $\bm{G}:=\bm{X}^\top \bm{X}$ is nothing but the sample correlation matrix without centering by the sample mean. We can then apply Lemma 2.1 in \cite{jiang2019determinant} and Theorem 5.1.3 in \cite{muirhead2009aspects} to get the density of $\bm{G}$ under $\mb{P}_{0n}$ as
\begin{align} \label{density P0}
f \lb \bm{G} \rb \propto  \mbox{det} \lb \bm{G} \rb^{(p-n-1)/2} d\bm{G}.
\end{align}
 We have $p-1$ in the above formula instead of $p-2$ in \cite{muirhead2009aspects} because there is no centering term in $\bm{G}$, and Lemma 2.1 in \cite{jiang2019determinant} asserts that such a difference is in fact equivalent to a one-unit shift in $p$. 
 
 Here $d\bm{G}$ denotes the joint density of the upper-diagonal entries of $\bm{G}$. Equivalently, it can also be regarded as a measure on $\mc{Y}$, defined as the pushforward measure of the Lebesgue measure on an open subset of $\mb{R}^{n(n-1)/2}$ to $\mc{Y}$ via the natural embedding.  

Similarly, under $\mb{P}_{1n}$, the density of $\bm{G}$ is given by
\begin{align} \label{density P1}
f \lb \bm{G} \rb \propto  \mbox{det} \lb \bm{G} \rb^{(k-n-1)/2} d\bm{G}.
\end{align}
It is easy to see that the two laws in \eqref{density P0} and \eqref{density P1} are mutually absolutely continuous. Also, these two densities are well-defined due to our assumption that $n+1 \leq k \leq p$. Thus, Lemma \ref{group action} gives
\begin{align} \label{mc{L}}
\mc{L}_n:=\frac{d \mb{P}_{1n}}{d \mb{P}_{0n}} \lb \bm{X}_1,\dots,\bm{X}_n \rb = \frac{d \mb{P}_{1n} \circ T^{-1}}{d \mb{P}_{0n} \circ T^{-1}} \lb \bm{G} \rb =   C \lb \frac{k-p}{2} \rb \cdot \mbox{det} \lb \bm{G} \rb^{(k-p)/2}.
\end{align}
where the normalizing constant $C(\theta)$ satisfies 
\[
C \lb \theta \rb := \left[ \mb{E}_{\mb{P}_{0n}} \lb \mbox{det} \lb \bm{G} \rb^{\theta}  \rb  \right]^{-1}.
\]
Similarly to the proof of Theorem \ref{information bound}, we only need to show that 
\[
\mb{E}_{\mb{P}_{0n}} \lb \mc{L}_n^2 \rb = 1 + o(1)
\]
whenever $n \lb 1 - k/p \rb \to 0$ with $\mc{L}_n$ defined in \eqref{mc{L}}.

From Lemma 5.2 in \cite{jiang2015likelihood}, we find that 
\begin{align*}
  \mb{E}_{\mb{P}_{0n}} \lb \mbox{det} \lb \bm{G} \rb^{\theta}  \rb =  \left[ \frac{\Gamma \lb \frac{p}{2} \rb}{\Gamma \lb \frac{p}{2} +\theta \rb} \right]^n \cdot \frac{\Gamma_n \lb \frac{p}{2}+\theta \rb}{\Gamma_n \lb \frac{p}{2} \rb}
\end{align*}
for $\theta > -\max \la 1, (p-n)/2  \ra$, where $\Gamma_n$ is the multivariate Gamma function defined as in (5.1) of \cite{jiang2015likelihood}. The specific form of $\Gamma_n$ is not relevant to our proof, as we only need the asymptotic result from Proposition 5.1 of \cite{jiang2015likelihood}. These asymptotic results are collected in Lemma \ref{gamma asymptotic} in Appendix \ref{appendix LLR} below.

Consequently, with $\Delta:=k-p$ we obtain 
\begin{align*}
        \mb{E}_{\mb{P}_{0n}} \lb \mc{L}_n^2 \rb =  \frac{C \lb \frac{\Delta}{2} \rb^2}{C \lb \Delta \rb} &= \left[ \frac{\Gamma \lb \frac{p}{2} \rb}{\Gamma \lb \frac{p}{2} +\Delta \rb} \right]^n \cdot \frac{\Gamma_n \lb \frac{p}{2}+\Delta \rb}{\Gamma_n \lb \frac{p}{2} \rb} \cdot \left[ \frac{\Gamma \lb \frac{p+\Delta}{2} \rb}{\Gamma \lb \frac{p}{2}  \rb} \right]^{2n} \cdot \frac{\Gamma^2_n \lb \frac{p}{2} \rb}{\Gamma^2_n \lb \frac{p+\Delta}{2} \rb} \\
        &= \exp \la  F_n(\Delta) - 2F_n \lb \frac{\Delta}{2} \rb  \ra.
\end{align*}
where
\begin{align} \label{Fn}
F_n(\Delta):= n \log \left[ \frac{\Gamma \lb \frac{p}{2} \rb}{\Gamma \lb \frac{p}{2} +\Delta \rb}  \right]  +\log \left[  \frac{\Gamma_n \lb \frac{p}{2}+\Delta \rb}{\Gamma_n \lb \frac{p}{2} \rb} \right].
\end{align}
The proof is completed by applying Proposition \ref{Fn asymptotic} in Appendix \ref{appendix LLR}, which states that 
\[
F_n(\Delta) - 2F_n \lb \frac{\Delta}{2} \rb \to 0
\]
as $n \to \infty$ such that $n(1-k/p) \to 0$. $\hfill$ $\square$


\section{Proof of Theorem \ref{watson minimax}}
\subsection{Technical preparation}

\begin{lemma}\label{lem:rank-one-full}
Fix $\eta>0$. Then there exist constants $
0<c_\eta<C_\eta<\infty
$
such that for all sufficiently large $p$ and all
\[
0\le t\le \frac p2-\eta\sqrt p,
\]
one has
\begin{align*}
c_\eta \sqrt p\,a^{-1/2}
&\le Z_p(t)\le
C_\eta \sqrt p\,a^{-1/2},\\
c_\eta \sqrt p\,a^{-3/2}
&\le Z_p'(t)\le
C_\eta \sqrt p\,a^{-3/2},\\
c_\eta \sqrt p\,a^{-5/2}
&\le Z_p''(t)\le
C_\eta \sqrt p\,a^{-5/2},
\end{align*}
where
\[
a:=\frac p2-t.
\]
As a consequence, with 
\begin{align} \label{log Zp}
f_p(t):=\log Z_p(t),
\end{align}
we have
\begin{align*}
\frac{c_\eta}{p-2t}
\le \frac{Z_p'(t)}{Z_p(t)}
\le \frac{C_\eta}{p-2t},
\qquad
0\le f_p''(t)\le \frac{C_\eta}{(p-2t)^2}.
\end{align*}
\end{lemma}

\noindent \textbf{Proof of Lemma \ref{lem:rank-one-full}.}
For $r\in\{0,1,2\}$ define
\[
M_r(t)
:=
c_p\int_0^1 e^{tu}u^{r-1/2}(1-u)^{(p-3)/2}\,du,
\qquad
c_p:=\frac{\Gamma(p/2)}{\sqrt\pi\,\Gamma((p-1)/2)}.
\]
Then
\[
Z_p(t)=M_0(t),\qquad Z_p'(t)=M_1(t),\qquad Z_p''(t)=M_2(t).
\]
By the gamma-ratio bound, there exist universal constants $0<c<C<\infty$ such that
\[
c\sqrt p\le c_p\le C\sqrt p
\]
for all $p\ge 3$.

Recall $a=p/2-t\ge \eta\sqrt p$. Split the integral into
\begin{align*}
M_r(t)
&=
c_p\int_0^1
u^{r-1/2}(1-u)^{-3/2}
\exp\!\left(
-au+\frac p2\bigl(u+\log(1-u)\bigr)
\right)\,du \\
&= M^{(1)}_r(t) + M^{(2)}_r(t)
\end{align*}
where
\begin{align*}
    M^{(1)}_r(t):&= c_p\int_0^{1/2}
u^{r-1/2}(1-u)^{-3/2}
\exp\!\left(
-au+\frac p2\bigl(u+\log(1-u)\bigr)
\right)\,du ; \\
M^{(2)}_r(t):&= c_p\int_{1/2}^{1}
u^{r-1/2}(1-u)^{-3/2}
\exp\!\left(
-au+\frac p2\bigl(u+\log(1-u)\bigr)
\right)\,du.
\end{align*}
We will show that for $r \in \la 0,1,2 \ra$,
\[
 M^{(2)}_r(t) = O(\sqrt p\,e^{-c_\eta' p}) \qquad \text{and} \qquad  M^{(1)}_r(t) \asymp \sqrt p\,a^{-(r+1/2)}.
\]
Given this, the proof of the first part of Lemma \ref{lem:rank-one-full} is complete since $a\leq p/2 \ll e^{c_\eta p}$.

\underline{\textit{Analysis of $M^{(2)}_r(t)$.}}
Recall the original representation
\[
M_r(t)=c_p\int_0^1 e^{tu}u^{r-1/2}(1-u)^{(p-3)/2}\,du.
\]
Consider
\[
h(u):=tu+\frac{p-3}{2}\log(1-u), \qquad u\in[1/2,1).
\]
Since $t\le p/2-\eta\sqrt p$, we have for all sufficiently large $p$,
\[
h'(u)=t-\frac{p-3}{2(1-u)}
\le
\frac p2-\eta\sqrt p-(p-3)
<0,
\qquad u\in[1/2,1].
\]
Thus $h$ is decreasing on $[1/2,1]$, and therefore
\[
h(u)\le h(1/2)
=
\frac t2+\frac{p-3}{2}\log\frac12
\le -c_\eta' p
\]
for some $c_\eta'>0$ and all sufficiently large $p$. Moreover, since
\[
u^{r-1/2}\le 1 \qquad \text{on } [1/2,1],
\]
it follows that
\[
\int_{1/2}^1 e^{tu}u^{r-1/2}(1-u)^{(p-3)/2}\,du
\le
e^{-c_\eta' p}.
\]
Using $c_p\asymp \sqrt p$, we conclude that the contribution of $[1/2,1]$ to $M_r(t)$ is
\[
O(\sqrt p\,e^{-c_\eta' p}),
\]
uniformly over $a\ge \eta\sqrt p$.

\underline{\textit{Analysis of $M^{(1)}_r(t)$.}} 
For $u\in[0,1/2]$, one has
\[
-u^2\le u+\log(1-u)\le -\frac{u^2}{2},
\qquad
1\le (1-u)^{-3/2}\le 2^{3/2},
\]
hence
\[
u^{r-1/2}e^{-au-pu^2/2}
\lesssim
u^{r-1/2}(1-u)^{-3/2}
\exp\!\left(
-au+\frac p2(u+\log(1-u))
\right)
\lesssim
u^{r-1/2}e^{-au-pu^2/4}
\]
uniformly on $[0,1/2]$.

Put 
\[
I_r(c,a,p):=\int_0^{1/2} u^{r-1/2}e^{-au-cpu^2}\,du.
\]
We claim that
\[
I_r(1/2, a,p)\gtrsim a^{-(r+1/2)} \qquad \text{and} \qquad I_r(1/4, a,p)\lesssim a^{-(r+1/2)}
\]
uniformly over $p/2 \geq a \geq \eta\sqrt{p}$.

For the upper bound, discard the quadratic term:
\[
I_r(1/4, a,p)\le \int_0^\infty u^{r-1/2}e^{-au}\,du
=
\Gamma\!\left(r+\frac12\right)a^{-(r+1/2)}.
\]

For the lower bound, restrict to $u\in[0,1/a]$. Since $a\ge \eta\sqrt p$, one has
\[
\frac1a\le \eta^{-1}p^{-1/2},
\]
and thus for $0\le u\le 1/a$,
\[
pu^2\le \frac{p}{a^2}\le \eta^{-2}.
\]
Hence
\[
I_r(1/2,a,p)
\ge
e^{-1/2\eta^2}\int_0^{1/a} u^{r-1/2}e^{-au}\,du.
\]
With the change of variables $s=au$,
\[
\int_0^{1/a} u^{r-1/2}e^{-au}\,du
=
a^{-(r+1/2)}\int_0^1 s^{r-1/2}e^{-s}\,ds.
\]
Therefore
\[
I_r(1/2,a,p)\gtrsim a^{-(r+1/2)}.
\]
Putting the analyses above together, we have
\[
M_r(t)\asymp \sqrt p\,a^{-(r+1/2)},
\qquad r=0,1,2.
\]
Next,
\[
\frac{Z_p'(t)}{Z_p(t)}=\frac{M_1(t)}{M_0(t)}\asymp a^{-1}\asymp \frac1{p-2t}.
\]
Finally,
\[
f_p''(t)
=
\frac{Z_p''(t)}{Z_p(t)}-\left(\frac{Z_p'(t)}{Z_p(t)}\right)^2
\le
\frac{Z_p''(t)}{Z_p(t)}
=
\frac{M_2(t)}{M_0(t)}
\lesssim a^{-2}\asymp \frac1{(p-2t)^2}.
\]
The proof is complete.
$\hfill$ $\square$

\begin{lemma}\label{lem:product-bound-full}
For all $a,b\ge 0$,
\[
Z_p^{(2)}(a,b)\le Z_p(a)Z_p(b),
\]
where
\[
Z_p^{(2)}(a,b)
:=
\int_{\mathbb S^{p-1}} e^{a x_1^2+b x_2^2}\,d\sigma(x).
\]
\end{lemma}
\noindent \textbf{Proof of Lemma \ref{lem:product-bound-full}.} By Lemma 2.7 in \cite{barthe2010generalized}, the vector \((x_1^2,x_2^2,\dots,x_p^2)\) has a 
\[\mbox{Dirichlet}(1/2,\dots,1/2)\] 
distribution, and Proposition 2.10 in the same paper implies that \((x_1^2,x_2^2,\dots,x_p^2)\) is negatively associated. This means that for increasing functions \(f_1,\dots,f_p\), we have
\[
\mb{E}_{\mu_0} \prod_{i=1}^p f_i(x_i^2) \leq \prod_{i=1}^p\mb{E}_{\mu_0}  f_i(x_i^2).
\]
By setting \(f_1(x)=\exp \lb ax \rb\), \(f_2(x)=\exp \lb bx \rb\), and taking the remaining functions to be equal to \(1\),
we get the result. $\hfill$ $\square$

\subsection{Proof of Theorem \ref{watson minimax}}

Write $\mb{P}_{\kappa,\bm{\theta}}$, $\mb{E}_{\kappa,\bm{\theta}}$ for the probability and expectation under the Watson distribution with concentration parameter $\kappa$ and location $\bm{\theta}$. Given the data $\bm{X}_1,\dots,\bm{X}_n$, take $\delta$ to be a sufficiently small constant to be chosen later. Define
\begin{align}
    L_{n,\bm{\theta}} :&= \prod_{i=1}^n \frac{\exp \lb \kappa \lb \bm{\theta}^\top \bm{X}_i \rb^2 \rb}{Z_p(\kappa)}; \qquad S_{n,\bm{\theta}}:= \sum_{i=1}^n \lb \bm{\theta}^\top \bm{X}_i \rb^2; \qquad  Z_p(\kappa):= \int_{\mb{S}^{p-1}} \exp(\kappa x_1^2) \,d\mu_0(\bm{x}); \label{L-S-Z}\\
    \mu_{\kappa}:&= \frac{Z_p'(\kappa)}{Z_p(\kappa)}; \qquad \Tilde{\mu}_\kappa:=(1+\delta)\cdot\mu_{\kappa}; \qquad G_n=G_{n,\bm{\theta}}:= \la S_{n,\bm{\theta}} \leq n \, \Tilde{\mu}_\kappa \ra; \label{mu-mu-G} \\ 
    \Tilde{L}_{n, \bm{\theta}}:&= L_{n,\bm{\theta}}\,\mathbf{1}_{G_n}; \qquad
    \Delta:= \frac{p-2\kappa}{2\kappa}. \label{L-Delta}
\end{align}
Let $\pi$ be the uniform prior on $\mb{S}^{p-1}$ and set 
\[
\mc{L}_n \lb \bm{X}_1,\dots, \bm{X}_n  \rb:= \int_{\mb{S}^{p-1}} L_{n,\bm{\theta}} \, \pi(d\bm{\theta});  \qquad \Bar{L}_n \lb \bm{X}_1,\dots, \bm{X}_n  \rb:= \int_{\mb{S}^{p-1}} \Tilde{L}_{n, \bm{\theta}}\, \pi(d\bm{\theta}). 
\]
We split the proof into several steps.

\noindent \underline{\textbf{Step 1: The events $G_n^c$ are unlikely.}}
Observe that for $s>0$,
\begin{align*}
    \mb{E}_{\kappa,\bm{\theta}} \lb \exp \lb sS_{n,\bm{\theta}} \rb  \rb &=     \mb{E}_{\kappa,\bm{\theta}} \left[  \exp \lb  s\lb \bm{\theta}^\top \bm{X}_1 \rb^2 \rb  \right]^n \\
    &= \left[ \int_{\mb{S}^{p-1}} \exp \lb  s\lb \bm{\theta}^\top \bm{X}_1 \rb^2 \rb \cdot  \frac{\exp \lb \kappa \lb \bm{\theta}^\top \bm{X}_1 \rb^2 \rb}{Z_p(\kappa)} d\mu_0(\bm{x})  \right]^n \\
    &= \left[  \frac{Z_p(\kappa+s)}{Z_p(\kappa)} \right]^n.
\end{align*}
Recall from \eqref{log Zp} that $f_p(t)=\log Z_p(t)$. From the mgf above, we find that
\[
\mb{E} S_{n,\bm{\theta}}:= n\cdot f_p'(\kappa)=n\cdot \mu_\kappa; \qquad  \mbox{Var} \lb S_{n,\bm{\theta}} \rb:= n\cdot  f_p''(\kappa).
\]
Consequently, Chebyshev's inequality gives
\[
\sup_{\bm{\theta} \in \mb{S}^{p-1}} \mb{P}_{\kappa,\bm{\theta}}^{\otimes n} \lb G_n^c \rb \leq \frac{f_p''(\kappa)}{\delta^2 \,n \mu_{\kappa}^2} \lesssim \frac{1}{n} 
\]
as long as $\delta$ is fixed away from $0$, since by Lemma \ref{lem:rank-one-full}, $\mu_\kappa= f_p'(\kappa) \asymp 1/(p-2\kappa)$ and $f_p''(\kappa) \lesssim (p-2\kappa)^{-2}$.

\noindent \underline{\textbf{Step 2: Risk lower bound via truncated likelihood.}}
 Define
\[
\mathcal K_n
:=
\Bigl\{
0<\kappa<p/2:\ \delta_1(n)\ge \varepsilon
\Bigr\},
\qquad
\mathcal K_n'
:=
\Bigl\{
0<\kappa<p/2:\ \varepsilon\le \delta_1(n)\le 2\varepsilon
\Bigr\}.
\]
Obviously,
\begin{align*}
\inf_{\varphi \in \mathcal{T}_{n}}
\left\{
\mathbb{P}_{\mu_0}(\varphi=1)
+
\sup_{\mu_n \in H^W_{1,n}:\,\delta_1(n) \ge \varepsilon}
\mathbb{P}_{\mu_n}(\varphi =0)
\right\}
&=
\inf_{\varphi\in\mc T_n}
\left\{
\mb P_{\mu_0}(\varphi=1)
+
\sup_{\substack{\kappa\in\mathcal K_n\\ \bm\theta\in\mb S^{p-1}}}
\mb P_{\kappa,\bm\theta}^{\otimes n}(\varphi=0)
\right\} \\
&\ge
\inf_{\varphi\in\mc T_n}
\left\{
\mb P_{\mu_0}(\varphi=1)
+
\sup_{\bm\theta\in\mb S^{p-1}}
\mb P_{\kappa_n,\bm\theta}^{\otimes n}(\varphi=0)
\right\}
\end{align*}
for any chosen sequence \(\kappa_n\in\mathcal K_n'\), because the minimax risk over the larger model is always larger. 

Therefore, it suffices to show that
\[
\inf_{\varphi\in\mc T_n}
\left\{
\mb P_{\mu_0}(\varphi=1)
+
\sup_{\bm\theta\in\mb S^{p-1}}
\mb P_{\kappa_n,\bm\theta}^{\otimes n}(\varphi=0)
\right\}
\]
is asymptotically greater than \(1/2\) for a sequence $\kappa_n \in \mc{K}_n'$ and with $\ve$ being small enough.

Consider a sequence \(\kappa_n\in\mathcal K_n'\). For any test \(\varphi\in\mc T_n\),
\begin{align*}
\mb P_{\mu_0}(\varphi=1)
+
\sup_{\bm{\theta}\in \mb S^{p-1}}
\mb P_{\kappa_n,\bm{\theta}}^{\otimes n}(\varphi=0) 
&\ge
\mb P_{\mu_0}(\varphi=1)
+
\int_{\mb S^{p-1}}
\mb P_{\kappa_n,\bm{\theta}}^{\otimes n}(\varphi=0)\,\pi(d\bm\theta) \\
&=
\mb E_{\mu_0}\!\left[\varphi+(1-\varphi)\mc L_n\right].
\end{align*}
Since \(\mc L_n\ge \bar L_n\), we conclude that
\begin{align*}
\inf_{\varphi\in\mc T_n}
\left\{
\mb P_{\mu_0}(\varphi=1)
+
\sup_{\bm{\theta}\in \mb S^{p-1}}
\mb P_{\kappa_n,\bm{\theta}}^{\otimes n}(\varphi=0)
\right\}
&\ge
\inf_{\varphi\in\mc T_n}
\mb E_{\mu_0}\!\left[\varphi+(1-\varphi)\bar L_n\right] \\
&\ge
\mb E_{\mu_0}\!\left[\min\{1,\bar L_n\}\right],
\end{align*}
where the last inequality follows from pointwise minimization.

Next,
\[
\mb E_{\mu_0}\bar L_n
=
\int_{\mb S^{p-1}}
\mb E_{\mu_0}\!\left[L_{n,\bm\theta}\mathbf 1_{G_n}\right]\pi(d\bm\theta)
=
\int_{\mb S^{p-1}}
\mb P_{\kappa_n,\bm\theta}^{\otimes n}(G_n)\,\pi(d\bm\theta)
=
1-o(1)
\]
by Step 1. Hence
\begin{align*}
\mb E_{\mu_0}\!\left[\min\{1,\bar L_n\}\right]
&=
\mb E_{\mu_0}[\bar L_n]-\mb E_{\mu_0}\!\left[(\bar L_n-1)_+\right] \\
&\ge
1-o(1)-\sqrt{\mb E_{\mu_0}\!\left[(\bar L_n-1)^2\right]}.
\end{align*}
It follows that
\begin{align*}
\inf_{\varphi\in\mc T_n}
\left\{
\mb P_{\mu_0}(\varphi=1)
+
\sup_{\bm\theta\in\mb S^{p-1}}
\mb P_{\kappa_n,\bm\theta}^{\otimes n}(\varphi=0)
\right\}
\ge
1-o(1)
-
\sqrt{\mb E_{\mu_0}\!\left[(\bar L_n-1)^2\right]}.
\end{align*}
Consequently, it suffices to show that
\begin{align} \label{truncate go to 0}
\mb E_{\mu_0}\!\left[\bar L_n^2\right] \leq 1+\ve
\end{align}
along some sequence \(\kappa_n\) satisfying
\[
\varepsilon \le \delta_1(n)\le 2\varepsilon
\]
and $\ve$ is small enough. It is easy to see that we only need to show \eqref{truncate go to 0} when $\delta_1(n) \to 0$.

\noindent \underline{\textbf{Step 3: Second moment representation.}}
Let $\bm{U}$ and $\bm{V}$ be i.i.d. uniformly distributed on $\mb{S}^{p-1}$.
Write
\begin{align*}
    \mb{E} \lb \bar{L}_n^2 \rb &= \mb{E} \left[ \int_{\mb{S}^{p-1}} \Tilde{L}_{n, \bm{U}}\, \pi(d\bm{U}) \cdot \int_{\mb{S}^{p-1}} \Tilde{L}_{n, \bm{V}}\, \pi(d\bm{V})  \right] \\
    &= \mb{E} \left[   \mb{E} \lb  \Tilde{L}_{n, \bm{U}}\, \Tilde{L}_{n, \bm{V}} \Big| \bm{U}, \bm{V} \rb \right] \\
    &= \mb{E} \left[  \mc{R} \lb \bm{U}, \bm{V}   \rb \right]
\end{align*}
where
\[
\mc{R} \lb \bm{U}, \bm{V} \rb:=    \mb{E} \lb  \Tilde{L}_{n, \bm{U}}\, \Tilde{L}_{n, \bm{V}} \Big| \bm{U}, \bm{V} \rb.
\]
Define the good set
\[
\mc{G}_n:= \la \lb \bm{U}, \bm{V} \rb: |\bm{U}^\top \bm{V}| \leq \frac{\Delta}{2} \ra
\]
and define its complement to be the bad set. We will bound the expectation of $\mc{R}$ over the good set and the bad set in the next two steps.

\noindent \underline{\textbf{Step 4: Bounding the contribution over the good set $\mc{G}_n$.}}
We first obtain a pointwise bound on $\mc{R}(\bm{U},\bm{V})$ when $(\bm{U},\bm{V})\in \mc{G}_n$.

Write
\[
\rho:=\bm{U}^\top \bm{V}, \qquad
a_+:=\kappa(1+\rho), \qquad a_-:=\kappa(1-\rho),
\]
Since truncation only decreases the likelihood ratio, with $\tilde{L}_{n,\bm{\theta}}$ as in \eqref{L-Delta}, we have
\[
\mc{R}(\bm{U},\bm{V})
=
\mb{E}\Bigl[\Tilde{L}_{n,\bm{U}}\Tilde{L}_{n,\bm{V}}\,\big|\,\bm{U},\bm{V}\Bigr]
\le
\mb{E}\Bigl[L_{n,\bm{U}}L_{n,\bm{V}}\,\big|\,\bm{U},\bm{V}\Bigr].
\]
Using rotational invariance under $\mu_0$ and the definition of $Z_p^{(2)}$ as in Lemma \ref{lem:product-bound-full}, this gives
\[
\mc{R}(\bm{U},\bm{V})
\le
\left[
\frac{Z_p^{(2)}(a_+,a_-)}{Z_p(\kappa)^2}
\right]^n.
\]
Thanks to Lemma~\ref{lem:product-bound-full},
$
Z_p^{(2)}(a_+,a_-)\le Z_p(a_+)Z_p(a_-),
$
and hence
\[
\mc{R}(\bm{U},\bm{V})
\le
\exp\bigl(nF(\rho)\bigr),
\qquad
F(\rho):=f_p(a_+)+f_p(a_-)-2f_p(\kappa).
\]
Clearly,
\[
F(0)=0,
\qquad
F'(0)=0,
\]
and
\[
F''(t)
=
\kappa^2 f_p''\bigl(\kappa(1+t)\bigr)
+
\kappa^2 f_p''\bigl(\kappa(1-t)\bigr).
\]
On $\mc{G}_n$, one has $|\rho|\le \Delta/2$.
Thus for every $|t|\le |\rho|$,
\[
p-2\kappa(1\pm t)
=
2\kappa(\Delta\mp t)
\ge
2\kappa(\Delta-|t|)
\ge
2\kappa(\Delta-|\rho|)
\ge
\kappa\Delta.
\]
Note that $\kappa\Delta= p/2-\kappa \gg \sqrt{p}$ because
\[
\frac{p}{2} -\kappa = \frac{p\sqrt{n}}{2(\sqrt{n}+\sqrt{\delta_1(n)\, p})} \qquad \text{and} \qquad \sqrt{pn} \gg \sqrt{p}+\sqrt{n}.
\]
Therefore Lemma~\ref{lem:rank-one-full} applies uniformly to
\(\kappa(1+t)\) and \(\kappa(1-t)\), and yields
\[
f_p''\bigl(\kappa(1+t)\bigr)\lesssim (\kappa\Delta)^{-2},
\qquad
f_p''\bigl(\kappa(1-t)\bigr)\lesssim (\kappa\Delta)^{-2},
\qquad |t|\le |\rho|.
\]
Hence
\[
\sup_{|t|\le |\rho|}F''(t)\lesssim \frac{1}{\Delta^2}.
\]
Consequently, Taylor's expansion around \(0\) gives
\[
F(\rho)
\le
\frac12 \rho^2 \sup_{|t|\le |\rho|}F''(t)
\lesssim
\frac{\rho^2}{\Delta^2}.
\]
Consequently, on $\mc{G}_n$,
\[
\mc{R}(\bm{U},\bm{V})
\le
\exp\!\left(
C n \frac{\rho^2}{\Delta^2}
\right)
\]
for some universal constant \(C>0\).

Note that we may rewrite the exponent in the last display as
\[
n\frac{\rho^2}{\Delta^2}
=
\delta_1(n)\,\zeta_p
\]
where 
\begin{align} \label{zeta}
    \zeta_p:= p\rho^2 \stackrel{d}{\to} \chi_1^2.
\end{align}
In other words,
\[
\mc{R}(\bm{U},\bm{V})\,\mathbf 1_{\mc{G}_n}
\le
\exp \bigl(C\,\delta_1(n)\,\zeta_p\bigr).
\]
Taking expectations and using the exact mgf of $\zeta_p$ (which is the scaled inner product under uniformity), we obtain for all sufficiently small \(\delta_1(n)\),
\[
\mb{E}\Bigl[\mc{R}(\bm{U},\bm{V})\,\mathbf 1_{\mc{G}_n}\Bigr]
\le
\mb{E}\Bigl[\exp \bigl(C\,\delta_1(n)\,\zeta_p\bigr)\Bigr]
\le
\bigl(1-2C\,\delta_1(n)\bigr)^{-1/2}.
\]
The bound above is well known and can be checked easily using the exact formula for the mgf of $\rho$. In particular,
\[
\mb{E}\Bigl[\mc{R}(\bm{U},\bm{V})\,\mathbf 1_{\mc{G}_n}\Bigr]
\leq 
1+o(1)
\qquad\text{whenever}\qquad
\delta_1(n)\to 0.
\]

\noindent \underline{\textbf{Step 5: Bounding the contribution over the bad set $\mc{G}_n^c$.}}
If $\Delta\ge 2$, then $\mc{G}_n^c=\varnothing$ since $|\rho|\le 1$, and there is nothing to prove. Thus we may assume throughout this step that
$
\Delta<2.
$

Fix $(\bm{U},\bm{V})$ with $|\rho|>\Delta/2$. We will consider two cases.

\smallskip
\noindent
\underline{\emph{Case 1: $\rho\ge 0$.}}
Define
\[
\bm v_+:=\frac{\bm U+\bm V}{\sqrt{2+2\rho}},
\qquad
\bm v_-:=\frac{\bm U-\bm V}{\sqrt{2-2\rho}}.
\]
Then $(\bm v_+,\bm v_-)$ is an orthonormal pair in $\mathbb R^p$. Let
\[
T_+:=\sum_{i=1}^n (\bm v_+^\top \bm X_i)^2,
\qquad
T_-:=\sum_{i=1}^n (\bm v_-^\top \bm X_i)^2.
\]
Observe that since
\[
\bm U=\sqrt{\frac{1+\rho}{2}}\,\bm v_+ + \sqrt{\frac{1-\rho}{2}}\,\bm v_-,
\qquad
\bm V=\sqrt{\frac{1+\rho}{2}}\,\bm v_+ - \sqrt{\frac{1-\rho}{2}}\,\bm v_-,
\]
we have
\[
S_{n,\bm U}+S_{n,\bm V}=(1+\rho)T_+ + (1-\rho)T_-,
\]
Hence on $G_{n,\bm U}\cap G_{n,\bm V}$,
\[
T_+\le \frac{2n\tilde{\mu}_\kappa}{1+\rho}.
\]
Therefore, for every $\lambda\ge 0$,
\[
\mathbf 1_{G_{n,\bm U}\cap G_{n,\bm V}}
\le
\exp\!\left\{
\lambda\Bigl(\frac{2n\tilde{\mu}_\kappa}{1+\rho}-T_+\Bigr)
\right\}.
\]
Write 
\begin{align*}
    \mc{R}(\bm U,\bm V) &= \mb{E} \lb L_{n,\bm{U}} \cdot L_{n,\bm{V}} \cdot \mathbf{1}_{G_{n,\bm{U}} \cap G_{n,\bm{V}}} \Big| \bm{U}, \bm{V} \rb \\
    &\leq \exp\!\left\{
\lambda\cdot\frac{2n\tilde{\mu}_\kappa}{1+\rho}
\right\} \cdot  \mb{E} \left[  L_{n,\bm{U}} \cdot L_{n,\bm{V}} \cdot \exp \lb -\lambda T_{+} \rb \Big| \bm{U}, \bm{V}    \right] \\
&= \exp\!\left\{
\lambda\cdot\frac{2n\tilde{\mu}_\kappa}{1+\rho}
\right\} \cdot  \mb{E} \left[  \frac{\exp \lb \kappa S_{n,\bm{U}} + \kappa S_{n, \bm{V}} - \lambda T_{+} \rb}{Z_p(\kappa)^{2n}} \Big| \bm{U}, \bm{V}   \right] \\
&= \exp\!\left\{
\lambda\cdot\frac{2n\tilde{\mu}_\kappa}{1+\rho}
\right\} \cdot  \mb{E} \left[  \frac{\exp \Big( (a_{+}-\lambda)T_{+} + a_{-}T_{-} \Big)}{Z_p(\kappa)^{2n}} \Big| \bm{U}, \bm{V}   \right]\\
&= \exp\!\left\{
\lambda\cdot\frac{2n\tilde{\mu}_\kappa}{1+\rho}
\right\} \cdot  \mb{E} \left[  \frac{\exp \Big( (a_{+}-\lambda)(\bm v_+^\top \bm X_1)^2 + a_{-}(\bm v_-^\top \bm X_1)^2 \Big)}{Z_p(\kappa)^{2}} \Big| \bm{U}, \bm{V}   \right]^n 
\end{align*}
Conditional on $\bm{U}$ and $\bm{V}$, the two vectors $\bm{v}_+$ and $\bm{v}_-$ are orthogonal, so by using the rotational invariance of $\mu_0$, we deduce that
\[
 \mb{E} \left[  \frac{\exp \Big( (a_{+}-\lambda)(\bm v_+^\top \bm X_1)^2 + a_{-}(\bm v_-^\top \bm X_1)^2 \Big)}{Z_p(\kappa)^2} \Big| \bm{U}, \bm{V}   \right] =  \mb{E} \left[  \frac{\exp \Big( (a_{+}-\lambda)x_1^2 + a_{-}x_2^2 \Big)}{Z_p(\kappa)^2}    \right] = 
\left[
\frac{Z_p^{(2)}(a_+-\lambda,a_-)}{Z_p(\kappa)^2}
\right].
\]
Here $(x_1,x_2)$ are the first two coordinates of a point under $\mu_0$. This yields
\begin{equation}\label{eq:exp-trick-pos}
\mc{R}(\bm U,\bm V)
\le
\exp\!\left(\frac{2\lambda n\tilde{\mu}_\kappa}{1+\rho}\right)
\left[
\frac{Z_p^{(2)}(a_+-\lambda,a_-)}{Z_p(\kappa)^2}
\right]^n.
\end{equation}
Note that \eqref{eq:exp-trick-pos} is precisely our gain via truncation: instead of having the raw $a_+$, which can be as large as $2\kappa$ or $p$ and prevent us from using the estimate from Lemma \ref{lem:rank-one-full}, we now have $a_{+}-\lambda$, which is smaller. We will tune $\lambda$ in the following analysis.

\smallskip
\noindent
\underline{\emph{Case 2: $\rho<0$.}}
Arguing similarly as in the previous case, we have 
\begin{equation}\label{eq:exp-trick-neg}
\mc{R}(\bm U,\bm V)
\le
\exp\!\left(\frac{2\lambda n\tilde{\mu}_\kappa}{1-\rho}\right)
\left[
\frac{Z_p^{(2)}(a_+,a_--\lambda)}{Z_p(\kappa)^2}
\right]^n.
\end{equation}

\smallskip
\noindent
We now unify the two cases. Define
\[
a_{\max}:=\kappa(1+|\rho|),
\qquad
a_{\min}:=\kappa(1-|\rho|).
\]
Then \eqref{eq:exp-trick-pos} and \eqref{eq:exp-trick-neg} imply that for every $\lambda\ge 0$,
\begin{equation}\label{eq:exp-trick-final}
\mc{R}(\bm U,\bm V)
\le
\exp\!\left(\frac{2\lambda n\tilde{\mu}_\kappa}{1+|\rho|}\right)
\left[
\frac{Z_p^{(2)}(a_{\max}-\lambda,a_{\min})}{Z_p(\kappa)^2}
\right]^n.
\end{equation}
Now choose
\[
t_0:=\kappa\left(1+\frac{\Delta}{2}\right),
\qquad
\lambda:=a_{\max}-t_0.
\]
Since $|\rho|>\Delta/2$, we have
\[
\lambda=\kappa\left(|\rho|-\frac{\Delta}{2}\right).
\]
Write
$
r:=\frac{|\rho|}{\Delta}.
$
Then on $\mc G_n^c$, we have
\[
r>\frac12,
\qquad
\lambda=\kappa\Delta\left(r-\frac12\right)
=
\left(\frac p2-\kappa\right)\left(r-\frac12\right).
\]

Note that  $a_{\max}-\lambda = t_0$, so Lemma~\ref{lem:product-bound-full} gives
\[
Z_p^{(2)}(a_{\max}-\lambda,a_{\min})
\le
Z_p(t_0)\,Z_p(a_{\min}).
\]
Also,
\[
\frac p2-t_0=\frac{\kappa\Delta}{2},
\qquad
\frac p2-a_{\min}=\kappa\Delta+\kappa|\rho|=\kappa\Delta(1+r)= \lb \frac{p}{2} -\kappa \rb(1+r).
\]
Therefore, Lemma~\ref{lem:rank-one-full} applies to \(t_0\), \(a_{\min}\), and \(\kappa\) for all sufficiently large \(p\) since $\kappa \Delta = p/2-\kappa \gg \sqrt{p}$. Thus
\[
Z_p(t_0)\le C\sqrt{\frac{p}{\kappa\Delta}};
\qquad
Z_p(a_{\min})\le C\sqrt{\frac{p}{\kappa\Delta(1+r)}};
\qquad
Z_p(\kappa)\ge c\sqrt{\frac{p}{\kappa\Delta}}.
\]
Hence
\[
\frac{Z_p^{(2)}(a_{\max}-\lambda,a_{\min})}{Z_p(\kappa)^2}
\le
C\frac{1}{\sqrt{1+r}}.
\]
Also, by Lemma~\ref{lem:rank-one-full},
\[
\tilde{\mu}_\kappa=(1+\delta)\mu_\kappa\le \frac{C}{\kappa\Delta}.
\]
Therefore
\[
\frac{2\lambda \tilde{\mu}_\kappa}{1+|\rho|} \leq \kappa \lb 2|\rho| - \Delta \rb \cdot \frac{C}{\kappa \Delta}
\le Cr.
\]
Using \eqref{eq:exp-trick-final}, we conclude that on \(\mc G_n^c\),
\[
\mc{R}(\bm U,\bm V)
\le
e^{Cnr}(1+r)^{-n/2}
\le
e^{Cnr}.
\]
Next, with $\zeta_p$ as in \eqref{zeta}, observe that
\[
r^2=\frac{\rho^2}{\Delta^2}
=\frac{\zeta_p}{p\Delta^2}
=\frac{\delta_1(n)\,\zeta_p}{n},
\]
and therefore
\[
e^{Cnr}
=
\exp\!\left(C\sqrt{n\,\delta_1(n)\,\zeta_p}\right).
\]
Moreover, on the bad set \(\mc G_n^c=\{r>1/2\}\),
\[
\zeta_p
=
p\rho^2
>
\frac{p\Delta^2}{4}
=
\frac{n}{4\delta_1(n)}.
\]
Thus
\[
\mc R(\bm U,\bm V)\mathbf 1_{\mc G_n^c}
\le
\exp\!\left(C\sqrt{n\,\delta_1(n)\,\zeta_p}\right)
\mathbf 1_{\{\zeta_p\ge n/(4\delta_1(n))\}}.
\]

Fix two constants $\alpha,\beta$ such that \(0<\alpha<\beta<1/2\). By Young's inequality,
\[
C\sqrt{n\,\delta_1(n)\,\zeta_p}
\le
\alpha \zeta_p + \frac{C^2 n\delta_1(n)}{4\alpha}.
\]
Hence
\[
\mc R(\bm U,\bm V)\mathbf 1_{\mc G_n^c}
\le
\exp\!\left(\frac{C^2 n\delta_1(n)}{4\alpha}\right)
e^{\alpha\zeta_p}\,
\mathbf 1_{\{\zeta_p\ge n/(4\delta_1(n))\}}.
\]
Taking expectations, we get
\[
\mb E\!\left[\mc R(\bm U,\bm V)\mathbf 1_{\mc G_n^c}\right]
\le
\exp\!\left(\frac{C^2 n\delta_1(n)}{4\alpha}\right)
\mb E\!\left[e^{\alpha\zeta_p}\mathbf 1_{\{\zeta_p\ge n/(4\delta_1(n))\}}\right].
\]
On the event \(\{\zeta_p\ge x\}\), we have
\[
e^{\alpha\zeta_p}\le e^{-(\beta-\alpha)x}e^{\beta\zeta_p}.
\]
Applying this with \(x=n/(4\delta_1(n))\), we obtain
\[
\mb E\!\left[e^{\alpha\zeta_p}\mathbf 1_{\{\zeta_p\ge n/(4\delta_1(n))\}}\right]
\le
\exp\!\left(-\frac{(\beta-\alpha)n}{4\delta_1(n)}\right)\mb E[e^{\beta\zeta_p}].
\]
Moreover, it is easy to see from \eqref{zeta} and the exact formula for the mgf of $\rho$ that if $\beta<1/10$,
\[
\mb E[e^{\beta\zeta_p}]
\le
(1-2\beta)^{-1/2}
\]
uniformly in $p \geq 3$. Therefore
\[
\mb E\!\left[\mc R(\bm U,\bm V)\mathbf 1_{\mc G_n^c}\right]
\le
(1-2\beta)^{-1/2}
\exp\!\left(
\frac{C^2 n\delta_1(n)}{4\alpha}
-
\frac{(\beta-\alpha)n}{4\delta_1(n)}
\right)
=o(1),
\]
since \(\delta_1(n)\to 0\). The proof is complete. $\hfill$ $\square$


\section{Proof of Theorem \ref{watson local distribution}}

Due to the technical nature of the proof, we make the convention that the constants in the estimates in the proof below might change from line to line, despite using the same notation. Without stating otherwise, $C_a$ denotes a constant depending on $a$, $C_{a,\tau}$ denotes a constant depending on $a, \tau$ and so on. The notation $\lesssim$ indicates an upper bound that holds up to a universal constant.

Let $\la \bm{U}_i; 1\leq i \leq n \ra$ be i.i.d. uniformly distributed on $\mb{S}^{p-2}$ and let $\la Y_i; 1\leq i \leq n \ra$ be i.i.d. with common density proportional to
\[
\exp(\kappa x^2)\, (1-x^2)^{(p-3)/2}\, \mathbf{1}_{[-1,1]}(x). 
\]
By rotational invariance, we have
\begin{align*}
    \la  \bm{X}_i^\top \bm{X_j} \ra_{1 \leq i<j \leq n} \stackrel{d}{=} \la Y_iY_j + \sqrt{(1-Y_i^2)(1-Y_j^2)} \cdot \bm{U}_i^\top \bm{U}_j  \ra
\end{align*}
where the inner products on the left-hand side are under the Watson distribution on $\mb{S}^{p-1}$ with concentration $\kappa$.

Put $N_n := \binom{n}{2}$ and define
\begin{align}
    \beta_{ij} :&= Y_i Y_j;  \qquad
    \alpha_{ij} := \sqrt{(1-Y_i^2)(1-Y_j^2)}; \label{beta ij}  \\
    R_{ij} :&=   \bm{U}_i^\top \bm{U}_j; \nonumber   \\
    D_{ij}(t):&= \mathbf 1_{\{\sqrt p(\beta_{ij}+\alpha_{ij}R_{ij})\le t\}}
-
\mathbf 1_{\{\sqrt p R_{ij}\le t\}}; \nonumber \\
\pi_{ij}(t):&=
\mathbb P\!\left(\sqrt p\,(\beta_{ij}+\alpha_{ij}R_{12})\le t
\,\middle|\,Y_i,Y_j\right). \nonumber
\end{align}
Decompose 
\[
\sqrt{\frac{2}{n(n-1)}}  \left[ \sum_{1\leq i<j \leq n} \mathbf{1}_{\la \sqrt{p}\,\bm{X}_i^\top \bm{X}_j  \leq t \ra} - \mb{P}_{\mu_n} \lb \sqrt{p}\,\bm{X}_1^\top \bm{X}_2 \leq t \rb  \right] \stackrel{d}{=} \mc{Z}_n(t) + \mc{A}_n(t) + \mc{B}_n(t) + g_n(t)
\]
where
\begin{align}
     \mc{Z}_n(t) &:=  N_n^{-1/2} \cdot \left[ \sum_{1\leq i<j \leq n} \mathbf{1}_{\la \sqrt{p}\,R_{ij}  \leq t \ra} - \mb{P} \lb \sqrt{p}\, R_{12} \leq t \rb  \right]; \nonumber \\
     \mc{A}_n(t) &:= N_n^{-1/2} \cdot \sum_{1\le i<j\le n}
\Bigl(
D_{ij}(t)-\mathbb E\bigl[D_{ij}(t)\mid Y_i,Y_j\bigr] 
\Bigr); \label{An(t)} \\
\mc{B}_n(t) &:= N_n^{-1/2} \cdot 
\sum_{1\le i<j\le n}
\Bigl(
\pi_{ij}(t)-\mathbb E\pi_{12}(t)
\Bigr); \label{Bn(t)} \\
g_n(t):&= N_n^{1/2}\lb \mb{P}_{\mu_n} \lb \sqrt{p}\,\bm{X}_1^\top \bm{X}_2 \leq t \rb - \mb{P}_{\mu_0} \lb \sqrt{p}\,\bm{X}_1^\top \bm{X}_2 \leq t \rb  \rb \nonumber.
\end{align}
Here is a short description of the above terms. The term $\mc{Z}_n$ behaves like $T_n$ under the null, $\mc{A}_n$ and $\mc{B}_n$ are  random perturbations and $g_n(t)$ is a smooth deterministic shift. The term $\mc{Z}_n$ is already handled by Theorem \ref{null-dist}, and we know that it converges weakly to a Brownian bridge. By Proposition \ref{watson d asymptotic}, we also know that $g_n(t) $ converges uniformly to a smooth function. 

To deduce the result, we thus only need to verify that
\begin{align} \label{mcA}
    \sup_{t \in \mb{R}} |\mc{A}_n(t)| = o_{\mb{P}}(1)
\end{align}
and 
\begin{align} \label{mcB}
    \sup_{t \in \mb{R}} |\mc{B}_n(t)| = o_{\mb{P}}(1).
\end{align}

\begin{remark}
The suprema in \eqref{mcA} and \eqref{mcB} are of a different nature. In \eqref{mcA}, we are dealing with the supremum of a degenerate U-process, although it may not be immediate that the process is indeed of this type. Empirical process concentration inequalities can be used to handle \eqref{mcA}, provided that its variance profile converges to zero but does not decay exponentially fast (the condition $\log p = o(n^{1/4})$ is used to rule out this possibility). Lemmas \ref{upper bound variance An} and \ref{lower bound variance An} are devoted to studying this variance profile. 

In \eqref{mcB}, by contrast, we have a non-degenerate U-process, which is also a smooth function for every $n$. Empirical process concentration inequalities are no longer effective, since one would need a careful analysis of the magnitudes of both the degenerate and the non-degenerate parts. We bypass this difficulty by exploiting the smoothness of the process and reducing the problem of bounding the supremum to that of bounding the variance of its derivatives. This is done via the Brascamp--Lieb inequality (see Lemma \ref{Brascamp-Lieb} below).
\end{remark}

\subsection{Technical estimates and tools}

Note that $Y_1$ has the same distribution as the random variable $T$ in Lemma \ref{moment asymptotic}. For a fixed integer $m \geq 1$, by using a similar computation as in Lemma \ref{moment asymptotic} we have, for all large $p$,
\begin{equation}\label{eq:simple-weighted-delta}
\mathbb E\!\left[\frac{Y_1^{2k}}{(1-Y_1^2)^m}\right]
\le \frac{C_m}{\delta_p}
\end{equation}
for $k \in \la 0,1,2,3\ra$ and some constant $C_m$ depending only on $m$.

We skip the proof of \eqref{eq:simple-weighted-delta} since the computation is almost identical to what is done in the proof of Lemma \ref{moment asymptotic}. We will also need a variant of the classical Brascamp-Lieb inequality to bound the fluctuation of $\mc{B}_n(t)$. The lemma below is well-known and can likely be used without proof. However, we were not able to locate an exact reference that covers it as a corollary, so we decided to give a proof below.
\begin{lemma}[Brascamp-Lieb's inequality] \label{Brascamp-Lieb}
    Denote by $\nu$ the probability measure on $(-1,1)$ with density proportional to 
    \[
    \exp(\kappa x^2)\, (1-x^2)^{(p-3)/2} \mathbf{1}_{(-1,1)}(x). 
    \]
   Assume that $(p-3)/2-\kappa>0$. If $Y_1,\cdots,Y_n$ are drawn i.i.d. from $\nu$, then
   \[
   \mbox{Var} \lb f \lb Y_1,\dots,Y_n \rb \rb \leq \frac{1}{p-3-2\kappa} \cdot \mb{E} \lb \|\nabla f  \lb Y_1,\dots,Y_n \rb\|^2   \rb.
   \]
   Here $f$ is a locally Lipschitz function such that $f \in L^2 \lb \nu \rb$. 
\end{lemma}

\noindent \textbf{Proof of Lemma \ref{Brascamp-Lieb}.}
It suffices to prove the result for $n=1$ due to the tensorization property of variance; see Chapter 2.1 in \cite{van2014probability}. Put $V(x)= -\kappa x^2 - \frac{p-3}{2}\, \log(1-x^2)$ and note that $\nu \propto \exp(-V(x))dx$. A direct computation yields that 
\[
V''(x) = -2\kappa + \frac{p-3}{2}\cdot \frac{2(1+x^2)}{(1-x^2)^2} \geq p-3-2\kappa. 
\]
Put $\phi_q(x)= \exp \lb V(x)/q \rb$. It is easy to see that $\phi_q$ is a $C^2$, strongly convex function for all $q>0$. Apply Theorem 8 in \cite{nguyen2014} to the law $\phi_q^{-q}$ (which is $\nu$) on the open, convex domain $(-1,1)$, we conclude that
\begin{align*}
\mbox{Var}_{\nu} \lb f \rb \leq \frac{1}{q-1} \int_{-1}^1 \frac{f'(x)^2}{\phi_q''(x)} d\nu &= \frac{1}{q-1} \int_{-1}^1 \frac{f'(x)^2\cdot e^{V/q}}{\frac{1}{q} \left[ e^{V/q} \cdot \frac{(V')^2}{q}  + e^{V/q}\cdot V^{''}(x) \right]   } d\nu \\
&= \frac{q}{q-1} \int_{-1}^1 \frac{f'(x)^2 \cdot e^{V/q}}{ e^{V/q} \cdot \frac{(V')^2}{q}  + e^{V/q}\cdot V^{''}(x)   } d\nu \\
&=  \frac{q}{q-1} \int_{-1}^1 \frac{f'(x)^2}{  \frac{(V')^2}{q}  +  V^{''}(x)   } d\nu \\
&\leq \frac{q}{q-1} \int_{-1}^1 \frac{f'(x)^2}{  V^{''}(x)   } d\nu \leq  \frac{q}{(q-1)(p-3-2\kappa)}  \int_{-1}^1 f'(x)^2 d\nu.
\end{align*}
The proof is completed by sending $q \to \infty$. $\hfill$ $\square$

The next result is a variance bound on $\mc{A}_n(t)$ that is needed for our analysis

\begin{lemma} \label{upper bound variance An}
    Recall $\mc{A}_n(t)$ in \eqref{An(t)}. Then for every fixed $a>0$, we have
    \begin{align}
    \sup_{t\in [-a,a]} \mbox{Var} \lb \mc{A}_n(t) \rb &\leq C_{a,\tau} \lb \frac{1}{\sqrt{p}} + \frac{1}{\sqrt{n}} \rb \label{variance upper bound mcA1}. 
\end{align}
\end{lemma}

\noindent \textbf{Proof of Lemma \ref{upper bound variance An}.}
Write
\[
W_{ij}(t)
:=
D_{ij}(t)-\mathbb E\!\left[D_{ij}(t)\mid Y_1,\dots,Y_n\right].
\]
Then
\begin{align} \label{degenerate}
\mathcal A_n(t)=\frac{1}{\sqrt{N_n}}\sum_{1\le i<j\le n}W_{ij}(t).
\end{align}
It is easy to see that conditional on $Y_1,\dots,Y_n$, the array
$\{W_{ij}(t):1\le i<j\le n\}$ is degenerate in the $\bm{U}$-variables.
Indeed, fix $i<j$. Conditional on $\bm{U}_i$ and $Y_1,\dots,Y_n$, the random variable
$R_{ij}=\bm{U}_i^\top \bm{U}_j$ has the same distribution as $R_{12}$, by rotational invariance of the uniform law on $\mathbb S^{p-2}$. Therefore
\[
\mathbb E\!\left[D_{ij}(t)\mid \bm{U}_i,Y_1,\dots,Y_n\right]
=
\mathbb E\!\left[D_{ij}(t)\mid Y_1,\dots,Y_n\right],
\]
and hence
\[
\mathbb E\!\left[W_{ij}(t)\mid \bm{U}_i,Y_1,\dots,Y_n\right]=0.
\]
Similarly,
\[
\mathbb E\!\left[W_{ij}(t)\mid \bm{U}_j,Y_1,\dots,Y_n\right]=0.
\]
Hence
\[
\operatorname{Var}\!\left(\mathcal A_n(t)\mid Y_1,\dots,Y_n\right)
=
\frac{1}{N_n}\sum_{1\le i<j\le n}
\operatorname{Var}\!\left(W_{ij}(t)\mid Y_1,\dots,Y_n\right).
\]
Since $\operatorname{Var}(W_{ij}(t)\mid Y_i,Y_j)\le \mathbb E(D_{ij}(t)^2\mid Y_i, Y_j)$, we obtain
\[
\operatorname{Var}\!\left(\mathcal A_n(t)\mid Y_1,\dots,Y_n\right)
\le
\frac{1}{N_n}\sum_{1\le i<j\le n}
\mathbb E\!\left[D_{ij}(t)^2\mid Y_1,\dots,Y_n\right].
\]
Taking expectation, using $\mb{E} \left[ \mc{A}_n(t)| Y_1,\dots,Y_n \right]=0$ and exchangeability , we have
\begin{equation}
\label{eq:An-second-moment}
\mathbb E\bigl[\mathcal A_n(t)^2\bigr]
\le
\mathbb E\!\left[D_{12}(t)^2\right]
=
\mathbb P\!\left(D_{12}(t)\neq 0\right).
\end{equation}
We now bound $\mathbb P(D_{12}(t)\neq 0)$.
Let
\[
G:=\{|Y_1|\le 1/2,\ |Y_2|\le 1/2\}.
\]
On $G$, we have $\alpha_{12}\ge \frac34$. Moreover, $D_{12}(t)\neq 0$ can occur only if $R_{12}$ lies between the two thresholds $t/\sqrt{p}$ and $\frac{t/\sqrt{p}-\beta_{12}}{\alpha_{12}}$.
Hence, on $G$, using the independence between $R_{12}$ and $(\alpha_{12}, \beta_{12})$, we have 
\[
\mathbb P\!\left(D_{12}(t)\neq 0\mid Y_1,Y_2\right)
\le
\sup_{u\in[-1,1]} f_p(u)\,
\left|
\frac{t/\sqrt{p}-\beta_{12}}{\alpha_{12}}- t/\sqrt{p}
\right|,
\]
where $f_p$ denotes the density of $R_{12}$.
Since
\[
f_p(u)=c_p(1-u^2)^{(p-4)/2}\mathbf 1_{\{|u|\le 1\}},
\qquad
c_p=\frac{\Gamma((p-1)/2)}{\sqrt\pi\,\Gamma((p-2)/2)},
\]
the gamma-ratio asymptotics imply
\[
\sup_{u\in[-1,1]} f_p(u)\le C\sqrt p.
\]
for a universal constant $C>0$.

Next, on $G$,
\[
\left|
\frac{t/\sqrt{p}-\beta_{12}}{\alpha_{12}}- t/\sqrt{p}
\right|
=
\left|
\frac{t(1-\alpha_{12})-\sqrt{p}\,\beta_{12}}{\sqrt{p}\,\alpha_{12}}
\right|
\le
\frac{4}{3}\lb|\beta_{12}|+ \frac{a}{\sqrt{p}}(1-\alpha_{12})\rb.
\]
Furthermore,
\[
1-\alpha_{12}
=
\frac{1-\alpha_{12}^2}{1+\alpha_{12}}
=
\frac{Y_1^2+Y_2^2-Y_1^2Y_2^2}{1+\alpha_{12}}
\le
Y_1^2+Y_2^2.
\]
Therefore, on $G$,
\[
\left|
\frac{t/\sqrt{p}-\beta_{12}}{\alpha_{12}}- t/\sqrt{p}
\right|
\le
C_a\Bigl(|Y_1Y_2|+p^{-1/2}(Y_1^2+Y_2^2)\Bigr).
\]
for some constant $C_a$ depending only on $a$. Combining the last two displays yields
\[
\mathbb P\!\left(D_{12}(t)\neq 0\mid Y_1,Y_2\right)\mathbf 1_G
\le
C_a\Bigl(\sqrt p\,|Y_1Y_2|+Y_1^2+Y_2^2\Bigr)\mathbf 1_G.
\]
On $G^c$, we have the trivial bound 
\[
\mathbb P\!\left(D_{12}(t)\neq 0\mid Y_1,Y_2\right)\le 1.
\]
Hence
\[
\mathbb P\!\left(D_{12}(t)\neq 0\right)
\le
C_a\sqrt p\,\mathbb E|Y_1Y_2|
+
C_a\,\mathbb E(Y_1^2)
+
\mathbb P(G^c).
\]
Now, bv using Cauchy--Schwarz, $\mathbb E|Y_1Y_2|
\le
\sqrt{\mathbb E(Y_1^2)\,\mathbb E(Y_2^2)}
=
\mathbb E(Y_1^2)
$, so Lemma \ref{moment asymptotic} (with $T$ replaced by $Y_1$) gives
\[
\sqrt p\,\mathbb E|Y_1Y_2|
\le
C\frac{1+r_p}{\sqrt p}.
\]
where $r_p$ is defined as in \eqref{r-p}. Also,
\[
\mathbb E(Y_1^2)\le C\frac{1+r_p}{p}.
\]
Finally, by the union bound, Markov's inequality and Lemma \ref{moment asymptotic}
\[
\mathbb P(G^c)
\le
2\mathbb P(|Y_1|>1/2)
\le
32\,\mathbb E(Y_1^4)
\le
C\frac{(1+r_p)^2}{p^2}.
\]
Therefore
\[
\mathbb P\!\left(D_{12}(t)\neq 0\right)
\le
C_a\left(
\frac{1+r_p}{\sqrt p}
+
\frac{1+r_p}{p}
+
\frac{(1+r_p)^2}{p^2}
\right)\to 0
\]
since $r_p \asymp \sqrt{p/n}$. Combine the above with \eqref{eq:An-second-moment}, we get \eqref{variance upper bound mcA1} since $r_p \asymp \sqrt{p/n}$. $\hfill$ $\square$

In addition to the upper bound, we also need the following variance lower bound.

\begin{lemma} \label{lower bound variance An}
      Recall $\mc{A}_n(t)$ in \eqref{An(t)}. Then for every fixed $a>0$, we have
    \begin{align}
    \mbox{Var} \lb \mc{A}_n(a) \rb &\geq C_{a}\,p^{-1/2} \label{variance lower bound mcA1}.
\end{align}
\end{lemma}

\noindent \textbf{Proof of Lemma \ref{lower bound variance An}.}
As in the proof of Lemma \ref{upper bound variance An}, recall $f_p$ being the density of $R_{12}$. Put $\tilde{F}_p$ to be its density.  Thus
\[
\tilde{F}_p(x)=\mb{P}(R_{12}\le x),
\qquad
f_p(x)=c_p(1-x^2)^{(p-4)/2}\mathbf 1_{\{|x|\le 1\}},
\]
where
\[
c_p=\frac{\Gamma((p-1)/2)}{\sqrt{\pi}\,\Gamma((p-2)/2)}.
\]
Define
\[
\Delta_{12}(a):=\frac{a/\sqrt p-\beta_{12}}{\alpha_{12}}-\frac{a}{\sqrt p}
=
-\frac{\beta_{12}}{\alpha_{12}}
+
\frac{a}{\sqrt p}\Bigl(\frac1{\alpha_{12}}-1\Bigr).
\]
Then $D_{12}(a)=\mathbf 1_{\{R_{12}\le a/\sqrt p+\Delta_{12}(a)\}}
-\mathbf 1_{\{R_{12}\le a/\sqrt p\}}$. Hence, conditional on \(Y_1,Y_2\),
\[
\Var\!\left(D_{12}(a)\mid Y_1,Y_2\right)
=
q_{12}(a)\bigl(1-q_{12}(a)\bigr),
\]
where
\[
q_{12}(a):=\bigl|F_p(a/\sqrt p+\Delta_{12}(a))-F_p(a/\sqrt p)\bigr|.
\]
Therefore, by \eqref{degenerate},
\[
\Var(\mc A_n(a))
=
\mb{E}\!\left[q_{12}(a)\bigl(1-q_{12}(a)\bigr)\right].
\]

We next construct an event on which \(|\Delta_{12}(a)|\) is of order \(\delta_p^{-1}\) and which has probability bounded away from \(0\), where $\delta_p$ is defined as in \eqref{delta-p}.
Let
\[
w_p(y):=\exp(\kappa y^2)(1-y^2)^{(p-3)/2}\mathbf 1_{[-1,1]}(y),
\qquad
W_p:=\int_{-1}^1 w_p(y)\,dy.
\]
Then \(Y_1\) has density \(w_p/W_p\).

Let us estimate \(W_p\) from above. Since \(\kappa=p/2-\delta_p\) and \(\log(1-u)\le -u\) for \(u\in[0,1)\),
\[
\kappa y^2+\frac{p-3}{2}\log(1-y^2)
\le
\Bigl(\frac p2-\delta_p\Bigr)y^2-\frac{p-3}{2}y^2
=
-\Bigl(\delta_p-\frac32\Bigr)y^2.
\]
Hence, for all sufficiently large \(p\),
\[
w_p(y)\le \exp\!\left(-\frac{\delta_p}{2}y^2\right),
\]
and therefore
\[
W_p\le \int_{\mb{R}}\exp\!\left(-\frac{\delta_p}{2}y^2\right)\,dy
\le \frac{C}{\sqrt{\delta_p}}
\]
for some universal constant \(C>0\).

Next, fix any constants \(0<c_1<c_2<\infty\). We claim that
\[
\mb{P}\!\left(\frac{c_1}{\sqrt{\delta_p}}\le |Y_1|\le \frac{c_2}{\sqrt{\delta_p}}\right)\ge \varrho
\]
for some \(\varrho>0\) and all sufficiently large \(p\).
Indeed, write \(y=u/\sqrt{\delta_p}\). For \(u\in[c_1,c_2]\) and all sufficiently large \(p\), we have
\(u^2/\delta_p\le 1/2\), so using \(\log(1-v)\ge -v-v^2\) for \(v\in[0,1/2]\),
\begin{align*}
\kappa \frac{u^2}{\delta_p}+\frac{p-3}{2}\log\!\left(1-\frac{u^2}{\delta_p}\right)
&\ge
\Bigl(\frac p2-\delta_p\Bigr)\frac{u^2}{\delta_p}
-\frac{p-3}{2}\frac{u^2}{\delta_p}
-\frac{p-3}{2}\frac{u^4}{\delta_p^2} \\
&=
-\Bigl(1-\frac{3}{2\delta_p}\Bigr)u^2
-\frac{p-3}{2}\frac{u^4}{\delta_p^2}.
\end{align*}
Since \(p/\delta_p^2\to 0\), the right-hand side is bounded below by a constant depending only on \(c_2\).
Hence there exists \(c_*>0\) such that for all sufficiently large \(p\),
\[
w_p\!\left(\frac{u}{\sqrt{\delta_p}}\right)\ge c_*
\qquad\text{for all }u\in[c_1,c_2].
\]
Therefore
\begin{align*}
\mb{P}\!\left(\frac{c_1}{\sqrt{\delta_p}}\le |Y_1|\le \frac{c_2}{\sqrt{\delta_p}}\right)
&=
\frac{2}{W_p}\int_{c_1/\sqrt{\delta_p}}^{c_2/\sqrt{\delta_p}} w_p(y)\,dy \\
&=
\frac{2}{W_p\sqrt{\delta_p}}\int_{c_1}^{c_2} w_p\!\left(\frac{u}{\sqrt{\delta_p}}\right)\,du \\
&\ge
\frac{2c_*(c_2-c_1)}{C}
=:\varrho>0.
\end{align*}
Consequently, the event
\[
E_p:=\left\{\frac{c_1}{\sqrt{\delta_p}}\le |Y_1|,|Y_2|\le \frac{c_2}{\sqrt{\delta_p}}\right\}
\]
satisfies
\[
\mb{P}(E_p)\ge \varrho^2.
\]
Now, on \(E_p\), we have $
|Y_1Y_2|\asymp \delta_p^{-1}$ and $
1-\alpha_{12}
=
1-\sqrt{(1-Y_1^2)(1-Y_2^2)}
=
O(\delta_p^{-1})
$, so
\[
\left|\frac1{\alpha_{12}}-1\right|=O(\delta_p^{-1}).
\]
Therefore, on \(E_p\),
\[
\left|\frac{a}{\sqrt p}\Bigl(\frac1{\alpha_{12}}-1\Bigr)\right|
=
O_a\!\left(\frac{1}{\sqrt p\,\delta_p}\right),
\]
whereas
\[
\left|\frac{\beta_{12}}{\alpha_{12}}\right|
\asymp \delta_p^{-1}.
\]
Consequently, for all sufficiently large \(p\), there exists \(c\in (0,1)\), indepdent of $p$ and $n$, such that on \(E_p\),
\[
c^{-1}\,\delta_p^{-1}\geq |\Delta_{12}(a)|\ge c\,\delta_p^{-1},
\]
Observe that on \(E_p\),
\[
\left|\frac{a}{\sqrt p}+\theta \Delta_{12}(a)\right|
\le \frac{|a|+1}{\sqrt p}
\qquad\text{for all }\theta\in[0,1]
\]
when \(p\) is large enough. 

Furthermore,  there exists a constant \(c_a' \in (0,1)\) such that for all sufficiently large \(p\),
\[
\frac{\sqrt{p}}{c_a'} \geq f_p(x)\ge c_a' \sqrt p
\qquad\text{for all }|x|\le \frac{|a|+1}{\sqrt p}.
\]
Therefore, by the mean value theorem, on \(E_p\),
\[
q_{12}(a)
=
\bigl|F_p(a/\sqrt p+\Delta_{12}(a))-F_p(a/\sqrt p)\bigr|
=
|f_p(\xi_{12})|\,|\Delta_{12}(a)|
\]
for some \(\xi_{12}\) between \(a/\sqrt p\) and \(a/\sqrt p+\Delta_{12}(a)\). Hence
\[
q_{12}(a)\ge c_a' \sqrt p\,|\Delta_{12}(a)|
\ge c_a' c\,\frac{\sqrt p}{\delta_p}
\qquad\text{on }E_p.
\]
Moreover,
\[q_{12}(a)=|f_p(\xi_{12})|\,|\Delta_{12}(a)| \lesssim  \frac{\sqrt{p}}{\delta_p}.  \]
Thus \(q_{12}(a)\le 1/2\) on \(E_p\) for all large \(p\) because $\delta_p \gg \sqrt{p}$. Consequently,
\[
q_{12}(a)\bigl(1-q_{12}(a)\bigr)\ge \frac12 q_{12}(a)
\ge c_a''\,\frac{\sqrt p}{\delta_p}
\qquad\text{on }E_p
\]
for some constant \(c_a''>0\).

Taking expectation, we conclude that for all sufficiently large \(p\),
\[
\Var(\mc A_n(a))
=
\mb{E}\!\left[q_{12}(a)\bigl(1-q_{12}(a)\bigr)\right]
\ge
\mb{E}\!\left[q_{12}(a)\bigl(1-q_{12}(a)\bigr)\mathbf 1_{E_p}\right]
\ge
c_a''\,\frac{\sqrt p}{\delta_p}\mb{P}(E_p).
\]
Since \(\mb{P}(E_p)\ge \varrho^2\), this gives \eqref{lower bound variance An}. $\hfill$ $\square$

The next Lemma gives a variance bound on $\mc{B}_n(t)$ over a fixed interval $[-a,a]$. Its proof is based on Lemma \ref{Brascamp-Lieb}. 

\begin{lemma} \label{variance bound Brascamp-Lieb}
    Recall $\mc{B}_n$ in \eqref{Bn(t)}. Then, for any fixed $a>0$, we have 
    \begin{align} \label{variance bound on derivatives}
        \sup_{t\in [-a,a]} \la  \Var \lb \mc{B}_n(t) \rb + \Var \lb \mc{B}'_n(t) \rb + \Var \lb \mc{B}''_n(t) \rb \ra \leq C_a \cdot \lb \frac{p}{\delta_p} + \frac{p^2}{\delta_p^3} \rb
    \end{align}
    where the derivatives in the expression above are taken with respect to $t$.
\end{lemma}

\noindent \textbf{Proof of Lemma \ref{variance bound Brascamp-Lieb}.}
Let us introduce some notation so that it is easier to keep track of the notation. Define
\begin{align} \label{alpha xy}
\alpha(x,y):=\sqrt{(1-x^2)(1-y^2)},
\qquad
s_t(x,y):=\frac{t-\sqrt p\,xy}{\alpha(x,y)},
\end{align}
and
\[
\varphi_t(x,y)=F_p(s_t(x,y))
\]
where we recall that $F_p$ is the CDF of $\sqrt{p}R_{12}$. 

The derivatives of $\varphi$ with respect to $t$ is then 
\[
\varphi_t^{(v)}(x,y)=\alpha(x,y)^{-v}F_p^{(v)}(s_t(x,y)),
\qquad v=0,1,2.
\]
It is easy to see that for \(v=0,1,2\),
\begin{equation}\label{eq:G-deriv-bdd}
\sup_{p\ge 3}\sup_{u\in\mb R}
\Bigl(
|F_p^{(v)}(u)|+|F_p^{(v+1)}(u)|+|F_p^{(v+2)}(u)|
\Bigr)
<\infty.
\end{equation}
By using Lemma \ref{Brascamp-Lieb}, for all $v \in \la 0,1,2 \ra$, we have the bound 
\begin{align} \label{poincare}
    \Var \lb \mc{B}^{(v)}_n(t) \rb \lesssim \frac{1}{\delta_p} \sum_{i=1}^n \mb{E} \left[ (\partial_i\mc B_n^{(v)}(t))^2 \right].
\end{align}
Here we should clarify that the derivative $\partial_i$ is taken with respect to $Y_i$. 

Conditionally on \(Y_i\), the variables \(\{\partial_x\varphi_t^{(v)}(Y_i,Y_j):j\ne i\}\) are i.i.d., so
\begin{align} \label{variance decomp}
\mb E\bigl[(\partial_i\mc B_n^{(v)}(t))^2\bigr]
&=
\frac{1}{N_n}
\mb E\!\left[
(n-1)\Var\bigl(\partial_x\varphi_t^{(v)}(Y_i,Y_2)\mid Y_i\bigr)
+
(n-1)^2
\Bigl(\mb E[\partial_x\varphi_t^{(v)}(Y_i,Y_2)\mid Y_i]\Bigr)^2
\right]  \nonumber \\
&\lesssim
\frac{1}{n}\,
\mb E\bigl[(\partial_x\varphi_t^{(v)}(Y_1,Y_2))^2\bigr]
+
\,
\mb E\bigl[m_{v,t}(Y_1)^2\bigr].
\end{align}
where
\[
m_{v,t}(x):=\mb E\bigl[\partial_x\varphi_t^{(v)}(x,Y_2)\bigr].
\]
The proof will be completed if we show that
\begin{align}
   \max_{v \in \la 0,1,2 \ra} \mathbb E\bigl[(\partial_x\varphi^{(v)}_t(Y_1,Y_2))^2\bigr] &\leq C_a \frac{p}{\delta_p}; \label{phi-estimate}\\
      \max_{v \in \la 0,1,2 \ra} \mathbb E[m_{v,t}(Y_1)^2] &\leq C_a\,\frac{p^2}{\delta_p^3}. \label{m-estimate}
\end{align}
uniformly in $|t| \leq a$. 

In fact, given \eqref{phi-estimate} and \eqref{m-estimate}, we have
\[
\max_{v \in \la 0,1,2 \ra} \Var\!\bigl(\mc{B}^{(v)}_n(t)\bigr)
\lesssim
\max_{v \in \la 0,1,2 \ra} \la \frac{1}{\delta_p}\sum_{i=1}^n
\mathbb E\bigl[(\partial_i \mc{B}^{(v)}_n(t))^2\bigr] \ra
\le
C_a\left(
\frac{p}{\delta_p^2}
+
\frac{n p^2}{\delta_p^4}
\right)
\]
by using \eqref{variance decomp}. The right-hand side in the display above goes to $0$ under the assumption $\delta_p=p/2-\kappa \gg p^{1/2}n^{1/4}$.

Let us start with \eqref{phi-estimate}. Since
\[
\alpha_x(x,y)=-\frac{x}{1-x^2}\alpha(x,y),
\qquad
\partial_x s_t(x,y)
=
\frac{tx-\sqrt p\,y}{(1-x^2)^{3/2}(1-y^2)^{1/2}},
\]
we obtain
\[
\partial_x\varphi_t^{(v)}(x,y)
=
\alpha(x,y)^{-v}
\left[
F_p^{(v+1)}(s_t(x,y))\,\partial_x s_t(x,y)
+
v \frac{x}{1-x^2}F_p^{(v)}(s_t(x,y))
\right].
\]
Therefore, by \eqref{eq:G-deriv-bdd} and \(|t|\le a\),
\begin{equation}\label{eq:varphi-x-pointwise}
|\partial_x\varphi_t^{(v)}(x,y)|^2
\le
C_{a,v}
\frac{x^2+p y^2}{(1-x^2)^{v+3}(1-y^2)^{v+1}}.
\end{equation}
Using independence and \eqref{eq:simple-weighted-delta}, we get
\begin{align}
& \max_{v \in \la 0,1,2 \ra}  \mb E\bigl[(\partial_x\varphi_t^{(v)}(Y_1,Y_2))^2\bigr] \nonumber \\
\le &
C_{a,v} \cdot \max_{v \in \la 0,1,2 \ra} \la
\underbrace{\mb E\!\left[\frac{Y_1^2}{(1-Y_1^2)^{v+3}}\right]
\mb E\!\left[\frac{1}{(1-Y_2^2)^{v+1}}\right]}_{O \lb \delta_p^{-1} \rb}
\nonumber
+
p\,
\underbrace{\mb E\!\left[\frac{1}{(1-Y_1^2)^{v+3}}\right]
\mb E\!\left[\frac{Y_2^2}{(1-Y_2^2)^{v+1}}\right]}_{O \lb \delta_p^{-1} \rb} \ra
\nonumber\\
&\le
C_{a,v}\frac{p}{\delta_p}.
\label{eq:varphi-x-L2}
\end{align}
This concludes the proof of \eqref{phi-estimate}. 

We now prove \eqref{m-estimate}. Let us first claim that
\begin{equation}\label{eq:m-zero}
m_{v,t}(0)=0.
\end{equation}
for all $v \in \la 0,1,2 \ra$. Indeed, we can get the explicit formula
\[
\partial_x \varphi_t^{(v)} \lb 0,Y_2 \rb = \frac{1}{(1-Y_2^2)^{-v/2}} \cdot F_p^{(v+1)} \lb \frac{t}{\sqrt{1-Y_2^2}} \rb \cdot \lb \frac{-\sqrt{p} Y_2}{\sqrt{1-Y_2^2}} \rb        
\]
 Since the law of \(Y_2\) is symmetric, \eqref{eq:m-zero} follows.

Next, differentiate once more in \(x\), we arrive at
\[
\partial_x^2 s_t(x,y)
=
\frac{t(1+2x^2)-3x\sqrt p\,y}{(1-x^2)^{5/2}(1-y^2)^{1/2}}.
\]
Thanks to \eqref{eq:G-deriv-bdd}, we obtain the rough pointwise bound
\begin{equation}\label{eq:varphi-xx-pointwise}
|\partial_x^2\varphi_t^{(v)}(x,y)|
\le
C_{a,v}
\frac{1+\sqrt p\,|y|+p y^2}{(1-x^2)^{v+3}(1-y^2)^{v+2}}.
\end{equation}
Therefore, using \eqref{eq:simple-weighted-delta} and Cauchy--Schwarz,
\begin{align}
|m_{v,t}'(x)|
&\le
\mb E\bigl[|\partial_x^2\varphi_t^{(v)}(x,Y_2)|\bigr]
\nonumber\\
&\le
\frac{C_{a,v}}{(1-x^2)^{v+3}}
\left(
\mb E\frac{1}{(1-Y_2^2)^{v+2}}
+
\sqrt p\,\mb E\frac{|Y_2|}{(1-Y_2^2)^{v+2}}
+
p\,\mb E\frac{Y_2^2}{(1-Y_2^2)^{v+2}}
\right)
\nonumber\\
&\le
\frac{C_{a,v}}{(1-x^2)^{v+3}}
\left(
1+\sqrt{\frac{p}{\delta_p}}+\frac{p}{\delta_p}
\right)
\nonumber\\
&\le
C_{a,v}\frac{p/\delta_p}{(1-x^2)^{v+3}}.
\label{eq:m-prime}
\end{align}
Combining \eqref{eq:m-zero} and \eqref{eq:m-prime},
\[
|m_{v,t}(x)|
=
\left|\int_0^x m_{v,t}'(u)\,du\right|
\le
C_{a,v} \cdot \frac{p}{\delta_p} \cdot
\int_0^{|x|}\frac{du}{(1-u^2)^{v+3}}
\le
C_{a,v}\cdot \frac{p}{\delta_p} \cdot \frac{|x|}{(1-x^2)^{v+3}}.
\]
Finally, by using \eqref{eq:simple-weighted-delta}, 
\begin{equation*}
\mb E\bigl[m_{v,t}(T_1)^2\bigr]
\le
C_{a,v}\cdot \frac{p^2}{\delta_p^2} \cdot 
\mb E\!\left[\frac{Y_1^2}{(1-Y_1^2)^{2v+6}}\right]
\le
C_{a,v}\frac{p^2}{\delta_p^3}.
\end{equation*}
This concludes the proof of \eqref{m-estimate}. The proof is completed. $\hfill$ $\square$

The following variance estimate of $\mc{B}'_n(t)$ will also be used in the proof of our main result. 

\begin{lemma} \label{tail variance Bn(t)}
Recall \(\mc{B}_n(t)\) in \eqref{Bn(t)}. Then there exist universal constants
\(C,c>0\) such that, for all \(a\ge 1\) and all \(p\) large enough,
\begin{equation}\label{eq:Bnprime-tail-var}
 \Var\!\bigl(\mc B_n'(t)\bigr)
\le
C \cdot \Bigl[ e^{-c t^2} + (1+|t|)^{-4} \Bigr]
\end{equation}
whenever $|t|>a$.
\end{lemma}

\noindent \textbf{Proof of Lemma \ref{tail variance Bn(t)}.} We prove the stronger estimate
\begin{equation}\label{eq:Bnprime-tail-var-stronger}
 \Var\!\bigl(\mc B_n'(t)\bigr)
\le
C\left(
\frac{p}{\delta_p^2}
+
\frac{n p^2}{\delta_p^4}
\right)
\cdot \Bigl[e^{-c t^2}+(1+|t|)^{-4}\Bigr]
\end{equation}
for universal constants \(C,c>0\), all \(a\ge 1\), and all \(p\) large enough.

As before, by \eqref{poincare}, it is enough to bound
\[
\mathbb E\bigl[(\partial_i \mc B_n'(t))^2\bigr]
\]
uniformly in \(|t|>a\). By \eqref{variance decomp}  with \(v=1\),
\[
\mathbb E\bigl[(\partial_i\mc B_n'(t))^2\bigr]
\le
\frac{2}{n}\,\mathbb E\bigl[(\partial_x\varphi_t^{(1)}(Y_1,Y_2))^2\bigr]
+
2\,\mathbb E[m_{1,t}(Y_1)^2],
\]
where
\[
m_{1,t}(x):=\mathbb E\bigl[\partial_x\varphi_t^{(1)}(x,Y_2)\bigr].
\]
Hence it is enough to show that, uniformly in \(|t|>a\),
\begin{align}
   \mathbb E\bigl[(\partial_x\varphi_t^{(1)}(Y_1,Y_2))^2\bigr]
   &\le C \frac{p}{\delta_p}
   \Bigl(e^{-c t^2}+(1+|t|)^{-4}\Bigr), \label{phi-estimate-tail}\\
   \mathbb E[m_{1,t}(Y_1)^2]
   &\le C \frac{p^2}{\delta_p^3}
   \Bigl(e^{-c t^2}+(1+|t|)^{-4}\Bigr). \label{m-estimate-tail}
\end{align}
Indeed, given \eqref{phi-estimate-tail} and \eqref{m-estimate-tail}, we get
\[
\mathbb E\bigl[(\partial_i\mc B_n'(t))^2\bigr]
\le
C\left[
\frac{1}{n}\cdot \frac{p}{\delta_p}
\Bigl(e^{-c t^2}+(1+|t|)^{-4}\Bigr)
+
\frac{p^2}{\delta_p^3}
\Bigl(e^{-c t^2}+(1+|t|)^{-4}\Bigr)
\right],
\]
and therefore, by \eqref{poincare},
\[
\Var\!\bigl(\mc B_n'(t)\bigr)
\le
\frac{C}{\delta_p}\sum_{i=1}^n
\mathbb E\bigl[(\partial_i\mc B_n'(t))^2\bigr]
\le
C\left(
\frac{p}{\delta_p^2}
+
\frac{n p^2}{\delta_p^4}
\right)
\Bigl(e^{-c t^2}+(1+|t|)^{-4}\Bigr),
\]
which is exactly \eqref{eq:Bnprime-tail-var-stronger}.

\medskip

\noindent
\underline{\textit{Proof of \eqref{phi-estimate-tail}.}}
Recall that
\[
\partial_x\varphi_t^{(1)}(x,y)
=
\alpha(x,y)^{-1}
\left[
F_p^{(2)}(s_t(x,y))\,\partial_x s_t(x,y)
+
\frac{x}{1-x^2}F_p^{(1)}(s_t(x,y))
\right].
\]
We will use the fact that there exist universal constants \(C,c>0\) such that
\begin{equation} \label{eq:gauss-tail-F}
|F_p^{(1)}(u)|+|F_p^{(2)}(u)|+|F_p^{(3)}(u)|
\le C e^{-c u^2},
\qquad u\in\mathbb R,
\end{equation}
uniformly in \(p\ge 3\). This follows from the explicit formula for the density of \(\sqrt p\,R_{12}\) (which is asymptotically a standard Gaussian).

Recall \(\beta_{ij}\), \(\alpha(x,y)\) and \(s_t(x,y)\) in \eqref{beta ij} and \eqref{alpha xy}. Put
\[
E_t:=\Bigl\{ \sqrt{p}\,|\beta_{12}|\le |t|/2\Bigr\}.
\]
On the event \(E_t\), we have
\[
|t-\sqrt{p}\,\beta_{12}|\ge |t|/2.
\]
Moreover, since \(\alpha(Y_1,Y_2)\le 1\), it follows that
\[
|s_t(Y_1,Y_2)|
=
\frac{|t-\sqrt{p}\,\beta_{12}|}{\alpha(Y_1,Y_2)}
\ge \frac{|t|}{2}
\]
on \(E_t\). Hence, by \eqref{eq:gauss-tail-F}, on \(E_t\),
\[
|F_p^{(1)}(s_t(Y_1,Y_2))|+|F_p^{(2)}(s_t(Y_1,Y_2))|
\le C e^{-c t^2}.
\]
Combining this with
\[
\partial_x s_t(x,y)
=
\frac{tx-\sqrt p\,y}{(1-x^2)^{3/2}(1-y^2)^{1/2}},
\]
we get, on \(E_t\),
\begin{align*}
|\partial_x\varphi_t^{(1)}(Y_1,Y_2)|^2
&\le
C e^{-c t^2}
\left(
\frac{(tY_1-\sqrt p\,Y_2)^2}{(1-Y_1^2)^4(1-Y_2^2)^2}
+
\frac{Y_1^2}{(1-Y_1^2)^3(1-Y_2^2)}
\right) \\
&\le
C e^{-c t^2}
\left(
\frac{Y_1^2+p Y_2^2}{(1-Y_1^2)^4(1-Y_2^2)^2}
+
\frac{Y_1^2}{(1-Y_1^2)^3(1-Y_2^2)}
\right),
\end{align*}
where we used that \((1+t^2)e^{-c t^2}\lesssim e^{-c t^2/2}\) for all \(t\in\mathbb R\). Taking expectation and using \eqref{eq:simple-weighted-delta}, we obtain
\begin{equation}\label{eq:phi-tail-good}
\mathbb E\Bigl[(\partial_x\varphi_t^{(1)}(Y_1,Y_2))^2 \mathbf 1_{E_t}\Bigr]
\le
C \frac{p}{\delta_p} e^{-c t^2}.
\end{equation}

On \(E_t^c\), we use boundedness of \(F_p^{(1)}\) and \(F_p^{(2)}\). Since \(|t|\le 2\sqrt{p}\,|\beta_{12}|\) on \(E_t^c\), we have
\[
|tY_1-\sqrt p\,Y_2|
\le |t|\,|Y_1|+\sqrt p\,|Y_2|
\le 2\sqrt p\,Y_1^2|Y_2|+\sqrt p\,|Y_2|
\le 3\sqrt p\,|Y_2|.
\]
Therefore,
\begin{equation}\label{eq:phi-tail-bad-pointwise}
|\partial_x\varphi_t^{(1)}(Y_1,Y_2)|^2 \mathbf 1_{E_t^c}
\le
C\left(
\frac{p Y_2^2}{(1-Y_1^2)^4(1-Y_2^2)^2}
+
\frac{Y_1^2}{(1-Y_1^2)^3(1-Y_2^2)}
\right)\mathbf 1_{E_t^c}.
\end{equation}
Next, by using
\[
\mathbf 1_{E_t^c}
=
\mathbf 1_{\{\sqrt{p}\,|\beta_{12}|>|t|/2\}}
\le
\frac{16p^2 \beta_{12}^4}{t^4}
=
\frac{16 p^2 Y_1^4Y_2^4}{t^4},
\]
we get
\begin{align*}
\mathbb E\Bigl[
\frac{pY_2^2}{(1-Y_1^2)^4(1-Y_2^2)^2}\mathbf 1_{E_t^c}
\Bigr]
&\le
\frac{C p^3}{t^4}
\mathbb E\!\left[\frac{Y_1^4}{(1-Y_1^2)^4}\right]
\mathbb E\!\left[\frac{Y_2^6}{(1-Y_2^2)^2}\right]
\le
\frac{C p^3}{\delta_p^5 t^4},
\\
\mathbb E\Bigl[
\frac{Y_1^2}{(1-Y_1^2)^3(1-Y_2^2)}\mathbf 1_{E_t^c}
\Bigr]
&\le
\frac{C p^2}{t^4}
\mathbb E\!\left[\frac{Y_1^6}{(1-Y_1^2)^3}\right]
\mathbb E\!\left[\frac{Y_2^4}{1-Y_2^2}\right]
\le
\frac{C p^2}{\delta_p^5 t^4}.
\end{align*}
Plugging these bounds into \eqref{eq:phi-tail-bad-pointwise}, we obtain
\begin{equation}\label{eq:phi-tail-bad}
\mathbb E\Bigl[(\partial_x\varphi_t^{(1)}(Y_1,Y_2))^2 \mathbf 1_{E_t^c}\Bigr]
\le
C\left(
\frac{p^3}{\delta_p^5}
+
\frac{p^2}{\delta_p^5}
\right)(1+|t|)^{-4}.
\end{equation}
Since \(\delta_p\gg p^{1/2}\), we have, for all large \(n\),
\[
\frac{p^3}{\delta_p^5}
+
\frac{p^2}{\delta_p^5}
\lesssim
\frac{p}{\delta_p}.
\]
Combining this with \eqref{eq:phi-tail-good} and \eqref{eq:phi-tail-bad}, we obtain \eqref{phi-estimate-tail}.

\medskip

\noindent
\underline{\textit{Proof of \eqref{m-estimate-tail}.}}
Recall that
\[
m_{1,t}(x):=\mathbb E\bigl[\partial_x\varphi_t^{(1)}(x,Y_2)\bigr].
\]
As in the proof above, we have \(m_{1,t}(0)=0\).

We first prove the bound
\begin{equation}\label{eq:phi-xx-rough-tail}
|\partial_x^2\varphi_t^{(1)}(x,y)|
\lesssim
\frac{1+\sqrt p\,|y|+p y^2}{(1-x^2)^4(1-y^2)^3},
\qquad x,y\in(-1,1),\ t\in\mathbb R.
\end{equation}
To do this, note that \eqref{eq:gauss-tail-F} gives
\begin{equation} \label{eq:Fp-weighted-derivatives}
\sup_{p\ge 100}\sup_{u\in\mathbb R}
(1+|u|+u^2)\Bigl(
|F_p^{(1)}(u)|+|F_p^{(2)}(u)|+|F_p^{(3)}(u)|
\Bigr)<\infty.
\end{equation}
Put \(a(x)=x/(1-x^2)\). Since
\[
\alpha_x(x,y)=-a(x)\alpha(x,y),
\]
we have
\[
(\alpha^{-1})_x=a(x)\alpha^{-1}.
\]
Also, recalling \(s_t\) from \eqref{alpha xy},
\[
t-\sqrt p\,xy=\alpha(x,y)s_t(x,y),
\]
and hence
\begin{equation}\label{eq:sx-rewrite}
\partial_x s_t(x,y)
=
-\sqrt p\,y\,\alpha(x,y)^{-1}+a(x)s_t(x,y).
\end{equation}
Differentiating once more gives
\begin{equation}\label{eq:sxx-rewrite}
\partial_x^2 s_t(x,y)
=
\bigl(a'(x)+a(x)^2\bigr)s_t(x,y)-2a(x)\sqrt p\,y\,\alpha(x,y)^{-1}.
\end{equation}
Since
\[
a'(x)=\frac{1+x^2}{(1-x^2)^2},
\qquad
a'(x)+a(x)^2=\frac{1+2x^2}{(1-x^2)^2},
\]
it follows from \eqref{eq:sx-rewrite} and \eqref{eq:sxx-rewrite} that
\begin{align}
|\partial_x s_t(x,y)|
&\lesssim
\frac{|s_t(x,y)|}{1-x^2}
+
\frac{\sqrt p\,|y|}{\sqrt{(1-x^2)(1-y^2)}},
\label{eq:sx-bound}
\\
|\partial_x^2 s_t(x,y)|
&\lesssim
\frac{|s_t(x,y)|}{(1-x^2)^2}
+
\frac{\sqrt p\,|y|}{(1-x^2)^{3/2}(1-y^2)^{1/2}}.
\label{eq:sxx-bound}
\end{align}
Now
\[
\varphi_t^{(1)}(x,y)=\alpha(x,y)^{-1}F_p^{(1)}(s_t(x,y)),
\]
so
\[
\partial_x\varphi_t^{(1)}(x,y)
=
\alpha(x,y)^{-1}
\left[
F_p^{(2)}(s_t(x,y))\,\partial_x s_t(x,y)
+
a(x)F_p^{(1)}(s_t(x,y))
\right].
\]
Differentiating once more and using \((\alpha^{-1})_x=a(x)\alpha^{-1}\), we get
\begin{align} \label{d^2/dx^2 varphi}
\partial_x^2\varphi_t^{(1)}(x,y)
=
\alpha(x,y)^{-1}\Big[
&F_p^{(3)}(s_t(x,y))\,(\partial_x s_t(x,y))^2
+
F_p^{(2)}(s_t(x,y))\,\partial_x^2 s_t(x,y) \nonumber
\\
&\qquad\qquad
+2a(x)F_p^{(2)}(s_t(x,y))\,\partial_x s_t(x,y)
+
\bigl(a'(x)+a(x)^2\bigr)F_p^{(1)}(s_t(x,y))
\Big].
\end{align}
We now bound the four terms on the right-hand side separately. For the first term, by \eqref{eq:Fp-weighted-derivatives},
\[
|F_p^{(3)}(u)|\le \frac{C}{1+u^2},
\]
hence by \eqref{eq:sx-bound},
\begin{align*}
|F_p^{(3)}(s_t)|\,|\partial_x s_t|^2
&\le
\frac{C}{1+s_t^2}
\left(
\frac{|s_t|}{1-x^2}
+
\frac{\sqrt p\,|y|}{\sqrt{(1-x^2)(1-y^2)}}
\right)^2
\\
&\le
C\left(
\frac{1}{(1-x^2)^2}
+
\frac{p y^2}{(1-x^2)(1-y^2)}
\right).
\end{align*}

For the second term, by \eqref{eq:Fp-weighted-derivatives},
\[
|F_p^{(2)}(u)|\le \frac{C}{1+|u|},
\]
and therefore, using \eqref{eq:sxx-bound},
\begin{align*}
|F_p^{(2)}(s_t)|\,|\partial_x^2 s_t|
&\le
\frac{C}{1+|s_t|}
\left(
\frac{|s_t|}{(1-x^2)^2}
+
\frac{\sqrt p\,|y|}{(1-x^2)^{3/2}(1-y^2)^{1/2}}
\right)
\\
&\le
C\left(
\frac{1}{(1-x^2)^2}
+
\frac{\sqrt p\,|y|}{(1-x^2)^{3/2}(1-y^2)^{1/2}}
\right).
\end{align*}

For the third term, using again \(|F_p^{(2)}(u)|\le C/(1+|u|)\), \eqref{eq:sx-bound}, and
\[
|a(x)|\le \frac{1}{1-x^2},
\]
we obtain
\begin{align*}
|a(x)F_p^{(2)}(s_t)\partial_x s_t|
&\le
\frac{C}{1-x^2}\cdot \frac{1}{1+|s_t|}
\left(
\frac{|s_t|}{1-x^2}
+
\frac{\sqrt p\,|y|}{\sqrt{(1-x^2)(1-y^2)}}
\right)
\\
&\le
C\left(
\frac{1}{(1-x^2)^2}
+
\frac{\sqrt p\,|y|}{(1-x^2)^{3/2}(1-y^2)^{1/2}}
\right).
\end{align*}

For the last term, \eqref{eq:Fp-weighted-derivatives} gives \(|F_p^{(1)}(u)|\le C\), hence
\[
\bigl(a'(x)+a(x)^2\bigr)|F_p^{(1)}(s_t)|
\le
C\frac{1}{(1-x^2)^2}.
\]

Putting the four bounds together, we find
\[
|\partial_x^2\varphi_t^{(1)}(x,y)|
\le
C\,\alpha(x,y)^{-1}
\left(
\frac{1}{(1-x^2)^2}
+
\frac{\sqrt p\,|y|}{(1-x^2)^{3/2}(1-y^2)^{1/2}}
+
\frac{p y^2}{(1-x^2)(1-y^2)}
\right).
\]
Since
\[
\alpha(x,y)^{-1}=\frac{1}{\sqrt{(1-x^2)(1-y^2)}},
\]
this implies
\[
|\partial_x^2\varphi_t^{(1)}(x,y)|
\le
C\left(
\frac{1}{(1-x^2)^{5/2}(1-y^2)^{1/2}}
+
\frac{\sqrt p\,|y|}{(1-x^2)^2(1-y^2)}
+
\frac{p y^2}{(1-x^2)^{3/2}(1-y^2)^{3/2}}
\right).
\]
Finally, since \(1-x^2\le 1\) and \(1-y^2\le 1\), each term on the right-hand side is bounded by
\[
C\frac{1+\sqrt p\,|y|+p y^2}{(1-x^2)^4(1-y^2)^3},
\]
which proves \eqref{eq:phi-xx-rough-tail}.

Fix \(x\in(-1,1)\), and define
\[
E_t(x):=\Bigl\{|\sqrt p\,xY_2|\le |t|/2\Bigr\}.
\]
On \(E_t(x)\), the same argument as above gives
\[
|s_t(x,Y_2)|\ge |t|/2,
\]
hence, by \eqref{eq:gauss-tail-F} and \eqref{eq:phi-xx-rough-tail},
\[
|\partial_x^2\varphi_t^{(1)}(x,Y_2)| \cdot \mathbf 1_{E_t(x)}
\le
C e^{-c t^2}
\frac{1+\sqrt p\,|Y_2|+pY_2^2}{(1-x^2)^4(1-Y_2^2)^3}.
\]
Taking expectation and using \eqref{eq:simple-weighted-delta} together with Cauchy--Schwarz,
\[
\left|
\mathbb E\bigl[\partial_x^2\varphi_t^{(1)}(x,Y_2)\mathbf 1_{E_t(x)}\bigr]
\right|
\le
C \frac{p}{\delta_p}\frac{e^{-c t^2}}{(1-x^2)^4}.
\]
On \(E_t(x)^c\), we use \eqref{eq:phi-xx-rough-tail} and
\[
\mathbf 1_{E_t(x)^c}
=
\mathbf 1_{\{|\sqrt p\,xY_2|>|t|/2\}}
\le
\frac{4 p x^2 Y_2^2}{t^2}.
\]
Therefore
\begin{align*}
\left|
\mathbb E\bigl[\partial_x^2\varphi_t^{(1)}(x,Y_2)\mathbf 1_{E_t(x)^c}\bigr]
\right|
&\le
\mathbb E\Bigl[|\partial_x^2\varphi_t^{(1)}(x,Y_2)|\mathbf 1_{E_t(x)^c}\Bigr] \\
&\le
\frac{C p x^2}{t^2(1-x^2)^4}
\mathbb E\!\left[
\frac{(1+\sqrt p\,|Y_2|+pY_2^2)Y_2^2}{(1-Y_2^2)^3}
\right].
\end{align*}
Using \eqref{eq:simple-weighted-delta} and Cauchy--Schwarz once more,
\begin{align*}
\mathbb E\!\left[\frac{Y_2^2}{(1-Y_2^2)^3}\right] &\le C\delta_p^{-1},\\
\sqrt p\,\mathbb E\!\left[\frac{|Y_2|^3}{(1-Y_2^2)^3}\right]
&\le
\sqrt p\left(
\mathbb E\!\left[\frac{Y_2^2}{(1-Y_2^2)^3}\right]
\mathbb E\!\left[\frac{Y_2^4}{(1-Y_2^2)^3}\right]
\right)^{1/2}
\le C \frac{p}{\delta_p^2},\\
p\,\mathbb E\!\left[\frac{Y_2^4}{(1-Y_2^2)^3}\right]
&\le C \frac{p}{\delta_p^2}.
\end{align*}
Hence
\[
\left|
\mathbb E\bigl[\partial_x^2\varphi_t^{(1)}(x,Y_2)\mathbf 1_{E_t(x)^c}\bigr]
\right|
\le
C \frac{p^2}{\delta_p^2}\cdot \frac{x^2}{(1+|t|)^2(1-x^2)^4}.
\]

Combining the two pieces, we have shown that
\[
|m_{1,t}'(x)|
\le
C \frac{p}{\delta_p}\frac{e^{-c t^2}}{(1-x^2)^4}
+
C \frac{p^2}{\delta_p^2}\cdot \frac{x^2}{(1+|t|)^2(1-x^2)^4}.
\]
Since \(m_{1,t}(0)=0\),
\begin{align*}
|m_{1,t}(x)|
&=
\left|\int_0^x m_{1,t}'(u)\,du\right| \\
&\le
C \frac{p}{\delta_p} e^{-c t^2}
\int_0^{|x|}\frac{du}{(1-u^2)^4}
+
C \frac{p^2}{\delta_p^2(1+|t|)^2}
\int_0^{|x|}\frac{u^2\,du}{(1-u^2)^4} \\
&\le
C \frac{p}{\delta_p} e^{-c t^2}\frac{|x|}{(1-x^2)^4}
+
C \frac{p^2}{\delta_p^2(1+|t|)^2}\frac{|x|^3}{(1-x^2)^4}.
\end{align*}
Therefore, using \((u+v)^2\le 2u^2+2v^2\) and \eqref{eq:simple-weighted-delta},
\begin{align*}
\mathbb E[m_{1,t}(Y_1)^2]
&\le
C \frac{p^2}{\delta_p^2} e^{-c t^2}
\mathbb E\!\left[\frac{Y_1^2}{(1-Y_1^2)^8}\right]
+
C \frac{p^4}{\delta_p^4(1+|t|)^4}
\mathbb E\!\left[\frac{Y_1^6}{(1-Y_1^2)^8}\right] \\
&\le
C \frac{p^2}{\delta_p^3} e^{-c t^2}
+
C \frac{p^4}{\delta_p^7}(1+|t|)^{-4}.
\end{align*}
Since \(\delta_p\gg p^{1/2}\), we have for all large \(n\),
\[
\frac{p^4}{\delta_p^7}
\lesssim
\frac{p^2}{\delta_p^3}.
\]
This proves \eqref{m-estimate-tail}. The proof of \eqref{eq:Bnprime-tail-var-stronger}, and therefore of \eqref{eq:Bnprime-tail-var}, is now complete. $\hfill$ $\square$

\subsection{Proof of \eqref{mcA}}

To show \eqref{mcA}, it suffices to show the following two statements:
\begin{align} \label{mcA1}
    \sup_{|t| \leq a} |\mc{A}_n(t)| = o_{\mb{P}}(1) \qquad \text{for any fixed $a>0$},
\end{align}
and 
\begin{align} \label{mcA2}
     \qquad \lim_{a \to \infty}\la \limsup_{n \to \infty} \mb{P} \lb \sup_{|t|>a} |\mc{A}_n(t)|>\ve \rb \ra =0 \qquad \text{for all $\ve>0$}.
\end{align}
\noindent \underline{\textbf{Proof of \eqref{mcA1}.}} We will first rewrite the supremum in the form of the supremum of a degenerate $U$-process.
To this end, define the random elements $\bm{\mc{W}}_i=\lb Y_i, \bm{U}_i \rb$ taking value in $\mc{X}_1=[-1,1]\times \mb{S}^{p-2}$. Let $K_t:\mc{X}_1 \times \mc{X}_1 \to \mb{R}$ be the symmetric kernels defined by 
\[
K_t\lb \bm{a}, \bm{b} \rb := \mathbf 1_{\la \sqrt p \Big( \,y_{\bm{a}}y_{\bm{b}} +\sqrt{\lb 1 - y^2_{\bm{a}} \rb \lb 1-y_{\bm{b}}^2 \rb}\, \bm{u}_{\bm{a}}^\top \bm{u}_{\bm{b}} \Big) \le t \ra}
\]
where $\bm{a}=\lb y_{\bm{a}}, \bm{u}_{\bm{a}} \rb$ and  $\bm{b}=\lb y_{\bm{b}}, \bm{u}_{\bm{b}} \rb$.

Put 
$$\tilde{K}_t(\bm{a}, \bm{b}):= K_t\lb \bm{a}, \bm{b} \rb - \mb{E} \lb K_t\lb \lb y_{\bm{a}}, \bm{U}_{\bm{a}} \rb, \lb y_{\bm{b}}, \bm{U}_{\bm{b}} \rb \rb \rb$$
where the expectation is taken over $\bm{U}_a, \bm{U}_b$ uniformly distributed on $\mb{S}^{p-2}$.

Now, observe that $\mc{A}_n$ in \eqref{mcA} can be rewritten as
\begin{align*}
    \mc{A}_n(t)= N_n^{-1/2} \sum_{1\leq i<j \leq n} \tilde{K}_t \lb \bm{\mc{W}}_i, \bm{\mc{W}}_j   \rb 
\end{align*}
Moreover, $\mb{E} \lb \tilde{K}_t \lb \bm{\mc{W}}_1, \bm{\mc{W}}_2   \rb |   \bm{\mc{W}}_1 \rb=\mb{E} \lb \tilde{K}_t \lb \bm{\mc{W}}_1, \bm{\mc{W}}_2   \rb |   \bm{\mc{W}}_2 \rb=0$ almost surely. To see this, fix \(t\in \mb{R}\) and condition on \(\bm{\mc W}_1=(Y_1,\bm U_1)\) and \(Y_2\). Observe that
\[
K_t(\bm{\mc W}_1,\bm{\mc W}_2)
=
\mathbf 1_{\left\{
\sqrt p \Big(
Y_1Y_2+\sqrt{(1-Y_1^2)(1-Y_2^2)}\,\bm U_1^\top \bm U_2
\Big)\le t
\right\}}.
\]
Given \((Y_1,\bm U_1,Y_2)\), the only remaining randomness comes from \(\bm U_2\), which is uniformly distributed on \(\mb S^{p-2}\). By rotational invariance of the uniform distribution on \(\mb S^{p-2}\), the law of
$\bm U_1^\top \bm U_2$
does not depend on the fixed vector \(\bm U_1\); indeed, it is the same as the law of
$\bm e_1^\top \bm U_2$,
and hence also the same as the law of \(\bm U_{\bm a}^\top \bm U_{\bm b}\), where \(\bm U_{\bm a},\bm U_{\bm b}\) are independent and uniformly distributed on \(\mb S^{p-2}\). 

Consequently,
\[
\mb E\!\left[K_t(\bm{\mc W}_1,\bm{\mc W}_2)\mid \bm{\mc W}_1,Y_2\right]
=
\mb E\!\left[
K_t\bigl((Y_1,\bm U_{\bm a}),(Y_2,\bm U_{\bm b})\bigr)
\right]
\]
where the expectation on the right-hand side is taken over \(\bm U_{\bm a},\bm U_{\bm b}\).

Therefore,
\[
\mb E\!\left[\tilde K_t(\bm{\mc W}_1,\bm{\mc W}_2)\mid \bm{\mc W}_1,Y_2\right]=0,
\]
and by taking conditional expectation once more,
\[
\mb E\!\left[\tilde K_t(\bm{\mc W}_1,\bm{\mc W}_2)\mid \bm{\mc W}_1\right]=0.
\]
By symmetry, the same argument also yields
\[
\mb E\!\left[\tilde K_t(\bm{\mc W}_1,\bm{\mc W}_2)\mid \bm{\mc W}_2\right]=0.
\]
Hence the kernel \(\tilde K_t\) is degenerate.

Obviously, $\tilde{K}_t$'s are also symmetric as well. By Lemma 24 in \cite{belloni2019conditional}, the class of functions $\la \tilde{K}_t; t \in [-a,a] \ra$ has its entropy number grows at most polynomially fast in terms of its radius. Therefore, Lemma \ref{U-process concentration} gives 
\begin{align*}
    \mb{E} \left[ \sup_{|t| \leq a} |\mc{A}_n(t)| \right] \lesssim \sigma_{n,a} \lb 1 + \log\lb \frac{1}{\sigma_{n,a}} \rb \rb +\frac{1}{\sqrt{n}} \lb 1 + \log^2 \lb \frac{1}{\sigma_{n,a}} \rb \rb
\end{align*}
where 
$$\sigma_{n,a}^2:= \sup_{|t|\leq a} \mbox{Var} \lb \tilde{K}_t \lb \bm{\mc{W}}_1,\bm{\mc{W}}_2 \rb \rb = \sup_{|t|\leq a} \mbox{Var} \lb \mc{A}_{n}(t) \rb \leq C_{a,\tau} \cdot \max \la p^{-1/2}, n^{-1/2} \ra.$$
Consequently,
\[
\lim_{n \to \infty} \mb{E} \left[ \sup_{|t| \leq a} |\mc{A}_n(t)| \right] = 0 
\]
for every fixed $a>0$ due to Lemma \ref{upper bound variance An} and Lemma \ref{lower bound variance An}. Note that we have used the assumption $\log(p)=o(n^{1/4})$ to deduce that $\log \lb 1/\sigma_{n,a} \rb^2 \ll \sqrt{n}$. This concludes the proof of \eqref{mcA1}. $\hfill$ $\square$

\noindent \underline{\textbf{Proof of \eqref{mcA2}.}} We use the same trick as in Step 2 of the proof of \eqref{mcA1}. Fix $\ve>0$. Using the same degenerate U-process representation, it suffices to estimate the variance profile 
\[
\tilde{\sigma}^2_{n,a}:= \sup_{|t|> a} \mbox{Var} \lb \tilde{K}_t \lb \bm{\mc{W}}_1,\bm{\mc{W}}_2 \rb \rb = \sup_{|t| > a} \mbox{Var} \lb \mc{A}_{n,t} \rb.
\]
To this end, we will introduce the high probability event $G_n= \la \max_{1 \leq i \leq 2} |Y_i| \leq \gamma_n  \ra$ for a sequence of thresholds $\la \gamma_n \ra$ to be determined later.  From \eqref{eq:An-second-moment}, we have 
\begin{align*}
\sup_{|t| > a} \mbox{Var} \lb \mc{A}_{n,t} \rb &\leq \sup_{|t| > a} \mb{P} \lb D_{12}(t) \neq 0 \rb \\
&\leq \sup_{|t| > a} \mb{P} \lb D_{12}(t) \neq 0, G_n \rb + \mb{P} \lb G_n^c \rb.
\end{align*}
Recall that $D_{12}(t)= \mathbf 1_{\{\sqrt p(\beta_{12}+\alpha_{12}R_{12})\le t\}}
-
\mathbf 1_{\{\sqrt p R_{12}\le t\}}$. Suppose $\gamma_n$ is chosen such that $\sqrt{p} \gamma_n^2 <1$, then for $a>1$, we have 
\begin{align*}
|\sqrt p(\beta_{12}+\alpha_{12}R_{12})| &= \left|\sqrt{p}Y_1Y_2 + \sqrt{p}\sqrt{(1-Y_1^2)(1-Y_2^2)}R_{12}\right| \\
&\leq \sqrt{p}\gamma_n^2 + \sqrt{p}|R_{12}| \\
&\leq 1 +  \sqrt{p}|R_{12}|.
\end{align*}
If $D_{12}(t) \neq 0$, then on $G_n$, at least one of $|\sqrt p(\beta_{12}+\alpha_{12}R_{12})|$ or $\sqrt{p}\,|R_{12}|$ has to exceed $a$. If $|\sqrt p(\beta_{12}+\alpha_{12}R_{12})|>a$, then $\sqrt{p}\,|R_{12}|>a-1$, and we have $\sqrt{p}\,|R_{12}|>a$ otherwise.  Therefore,
\begin{align*}
     \sup_{|t| > a} \mb{P} \lb D_{12}(t) \neq 0, G_n \rb + \mb{P} \lb G_n^c \rb \leq \mb{P} \lb \sqrt{p}\,|R_{12}| > a-1 \rb + \mb{P} \lb G_n^c \rb.
\end{align*}
Pick $\gamma_n:=p^{-1/4}(\log p)^{-1}$,
we have $\sqrt p\,\gamma_n^2=(\log p)^{-2} \ll 1$
for all sufficiently large \(n\). Recall $\delta_p$ from \eqref{delta-p}. By Markov's inequality and \eqref{eq:simple-weighted-delta},
\[
\mb P(G_n^c)
\le 2\,\mb P(|Y_1|>\gamma_n)
\le \frac{2\mb E(Y_1^2)}{\gamma_n^2}
\lesssim \frac{\delta_p^{-1}}{p^{-1/2}(\log p)^{-2}}
= \frac{\sqrt{p}\,\log^2 p}{\delta_p}
=o(1).
\]
since $\delta_p=p/2-\kappa \gg p^{1/2}n^{1/4}$ and $\log(p) \ll n^{1/4}$. Hence, for every fixed \(a>1\),
\[
\tilde{\sigma}_{n,a}^2
\le \mb P\!\left(\sqrt p\,|R_{12}|>a-1\right)+o(1).
\]
Also, we have the lower bound 
\[
\tilde{\sigma}_{n,a}^2 \geq \Var \lb \mc{A}_n(a+1) \rb \gtrsim p^{-1/2}.
\]
Therefore, Lemma \ref{U-process concentration} yields
\begin{align*}
    \mb{E} \left[ \sup_{|t| > a} |\mc{A}_n(t)| \right] \lesssim \tilde{\sigma}_{n,a} \lb 1 + \log\lb \frac{1}{ \tilde{\sigma}_{n,a}} \rb \rb +\frac{1}{\sqrt{n}} \lb 1 + \log^2 \lb \frac{1}{\tilde{\sigma}_{n,a}} \rb \rb. 
\end{align*}
The proof is completed by first sending $n \to \infty$ and then sending $a \to \infty$. Note that we use Lemma \ref{upper bound variance An} and \ref{lower bound variance An} to guarantee that the variance term $\tilde{\sigma}_{n,a}$ goes to zero, but does not decay too quickly. $\hfill$ $\square$

\subsection{Proof of \eqref{mcB}}

To show \eqref{mcB}, it suffices to show the following two statements:
\begin{align} \label{mcB1}
    \sup_{|t| \leq a} |\mc{B}_n(t)| = o_{\mb{P}}(1) \qquad \text{for any fixed $a>0$},
\end{align}
and 
\begin{align} \label{mcB2}
     \qquad \lim_{a \to \infty}\la \limsup_{n \to \infty} \mb{P} \lb \sup_{|t|>a} |\mc{B}_n(t)|>\ve \rb \ra =0 \qquad \text{for all $\ve>0$}.
\end{align}

\noindent \textbf{\underline{Proof of \eqref{mcB1}.}} 
 Observe that by the fundamental theorem of calculus,
\begin{align*}
\sup_{|t| \leq a} |\mc{B}_n(t)| \leq |\mc{B}_n(0)|+ \left| \int_{-a}^{a} \mc{B}'_n(u)\, du   \right| \leq \int_{-a}^a \left| \mc{B}'_n(u) \right|\,du.
\end{align*}
Taking expectation in both sides and using Cauchy--Schwartz inequality, we get
\[
\mb{E} \sup_{|t| \leq a} |\mc{B}_n(t)| \leq \mb{E}|\mc{B}_n(0)|+ \int_{-a}^a \mb{E} \lb \left| \mc{B}'_n(u) \right| \rb\, du \leq \sqrt{\Var \lb \mc{B}_n(0) \rb} +  \int_{-a}^{a} \sqrt{\Var \lb \mc{B}'_n(u) \rb}\, du. 
\]
By using Lemma \ref{variance bound Brascamp-Lieb}, we get $\mb{E} \sup_{|t| \leq a} |\mc{B}_n(t)|=o(1)$ and this finishes the proof of \eqref{mcB1}.

\noindent \textbf{\underline{Proof of \eqref{mcB2}.}} 
As in the proof of \eqref{mcB1}, we estimate 
\begin{align*}
    \sup_{|t| > a} |\mc{B}_n(t)| \leq \sup_{t < -a} |\mc{B}_n(t)| + \sup_{t > a} |\mc{B}_n(t)| &\leq 2|\mc{B}_n(0)|+\int_{-\infty}^{-a} \left| \mc{B}'_n(t) \right|\,dt + \int_{a}^{\infty}  \left| \mc{B}'_n(t) \right|\,dt \\
    &= 2|\mc{B}_n(0)|+ \int_{\mb{R}\setminus [-a,a]} \left| \mc{B}'_n(t) \right|\,dt.
\end{align*}
Taking expectation in both sides, using Lemma \ref{tail variance Bn(t)} and Cauchy--Schwartz inequality, we get
\begin{align*}
\mb{E} \sup_{|t| > a} |\mc{B}_n(t)| &\leq 2\mb{E}|\mc{B}_n(0)|+ \int_{\mb{R}\setminus[-a,a]} \mb{E} \lb \left| \mc{B}'_n(u) \right| \rb\, du\\
&\leq  2\sqrt{\Var \lb \mc{B}_n(0) \rb}+ \int_{\mb{R}\setminus[-a,a]} \sqrt{\Var \lb \mc{B}'_n(u) \rb}\, du \\
&\lesssim  \int_{\mb{R}\setminus[-a,a]} \left[ \exp \lb -\frac{ct^2}{2} \rb + \frac{1}{(1+|t|)^2} \right] dt +o(1)
\end{align*}
where $c$ is the universal constant in Lemma \ref{tail variance Bn(t)}. 

Consequently,
\[
\mb{P} \lb   \sup_{|t| > a} |\mc{B}_n(t)| > \ve \rb \lesssim \frac{1}{\ve} \cdot    \int_{\mb{R}\setminus[-a,a]} \left[ \exp \lb -\frac{ct^2}{2} \rb + \frac{1}{(1+|t|)^2} \right] dt.
\]
The proof is completed by first letting $n \to \infty$ and then letting $a \to \infty$.

\section{Proof of Theorem \ref{inconsistency}}

We first give some quantitative estimates on the geometry of pairwise inner products. Let \( \bm{x}\in C_{r,p}\) and \( \bm{y}\in C_{s,p}\) with \(r\neq s\). Denote by
\[
\mathrm{dist}(\bm{a},\bm{b}):=\arccos(\bm{a}^\top \bm{b}),\qquad \bm{a},\bm{b}\in \mathbb S^{p-1},
\]
the geodesic distance on the sphere. Since \(\mathrm{dist}(\bm{v}_r,\bm{v}_s)=\arccos(-1/p)\), the triangle inequality gives
\[
\arccos(-1/p)-2\varepsilon_p
\le
\mathrm{dist}(\bm{x},\bm{y})
\le
\arccos(-1/p)+2\varepsilon_p.
\]
Observe that \(u\mapsto \cos u\) is \(1\)-Lipschitz on \([0,\pi]\), hence
\[
\left|\bm{x}^\top \bm{y}+\frac1p\right|
=
\left|\cos(\mathrm{dist}(\bm{x},\bm{y}))-\cos(\arccos(-1/p))\right|
\le 2\varepsilon_p.
\]
Therefore, for all large \(p\),
\[
\frac{-3}{2p}=-\frac1p-2\varepsilon_p
\le
\bm{x}^\top \bm{y}
\le
-\frac1p+2\varepsilon_p=
-\frac{1}{2p}.
\]
In particular,
\[
|\bm{x}^\top \bm{y}|
\le
\frac{3}{2p}.
\]
On the other hand, if \(\bm{x},\bm{y}\in C_{r,p}\) belong to the same cap, then
\[
\mathrm{dist}(\bm{x},\bm{y})\le 2\varepsilon_p,
\]
so
\[
\bm{x}^\top \bm{y}\ge \cos(2\varepsilon_p)>0.
\]
We now prove the first statement. Let \(\bm{X},\bm{Y}\stackrel{i.i.d.}{\sim}\mu_n\). If the cap labels of \(\bm{X}\) and \(\bm{Y}\) are different, then \(\bm{X}^\top \bm{Y}<0\); if the labels are the same, then \(\bm{X}^\top \bm{Y}>0\). Since the labels are i.i.d. uniform on \(\{1,\dots,p+1\}\),
\[
\mathbb P_{\mu_n}(\bm{X}^\top \bm{Y} \le 0)=\frac{p}{p+1}.
\]
Under \(\mu_0\), suppose $\bm{U}, \bm{V}$ are uniformly distributed. By symmetry,
\[
\mathbb P_{\mu_0}(\bm{U}^\top \bm{V} \le 0)=\frac12.
\]
Hence
\[
d(\mu_n,\mu_0)
\ge
\left|
\mathbb P_{\mu_n}(\bm{X}^\top \bm{Y} \le 0)
-
\mathbb P_{\mu_0}(\bm{U}^\top \bm{V} \le 0)
\right|
=
\frac{p}{p+1}-\frac12,
\]
which proves the first statement.

Regarding the other three statements, consider the sample \( \bm{X}_1,\dots,\bm{X}_n\) and let \(L_i\in\{1,\dots,p+1\}\) denote the label of the cap from which \(\bm{X}_i\) is drawn. Define the event
\[
E_n:=\{L_1,\dots,L_n\ \text{are all distinct}\}.
\]
Since the labels are i.i.d. uniform on \(\{1,\dots,p+1\}\),
\[
\mathbb P(E_n^c)
\le
\binom{n}{2} \cdot \frac{1}{p+1}
=
o(1)
\]
because \(p/n^2\to\infty\). Thus \(E_n\) occurs with probability \(1-o(1)\).

On \(E_n\), every pair \(\bm{X}_i,\bm{X}_j\) comes from different caps, and therefore for all \(i<j\),
\[
-\frac{3}{2p}\le \bm{X}_i^\top \bm{X}_j\le -\frac{1}{2p}.
\]
We now prove the second statement. On \(E_n\),
\[
|R_n|
\le
\frac{\sqrt{2p}}{n} \cdot \binom{n}{2} \cdot \frac{3}{2p}
\le
\frac{3(n-1)}{2\sqrt{2p}}
\to 0.
\]
Since \(\mathbb P(E_n)\to 1\), it follows that \(R_n\to 0\) in probability. This proves the second statement.

Let us prove the third statement. On \(E_n\),
\[
(\bm{X}_i^\top \bm{X}_j)^2\le \frac{9}{4p^2},
\]
hence
\[
B_n
\le
\frac{p}{n} \cdot \binom{n}{2} \cdot \left(\frac{9}{4p^2}-\frac1p\right)
=
\frac{n-1}{2}\left(\frac{9}{4p}-1\right)
\to -\infty.
\]
Therefore \(B_n\to -\infty\) in probability, which proves the third statement.

Finally, again on \(E_n\),
\[
\max_{1\le i<j\le n}(\bm{X}_i^\top \bm{X}_j)^2\le \frac{9}{4p^2},
\]
so
\[
P_n
\le
\frac{9}{4p}-4\log n+\log\log n
\to -\infty.
\]
Thus \(P_n\to -\infty\) in probability, proving the last statement. The proof is completed. $\hfill$ $\square$


\section{Proof of Proposition \ref{identifiable}} \label{proof-identifiable}
Before presenting the proof, we first state a version of Lebesgue's differentiation theorem for smooth, complete Riemannian manifolds. The version provided below is not the most general result, but it is sufficient for our purposes.

\begin{lemma} \label{Lebesgue-diff}
    Let $(M,g)$ be a smooth, complete Riemannian manifold with the corresponding geodesic distance $d$. Suppose $\nu$ is a non-negative, finite Borel measure on the metric space $(M,d)$ and $f$ is a non-negative integrable function with respect to $\nu$. Then, we have
    \begin{align} \label{Lebesgue-diff2}
        \lim_{r \to 0} \frac{\int_{B(x,r)} f(y) d \nu(y)}{\nu \lb B(x,r) \rb} = f(x)
    \end{align}
    for $\nu$-almost everywhere $x \in M$. Here $B(x,r)$ is the open ball with respect to the geodesic distance $d$.
\end{lemma}

{\noindent \textbf{Proof of Lemma \ref{Lebesgue-diff}}.} Define 
\begin{align*}
    \nu_1(A) = \int_A f(x) d\nu(x).
\end{align*}
It is easy to see that $\nu_1 \ll \nu$. By Theorem A.1 in \cite{jost2021probabilistic}, there exists a measurable set $S_0 \subset M$ such that $\nu(S_0)=0$ and 
\begin{align} \label{D(x)}
    D(x)= \lim_{r \to 0} \frac{\nu_1 \lb B(x,r) \rb}{\nu \lb B(x,r) \rb}
\end{align}
exists and is finite for all $x \in M \setminus S_0$. Also by Theorem A.1 in \cite{jost2021probabilistic}, $D(x)$ is the Radon--Nikodym derivative $d \nu_1/d\nu$ (up to a null set) whenever $(M,d)$ is complete. Thus, $D(x)=f(x)$ $\nu$-almost everywhere. The proof is complete.
$\hfill$ $\square$

In the argument of Lemma \ref{Lebesgue-diff}, the completeness of $(M,d)$ is needed only to deduce that $D(x)$ in \eqref{D(x)} is equal to the Radon--Nikodym derivative $f(x)$ $\nu$-almost everywhere. Results of this type hold for various metric spaces with different structures; see the classical monograph \cite{GMT} for more details.

\noindent \textbf{Proof of Proposition \ref{identifiable}.} It is easy to see that the assumption \eqref{characterization} is equivalent to
\begin{align} \label{identity2}
    \int_{\mb{S}^{p-1}} \int_{\mb{S}^{p-1}} g(\bm{x}^{\top} \bm{y}) d\nu(\bm{x}) d\nu(\bm{y}) =  \int_{\mb{S}^{p-1}} \int_{\mb{S}^{p-1}} g(\bm{x}^{\top} \bm{y}) d\mu_0(\bm{x}) d\mu_0(\bm{y})
\end{align}
for any bounded, measurable function $g: [-1,1] \to \mb{R}$. We will show that (\ref{identity2}) implies $\nu \equiv \mu_0$. To see this, fix $\eta \in (0,2]$ and define 
$$g_{\eta}(t) :=  \frac{ \mathbf{1}_{(1-\eta,1]}(t)}{\int_{\mb{S}^{p-1}} \int_{\mb{S}^{p-1}} \mathbf{1}_{(1-\eta,1]}(\bm{x}^{\top} \bm{y}) d\mu_0(\bm{x}) d\mu_0(\bm{y})}.$$
Let $\mu_1 := (\nu+\mu_0)/2$, and define
\[
f := \frac{d\nu}{d\mu_1}, \qquad h := \frac{d\mu_0}{d\mu_1}.
\]
Then
\[
d\nu = f\,d\mu_1,\qquad d\mu_0 = h\,d\mu_1,
\]
so
\[
2\,d\mu_1 = d\nu + d\mu_0 = (f+h)\,d\mu_1.
\]
Hence
\[
f+h=2
\]
$\mu_1$-almost surely. Plugging $g_{\eta}$ into (\ref{identity2}) gives
\begin{align*} 
   \frac{  \int_{\mb{S}^{p-1}} \int_{\mb{S}^{p-1}} \mathbf{1}_{(1-\eta,1]}( \bm{x}^{\top} \bm{y}) d\nu(\bm{x}) d\nu(\bm{y})}{\int_{\mb{S}^{p-1}} \int_{\mb{S}^{p-1}}  \mathbf{1}_{(1-\eta,1]}( \bm{x}^{\top} \bm{y}) d\mu_0(\bm{x}) d\mu_0(\bm{y})}=  \underbrace{\frac{\int_{\mb{S}^{p-1}} \int_{\mb{S}^{p-1}}  \mathbf{1}_{(1-\eta,1]}( \bm{x}^{\top} \bm{y}) d\mu_0(\bm{x}) d\mu_0(\bm{y})}{\int_{\mb{S}^{p-1}} \int_{\mb{S}^{p-1}}  \mathbf{1}_{(1-\eta,1]}( \bm{x}^{\top} \bm{y}) d\mu_0(\bm{x}) d\mu_0(\bm{y})}}_{=1}.
   \end{align*}
Thus,
\begin{align} \label{fact2}
    \int_{\mb{S}^{p-1}} \int_{\mb{S}^{p-1}} \mathbf{1}_{(1-\eta,1]}( \bm{x}^{\top} \bm{y}) d\nu(\bm{x}) d\nu(\bm{y}) = \int_{\mb{S}^{p-1}} \int_{\mb{S}^{p-1}}  \mathbf{1}_{(1-\eta,1]}( \bm{x}^{\top} \bm{y}) d\mu_0(\bm{x}) d\mu_0(\bm{y}).
\end{align}
Moreover, by Fubini's theorem,
\begin{align*}
    & \int_{\mb{S}^{p-1}} \int_{\mb{S}^{p-1}} \mathbf{1}_{(1-\eta,1]}( \bm{x}^{\top} \bm{y}) d\nu(\bm{x}) d\nu(\bm{y}) \\
    =& \int_{\mb{S}^{p-1}} \int_{\mb{S}^{p-1}} \mathbf{1}_{(1-\eta,1]}( \bm{x}^{\top} \bm{y}) f(\bm{x})f(\bm{y}) d\mu_1(\bm{x}) d\mu_1(\bm{y}) \\
    = &  \int_{\mb{S}^{p-1}} f(\bm{x}) \lb \int_{1 -\eta < \bm{x}^{\top} \bm{y} \leq 1} f(\bm{y}) d\mu_1(\bm{y}) \rb d \mu_1(\bm{x}).
\end{align*}
Additionally, 
\begin{align*}
    \int_{\mb{S}^{p-1}} \int_{\mb{S}^{p-1}}  \mathbf{1}_{(1-\eta,1]}( \bm{x}^{\top} \bm{y}) d\mu_0(\bm{x}) d\mu_0(\bm{y}) &= \mb{P}_{\mu_0} \lb 1-\eta < \bm{x}^{\top} \bm{y} \leq 1  \rb.
\end{align*}
Therefore, from \eqref{fact2} and the two equalities above, we get
\begin{align} \label{fact3}
     \int_{\mb{S}^{p-1}} f(\bm{x}) \lb \int_{1 -\eta < \bm{x}^{\top} \bm{y} \leq 1} f(\bm{y}) d\mu_1(\bm{y}) \rb d \mu_1(\bm{x}) = \mb{P}_{\mu_0} \lb 1-\eta < \bm{x}^{\top} \bm{y} \leq 1  \rb
\end{align}
for all $\eta \in (0,2]$. 

Let $d(\bm{x},\bm{y})=\arccos\lb \bm{x}^{\top} \bm{y} \rb$ be the geodesic distance on $\mb{S}^{p-1}$. It is easy to check that $\lb \mb{S}^{p-1},d \rb$ is a Polish space and that for all $\bm{x}\in \mb{S}^{p-1}$,
\[
\la \bm{y}: 1-\eta < \bm{x}^{\top} \bm{y} \leq 1  \ra = B(\bm{x},f_{\eta})
\]
where $f_{\eta}=\arccos(1-\eta)$ and $B(\bm{x},r)$ is the open ball with center $\bm{x}$ and radius $r$ with respect to $d$. Hence, for all $\eta \in (0,2]$, one can rewrite \eqref{fact3} as
\begin{align} \label{Holder}
    \int_{\mb{S}^{p-1}} f(\bm{x}) \lb  \frac{\int_{B(\bm{x},f_{\eta})} f(\bm{y}) d\mu_1(\bm{y})}{{\int_{B(\bm{x},f_{\eta})} h(\bm{y}) d\mu_1(\bm{y})}} \rb d\mu_1(\bm{x}) =1
\end{align}
where we have used Lemma \ref{degenerate-Ustat} to write 
$$\mb{P}_{\mu_0} \lb 1-\eta < \bm{x}^{\top} \bm{y} \leq  1 \rb =\mu_0 \lb B(\bm{x}, f_{\eta}) \rb = \int_{B(\bm{x},f_{\eta})} h(\bm{y}) d\mu_1(\bm{y}) $$
for all $\bm{x} \in \mb{S}^{p-1}$. Note that \eqref{Holder} holds since the right-hand side in the expression above is constant across $\bm{x}$. 

Since $\lb \mb{S}^{p-1}, d \rb$ is a smooth, complete Riemannian manifold with respect to the canonical Riemannian metric, Lemma \ref{Lebesgue-diff} can be applied to $\mu_1$ to deduce that 
\begin{align*}
    \frac{\int_{B(\bm{x},f_{\eta})} f(\bm{y}) d\mu_1(\bm{y})}{{\int_{B(\bm{x},f_{\eta})} h(\bm{y}) d\mu_1(\bm{y})}} =  \frac{\int_{B(\bm{x},f_{\eta})} f(\bm{y}) d\mu_1(\bm{y}) }{ \mu_1 \lb B(\bm{x}, f_{\eta}) \rb} \cdot \lb \frac{\int_{B(\bm{x},f_{\eta})} h(\bm{y}) d\mu_1(\bm{y}) }{ \mu_1 \lb B(\bm{x}, f_{\eta}) \rb} \rb^{-1} \to \frac{f(\bm{x})}{h(\bm{x})},
\end{align*}
as $\eta \to 0$, for $\mu_1$-almost every $\bm{x}$, since $f_{\eta} \to 0$ as $\eta \to 0$. Therefore, the above display together with (\ref{Holder}) and Fatou's lemma yields
\begin{align*}
    \int_{\mb{S}^{p-1}} \frac{f(\bm{x})^2}{2-f(\bm{x})} d\mu_1(\bm{x})  &= \int_{\mb{S}^{p-1}} \frac{f(\bm{x})^2}{h(\bm{x})} d\mu_1(\bm{x}) \\
    &= \int_{\mb{S}^{p-1}} \lim_{\eta \to 0} f(\bm{x}) \lb  \frac{\int_{B(\bm{x},f_{\eta})} f(\bm{y}) d\mu_1(\bm{y})}{{\int_{B(\bm{x},f_{\eta})} h(\bm{y}) d\mu_1(\bm{y})}} \rb d\mu_1(\bm{x})  \\
    &\leq \liminf_{\eta \to 0}   \int_{\mb{S}^{p-1}} f(\bm{x}) \lb  \frac{\int_{B(\bm{x},f_{\eta})} f(\bm{y}) d\mu_1(\bm{y})}{{\int_{B(\bm{x},f_{\eta})} h(\bm{y}) d\mu_1(\bm{y})}} \rb d\mu_1(\bm{x})= 1
\end{align*}
due to \eqref{Holder}. Moreover, Hölder's inequality gives
\begin{align*}
   \int_{\mb{S}^{p-1}} \frac{f(\bm{x})^2}{2-f(\bm{x})} d\mu_1(\bm{x}) &= \int_{\mb{S}^{p-1}} \frac{f(\bm{x})^2}{2-f(\bm{x})} d\mu_1(\bm{x}) \cdot \int_{\mb{S}^{p-1}} \left[ 2-f(\bm{x}) \right]d\mu_1(\bm{x}) \\
   & \geq \lb  \int_{\mb{S}^{p-1}} f(\bm{x}) d\mu_1(\bm{x}) \rb^2 =1.
\end{align*}
The two bounds above imply that the integral $\int_{\mb{S}^{p-1}} f(\bm{x})^2 \cdot \lb 2-f(\bm{x})\rb^{-1} d\mu_1(\bm{x})$ is exactly $1$ and $f(\bm{x})^2=(2-f(\bm{x}))^2$ for $\mu_1$-almost every $\bm{x}$. This in turn yields $f \equiv 1$ almost surely with respect to $\mu_1$. Therefore, we also get $h\equiv1$, and the conclusion follows. $\hfill$ $\square$


\section{Proof of Proposition \ref{watson d asymptotic}}

Let us start by explaining the idea. Suppose $\bm{U}, \bm{V}$ are i.i.d. uniformly distributed on $\mb{S}^{p-2}$ and put $H_{p}(x) := \mb{P} \lb \sqrt{p} \bm{U}^\top \bm{V} \leq x \rb$. By using either the Edgeworth expansion \cite{Hall2} or the Student density ratio expansion \cite{ouimet2022refined}, we get 
\begin{align} \label{edgeworth1}
    H_p(x)  = \Phi(x) + \frac{\Psi(x)}{p} + O(p^{-2})
\end{align}
uniformly in $x \in \mb{R}$. Here $\Phi$ is the CDF of the standard normal, and $\Psi$ is a smooth function such that 
\[
\sup_{x \in \mb{R}} (1+x^{2k}) |\Psi^{(l}(x)| < \infty
\]
for all $k,l \geq 1$. In fact, $\Psi$ is the product of the standard Gaussian density and a polynomial of degree no greater than three. The specific form of $\Psi$ is not relevant for the proof below. 

By rotational invariance, we may assume that the location parameter is $\bm{e}_1$. For $\bm{X},\bm{Y}$ having the Watson distribution with location $\bm{e}_1$ and concentration $\kappa$, decompose 
\begin{align*}
   \sqrt{p} \bm{X}^\top \bm{Y} = A + (1+B) \sqrt{p}\,\bm{U}^\top \bm{V}
\end{align*}

where 
\begin{align} \label{A,B}
T:&=\bm{e}_1^\top \bm{X};\qquad T':=\bm{e}_1^\top \bm{Y} \nonumber; \\
    A:&= \sqrt{p}TT^{'}; \qquad B:= \sqrt{(1-T^2)(1-{T'}^2)}-1.
\end{align}
Similarly, under the null (which is a Watson distribution with $\kappa=0$), let $\bm{X}_0, \bm{Y}_0$ be i.i.d. from $\mu_0$ and define analogously
\[
T_0:=\bm{e}_1^\top \bm{X}_0,\qquad T_0':=\bm{e}_1^\top \bm{Y}_0,
\]
\[
   \sqrt{p} \bm{X}_0^\top \bm{Y}_0 = A_0 + (1+B_0) \sqrt{p}\,\bm{U}_0^\top \bm{V}_0,
\]
\[
A_0:=\sqrt{p}\,T_0T_0',
\qquad
B_0:=\sqrt{(1-T_0^2)(1-{T_0'}^2)}-1.
\]
In the argument below, we assume all random variables are defined in the same probability space. Observe that
\begin{align*}
    d \lb \mu_n, \mu_0 \rb &= \sup_{u \in \mb{R}} \left| \mb{P}_{\mu_n} \lb  A + (1+B) \sqrt{p}\,\bm{U}^\top \bm{V} \leq u  \rb -  \mb{P}_{\mu_0} \lb  A_0 + (1+B_0) \sqrt{p}\,\bm{U}_0^\top \bm{V}_0 \leq u  \rb   \right| \\
    &= \sup_{u \in \mb{R}} \left| \mb{E} \left[ H_p \lb \frac{u-A}{1+B} \rb - H_p \lb \frac{u-A_0}{1+B_0} \rb \right] \right| \\
    &= \sup_{u \in \mb{R}} \left| \mb{E} \left[ \Phi \lb \frac{u-A}{1+B} \rb - \Phi \lb \frac{u-A_0}{1+B_0} \rb +  \frac{1}{p}\Psi \lb \frac{u-A}{1+B} \rb - \frac{1}{p}\Psi \lb \frac{u-A_0}{1+B_0} \rb \right]  \right| + O(p^{-2}). 
\end{align*}
One can see right away that the conclusion follows if one can show that $A,B$ and $A_0,B_0$ are close in the first few moments and that moments of the same orders are sufficiently close, since $\Phi$ and $\Psi$ have derivatives that decay faster than any fixed-degree polynomial. 

We will need the following asymptotic formulas for the moments of $A,B$ up to first order. 
\begin{lemma} \label{moment asymptotic}
    Recall $A,B$ in \eqref{A,B} and their corresponding null version $A_0, B_0$. Define
   \begin{align}
    \delta_p :&= p/2-\kappa \label{delta-p};  \\
    r_p :&= \kappa/\delta_p  \label{r-p}.
\end{align}
    Then,
    \begin{itemize}
        \item All odd moments of $A, A_0$ vanish and 
        \begin{align*}
            \mb{E} A^2 &= \frac{(1+r_p)^2}{p} + O \lb \frac{1}{p^2}+\frac{1}{n^2} \rb; \qquad \mb{E} A_0^2 = \frac{1}{p}; \\
            \mb{E} |A_0|^3 + \mb{E} |A|^3  &= O \lb \frac{1}{p^{3/2}} + \frac{1}{n^{3/2}} \rb.
        \end{align*}
        \item $B, B_0$ take values in $[-1,0]$ and 
        \begin{align*}
        \mb{E} B &= -\frac{1+r_p}{p} + O \lb \frac{1}{p^2} + \frac{1}{n\sqrt{pn}} \rb; \qquad \mb{E}B_0= -\frac{1}{p} + O \lb \frac{1}{p^2} \rb;    \\
            \mb{E} B^2 &=O \lb \frac{1}{p^2} + \frac{1}{pn} \rb; \qquad  \mb{E} B_0^2 =O(p^{-2}).
        \end{align*}
    \end{itemize}
\end{lemma}

The proof of Lemma \ref{moment asymptotic} is given in Section \ref{moment compute}. The following technical Taylor expansion will also be used.

\begin{lemma}\label{Taylor}
For a function \(g\in C^3(\mathbb R)\), define
\[
F_u^g(a,b):=g\!\left(\frac{u-a}{1+b}\right),
\qquad |b|\le \frac12.
\]
Assume
\[
M_g:=
\sup_{x\in\mathbb R}
(1+|x|+x^2+|x|^3)\bigl(|g'(x)|+|g''(x)|+|g'''(x)|\bigr)
<\infty.
\]
Then
\begin{align*}
F_u^g(a,b)
&=
g(u)-a g'(u)-u b\, g'(u)
+\frac12 a^2 g''(u) 
+\frac12 b^2\Bigl(2u g'(u)+u^2 g''(u)\Bigr)
+R_u^g(a,b),
\end{align*}
where
\[
\sup_{u,a\in\mathbb R}|R_u^g(a,b)|\cdot \mathbf{1}_{\la |b| \leq 1/2  \ra}
\le
C_g \lb |a|^3+ |ab|+ |b|^3 \rb
\]
for some constant \(C_g\) depending only on \(M_g\).
\end{lemma}

Applying Lemma \ref{Taylor} to $g \equiv \Phi$ and noting that $\mb{E} A=\mb{E} A_0=0$, we get
\begin{align*}
& \mb{E} \left[  \Phi \lb \frac{u-A}{1+B} \rb  - \Phi \lb \frac{u-A_0}{1+B_0} \rb  \right] \\
  =&  \mb{E} \left[  \Phi \lb \frac{u-A}{1+B} \rb \cdot \mathbf{1}_{\la |B| \leq 1/2 \ra} - \Phi \lb \frac{u-A_0}{1+B_0} \rb \cdot \mathbf{1}_{\la |B_0| \leq 1/2 \ra} \right] + \underbrace{O \lb \mb{P} \lb |B|>\frac{1}{2} \rb+  \mb{P} \lb |B_0|> \frac{1}{2} \rb  \rb }_{O \lb \mb{E}(B^2+B_0^2) \rb}\\
  =& -u\phi(u) \cdot \mb{E}\left[ B \cdot \mathbf{1}_{\la |B| \leq 1/2 \ra} - B_0 \cdot \mathbf{1}_{\la |B_0| \leq 1/2 \ra} + \frac{1}{2} \lb A^2-A_0^2 \rb \right] + O \lb \mb{E}(B^2+B_0^2) \rb \\
  +& \lb u \phi(u)+ \frac{u^2}{2} \phi'(u)  \rb \mb{E} \left[ B^2 \cdot \mathbf{1}_{\la |B| \leq 1/2 \ra} - B_0^2 \cdot \mathbf{1}_{\la |B_0| \leq 1/2 \ra} \right] \\
  +& O \lb \mb{E} \left[ |A|^3 + |AB| + |B|^3 + |A_0|^3 + |A_0B_0| + |B_0|^3   \right] \rb.
\end{align*}
Simplifying the above expression using $ B \cdot \mathbf{1}_{\la |B| \leq 1/2 \ra} - B_0 \cdot \mathbf{1}_{\la |B_0| \leq 1/2 \ra} = B-B_0 +  O \lb \mb{E}(B^2+B_0^2) \rb$, we get
\begin{align} 
    & \mb{E} \left[  \Phi \lb \frac{u-A}{1+B} \rb  - \Phi \lb \frac{u-A_0}{1+B_0} \rb  \right] \nonumber \\
     = & -u\phi(u) \cdot \mb{E}\left[ B - B_0  + \frac{1}{2} \lb A^2-A_0^2 \rb \right] + O \lb \mb{E}(B^2+B_0^2) \rb \nonumber \\
     + & O \lb \mb{E} \left[ |A|^3 + |AB| + |B|^3 + |A_0|^3 + |A_0B_0| + |B_0|^3   \right] \rb. \label{Phi taylor}
\end{align}
Now, by Lemma \ref{moment asymptotic}, we have 
\begin{align*}
    \mb{E}\left[ B - B_0  + \frac{1}{2} \lb A^2-A_0^2 \rb \right] &= -\frac{1+r_p}{p} + \frac{1}{p} + O \lb \frac{1}{p^2} + \frac{1}{n\sqrt{pn}} \rb \\
    &+\frac{1}{2} \lb \frac{1+2r_p+r_p^2}{p} - \frac{1}{p} \rb + O \lb \frac{1}{p^2}+\frac{1}{n^2} \rb \\
    &= \frac{r_p^2}{2p} + o(1/n) = \frac{\tau}{2n} + o(1/n).
\end{align*}
For the remaining terms in \eqref{Phi taylor}, note that Lemma \ref{moment asymptotic} implies that $\mb{E} A^2 = O(1/p + 1/n)$. 
\begin{align*}
    & \mb{E}(B^2+B_0^2) +  \mb{E} \left[ |A|^3 + |AB| + |B|^3 + |A_0|^3 + |A_0B_0| + |B_0|^3   \right] \\
    \leq & \mb{E}(B^2+B_0^2)  +   \mb{E} \left[ |A|^3 + |B|^3 + |A_0|^3 + |B_0|^3   \right] + O \lb \sqrt{\mb{E}A^2 \cdot \mb{E} B^2} +  \sqrt{\mb{E}A_0^2 \cdot \mb{E} B_0^2}\rb\\
    = & O \lb  \frac{1}{p^2} + \frac{1}{pn} \rb + O \lb \frac{1}{p^{3/2}} + \frac{1}{n^{3/2}} \rb + O \lb \frac{1}{p^{3/2}} + \frac{1}{p\sqrt{n}} + \frac{1}{n\sqrt{p}} \rb.
\end{align*}
Therefore,
\[
\mb{E} \left[  \Phi \lb \frac{u-A}{1+B} \rb  - \Phi \lb \frac{u-A_0}{1+B_0} \rb  \right] = -u\phi(u)\cdot \frac{\tau}{2n} + O \lb \frac{1}{p^{3/2}} + \frac{1}{n^{3/2}} + \frac{1}{pn} +\frac{1}{p\sqrt{n}} + \frac{1}{n\sqrt{p}} \rb.
\]
Similarly, by applying Lemma \ref{Taylor} to $g\equiv \Psi$ and performing the same moment-matching computation, we get
\[
\mb{E} \left[  \frac{1}{p}\Psi \lb \frac{u-A}{1+B} \rb  - \frac{1}{p}\Psi \lb \frac{u-A_0}{1+B_0} \rb  \right] =o(1/n).
\]
This is because the extra factor $1/p$ forces all the terms to be $o(1/n)$. The dominant term in the approximation is $O(1/p^{3/2})$, and this is precisely why we need $p/n^{2/3} \to \infty$. The proof is completed by combining the two expansions for $\Phi$ and $\Psi$. $\hfill$ $\square$


\section{Proof of Proposition \ref{local power FvML}} As in the proof of Theorem \ref{null-dist}, we need to check three conditions:

\begin{itemize}
    \item {\it Convergence in finite-dimensional distributions}. Recall $S_n(t)$ in \eqref{S_nt}. We need to show that for all $t_1<t_2<\dots<t_k$, 
    \begin{align} \label{fvml-fd conv}
        \lb S_n(t_1),\dots,S_n(t_k) \rb \stackrel{d}{\to} \lb B_{\Phi(t_1)} - \frac{\tau^2 \exp \lb -t_1^2/2 \rb}{2\sqrt{\pi}}, \dots, B_{\Phi(t_k)} -  \frac{\tau^2 \exp \lb -t_k^2/2 \rb}{2\sqrt{\pi}} \rb
    \end{align}
    under the FvML distributions with concentration parameter $\kappa_n$.
    \item {\it Tightness}. This condition is equivalent to \eqref{tightness} for all spaces $D[a,b]$ with $a<b$.

    \item {\it Negligibility of the tail}. This condition is \eqref{tail-negligible}.
\end{itemize}
To show \eqref{fvml-fd conv}, recall that by Proposition \ref{clt for Ustat} and the Cr\'amer--Wold device, we have 
\[
        \lb S_n(t_1),\dots,S_n(t_k), R_n \rb \stackrel{d}{\to} \lb \bm{B}_k, Z \rb
\]
under $\mu_0$, where $Z$ is a standard normal, $R_n$ is the Rayleigh test as in \eqref{Rayleigh}, and $\bm{B}_k= \lb B_{\Phi(t_1)},\dots,B_{\Phi(t_k)} \rb$ has a distribution equal to the joint distribution of the discretized Brownian bridge at $\Phi(t_1), \Phi(t_2),\dots,\Phi(t_k)$. Moreover, the correlation between $Z$ and $\bm{B}_k$ can be specified as follows (this is also the covariance limit in the proof of Proposition \ref{Kolmogorov expansion fvml}):
\[
\mb{E} \lb B_{\Phi(t_i)} Z  \rb =  \mb{E} \lb Z^{*}\cdot \mathbf{1}_{\la Z^{*} \leq t_i \ra} \rb =  \frac{-\exp \lb -t_i^2/2 \rb}{\sqrt{2\pi}}
\]
with $Z^* \sim N(0,1)$, for all $1 \leq i \leq k$.

Consequently, we obtain \eqref{fvml-fd conv} from the convergence above by using the LAN expansion \eqref{fvml LAN} (see \cite{Cutting-P-V} for a proof) and Le Cam's third lemma.

We now show \eqref{tightness} and \eqref{tail-negligible}. Since their proofs are similar, we will only show \eqref{tightness}. Assume the contrary. Then there exist $\ve, \ve_1>0$ and a sequence $\la \delta_k, n_k \ra_{k \geq1}$ such that
\begin{align} \label{contra 1}
\liminf_{k \to \infty} \la \mb{P}_{\mu_{n_k}} \lb \sup_{|t-s| \leq \delta_k} \Big| S_{n_k}(t)-S_{n_k}(s)  \Big| > \ve \rb  \ra \geq \ve_1. 
\end{align}
where $\mu_{n_k}$ is the corresponding subsequence of FvML alternatives. Put
\[
\mathcal{A}_{k}= \la  \sup_{|t-s| \leq \delta_k} \Big| S_{n_k}(t)-S_{n_k}(s)  \Big| > \ve \ra.
\]
Recall $L_n$ in \eqref{fvml LLR}. Thanks to Proposition \ref{2nd moment LLR}, we have
\begin{align}
    \mb{P}_{\mu_{n_k}} \lb \mathcal{A}_{k} \rb = \int_{\mathcal{A}_{k}} 1 d\mb{P}_{\mu_{nk}} 
    &= \int_{\mathcal{A}_{k}} L_n d \mb{P}_0 \nonumber \\
    &\leq \sqrt{\mb{P}_0 \lb \mathcal{A}_k \rb} \cdot \sqrt{\mb{E}_{\mb{P}_0} \lb L_n^2 \rb} \to 0 \label{contra 2}
\end{align}
since ${\mb{E}_{\mb{P}_0} \lb L_n^2 \rb} < \infty$ and $\mb{P}_0 \lb \mathcal{A}_{k} \rb \to 0$ under uniformity (which is due to the fact that \eqref{tightness} holds under uniformity). Note that the second equality in the display above follows from the fact that the distribution of $S_n(t)$ is invariant under rotations:
\[
S_n(t) \lb \bm{X}_1,\dots,\bm{X}_n  \rb \stackrel{d}{=} S_n(t) \lb \bm{O}\bm{X}_1,\dots,\bm{O}\bm{X}_n  \rb
\]
for all orthogonal matrices $\bm{O}$. 

Since \eqref{contra 1} contradicts \eqref{contra 2}, \eqref{tightness} must hold. This finishes the proof. $\hfill$ $\square$


\begin{appendix}

 \section{Technical results, discussions and other proofs} \label{aux}

\subsection{Proof of Proposition \ref{clt for Ustat}} \label{proof-clt Ustat}

To prove (\ref{reduced-fd}), we employ a martingale central limit theorem. We will use Corollary 3.1 in \cite{Hall}. Define 
\begin{align}
    Z_n: &= \sum_{1 \leq i<j \leq n} h_n \lb \bm{X}^{\top}_i \bm{X}_j \rb \label{def-Zn} \\
    s_n^2: &= \mbox{Var}(Z_n), \nonumber\\
    Y_{n,i} :&= \sum_{j=1}^{i-1} h_n(\bm{X}^{\top}_i \bm{X}_j), \nonumber \\
    Q_n :&= \sum_{i=2}^{n} \mb{E} \lb Y_{n,i}^2 \Big| \bm{X}_1, \bm{X}_2, \dots, \bm{X}_{i-1} \rb. \label{def-Q}
\end{align}
Thanks to Lemma \ref{degenerate-Ustat} in Section \ref{aux}, the sequence $\la Y_{n,i}; 2 \leq i \leq n \ra$ is a martingale difference sequence with respect to the natural sigma-fields $\mathcal{F}_i = \sigma \lb \bm{X}_1, \bm{X}_2,\dots, \bm{X}_{i} \rb$. This means that $\mb{E} \lb Y_{n,i} | \mathcal{F}_{i-1} \rb = 0$. To obtain (\ref{reduced-fd}) via the martingale CLT from Corollary 3.1 in \cite{Hall}, we need to verify the following two conditions:
\begin{align} \label{lindeberg-cond}
   s_n^{-2} \sum_{i=2}^{n} \mb{E} \Big[ Y_{n,i}^2 \cdot \mathbf{1}_{\la |Y_{n,i}|>\ve s_n \ra} \Big] \to 0,
\end{align}
for every fixed $\ve >0$, and
\begin{align} \label{cond-var}
    s_n^{-2} Q_n \xrightarrow{\mb{P}} 1.
\end{align}
To verify the Lindeberg condition (\ref{lindeberg-cond}), it suffices to show that
\begin{align} \label{Lindeberg-4th moment}
    s_n^{-4} \sum_{i=2}^{n} \mb{E} Y_{n,i}^4 \to 0.
\end{align}
By the pairwise independence property (see Lemma \ref{degenerate-Ustat}), we get $s_n^2 = n^2(1+o(1)) \cdot \mb{E} h_n^2(\bm{X}^{\top}_1 \bm{X}_2 )$ and thus, 
\begin{align} 
    s_n^4 = n^4(1+o(1)) \cdot  \sigma^4, \label{s_n}
\end{align}
where $\sigma^2$ is the limit in (\ref{conv-var}). To bound the fourth-moment terms in (\ref{Lindeberg-4th moment}), we first note that
\begin{align*}
    \mb{E} \Big[  h_n^2(\bm{X}^{\top}_1 \bm{X}_2 ) \cdot  h_n(\bm{X}^{\top}_1 \bm{X}_3 ) \cdot h_n(\bm{X}^{\top}_1 \bm{X}_4 ) \Big] & = \mb{E} \mb{E} \Big[  h_n^2(\bm{X}^{\top}_1 \bm{X}_2 ) \cdot  h_n(\bm{X}^{\top}_1 \bm{X}_3) \cdot h_n(\bm{X}^{\top}_1 \bm{X}_4) 
 \Big| \bm{X}_1 \Big] \\
 &= \mb{E} \Big[ \mb{E} \lb h_n^2(\bm{X}^{\top}_1 \bm{X}_2 ) \Big| \bm{X}_1 \rb \cdot \mb{E} \lb h_n(\bm{X}^{\top}_1 \bm{X}_3 ) \Big| \bm{X}_1 \rb\\
 & \cdot  \mb{E} \lb h_n(\bm{X}^{\top}_1 \bm{X}_4 ) \Big| \bm{X}_1 \rb \Big] \\
 &= 0,
 \end{align*}
 due to the assumption $\mb{E} h_n (\bm{X}_1^{\top} \bm{X}_2)=0 $ and Lemma \ref{degenerate-Ustat}. Similarly,
 \begin{align*}
      \mb{E} \Big[  h_n^3(\bm{X}^{\top}_1 \bm{X}_2) \cdot  h_n(\bm{X}^{\top}_1 \bm{X}_3 ) \Big] & = \mb{E} \mb{E} \Big[  h_n^3(\bm{X}^{\top}_1 \bm{X}_2) \cdot  h_n(\bm{X}^{\top}_1 \bm{X}_3 ) \Big| \bm{X}_1 \Big] \\
 & = \mb{E} \Big[ \mb{E} \lb h_n^3(\bm{X}^{\top}_1 \bm{X}_2 ) \Big| \bm{X}_1 \rb \cdot \mb{E} \lb h_n(\bm{X}^{\top}_1 \bm{X}_2 ) \Big| \bm{X}_1 \rb \Big]\\
 & = 0.
\end{align*}
Consequently, for some universal constant $C$, we get 
\begin{align*}
    \mb{E} Y_{n,i}^4 & \leq  \sum_{j=1}^{i-1} \mb{E} h_n^4 ( \bm{X}^{\top}_i \bm{X}_j ) + C \cdot \sum_{1 \leq r \neq t \leq i-1} \mb{E} \Big[ h_n^2(\bm{X}^{\top}_i \bm{X}_r ) \cdot h_n^2(\bm{X}^{\top}_i \bm{X}_t ) \Big] \\
    &= (i-1) \cdot \mb{E} h_n^4 (\bm{X}^{\top}_1 \bm{X}_2 ) + Ci^2 \cdot \mb{E} \Big[ h_n^2(\bm{X}^{\top}_1 \bm{X}_2) \cdot h_n^2(\bm{X}^{\top}_1 \bm{X}_3) \Big]\\
    &= (i-1) \cdot \mb{E} h_n^4 (\bm{X}^{\top}_1 \bm{X}_2) + Ci^2 \cdot \lb \mb{E} \Big[ h_n^2(\bm{X}^{\top}_1 \bm{X}_2 ) \Big] \rb^2. 
\end{align*}
where the last assertion follows from the second statement of Lemma \ref{degenerate-Ustat}. Summing the above display over all $2 \leq i \leq n$, we arrive at
\begin{align*}
    \sum_{i=2}^{n} \mb{E} Y_{n,i}^4 &\leq C_1 \cdot \Big(  n^2 \cdot \mb{E} h_n^4 (\bm{X}^{\top}_1 \bm{X}_2) + n^3 \cdot \mb{E}^2 \Big[ h_n^2(\bm{X}^{\top}_1 \bm{X}_2) \Big] \Big) \\
    & \leq (C_1+1) \cdot n^3 \cdot \mb{E} h_n^4 (\bm{X}^{\top}_1 \bm{X}_2),
\end{align*}
for some universal constant $C_1$. This in turn yields
$$s_n^{-4} \sum_{i=2}^{n} \mb{E} Y_{n,i}^4 \leq (C_1+1) \cdot \frac{\mb{E} \Big[ h_n^4( \bm{X}^{\top}_1 \bm{X}_2 ) \Big]}{n \cdot \mb{E}^2 \Big[  h_n^2(\bm{X}^{\top}_1 \bm{X}_2 ) \Big] }$$
by \eqref{s_n}. The last term in the above display goes to $0$ by assumption (\ref{lindeberg}) and thus implies (\ref{Lindeberg-4th moment}). This in turn concludes the proof of (\ref{lindeberg-cond}). 

We now prove (\ref{cond-var}). Recall $Q_n$ defined in (\ref{def-Q}), which can be rewritten as
\begin{align*}
    Q_n &= \sum_{i=2}^{n} \Big[  \sum_{j=1}^{i-1} \mb{E} \lb h_n^2(\bm{X}^{\top}_i \bm{X}_j) \Big| \bm{X_j} \rb + \sum_{1 \leq r \neq t \leq i-1 } \mb{E} \lb h_n(\bm{X}^{\top}_i \bm{X}_r) \cdot h_n(\bm{X}^{\top}_i \bm{X}_t) \Big| \bm{X}_{r}, \bm{X}_t \rb \Big] \\
    &= \sum_{i=2}^{n} (i-1) \cdot \mb{E} h_n^2(\bm{X}^{\top}_1 \bm{X}_2) + \sum_{i=2}^{n} \sum_{1 \leq r \neq t \leq i-1} \mb{E} \lb h_n(\bm{X}^{\top}_i \bm{X}_r) \cdot h_n(\bm{X}^{\top}_i \bm{X}_t) \Big| \bm{X}_{r}, \bm{X}_t \rb.
\end{align*}
Let $\bm{X}, \bm{Y}, \bm{Z}$ be i.i.d. realizations of the uniform distribution on $\mb{S}^{p-1}$. Define 
\begin{align} \label{def-H}
    H_n(\bm{X}, \bm{Y}):= \mb{E} \Big[  h_n \lb \bm{X}^T \bm{Z} \rb \cdot h_n \lb  \bm{Y}^T \bm{Z} \rb \Big| \bm{X}, \bm{Y} \Big].
\end{align}
 It is easy to see that for any $1 \leq r \neq t \leq i-1$, we have 
$$H_n(\bm{X}_r, \bm{X}_t) =  \mb{E} \lb h_n(\bm{X}^{\top}_i \bm{X}_r) \cdot h_n(\bm{X}^{\top}_i \bm{X}_t) \Big| \bm{X}_{r}, \bm{X}_t \rb. $$
Thus, we have 
\begin{align*}
    Q_n &= \sum_{i=2}^{n} (i-1) \cdot \mb{E} h_n^2(\bm{X}^{\top}_1 \bm{X}_2) + \sum_{i=2}^{n} \sum_{1 \leq r \neq t \leq i-1} H_n(\bm{X}_r, \bm{X}_t)\\
    &= s_n^2 \cdot (1+o(1))  +  Q_n^{*},
\end{align*}
where 
\begin{align*}
    Q_n^{*} &:= \sum_{i=2}^{n} \sum_{1 \leq r \neq t \leq i-1} H_n(\bm{X}_r, \bm{X}_t).
\end{align*}
To prove (\ref{cond-var}), we only need to show that $\mbox{Var} \lb Q_n^{*}/s_n^2 \rb \to 0$. Let $A = \la (r,t): 1 \leq r \neq t \leq n-1 \ra$. Note that for any two pairs $(r,t) \in A$ and $(r',t') \in A$, we have 
$$\mb{E} \Big[  H_n(\bm{X}_r, \bm{X}_t) \cdot H_n(\bm{X}_{r'}, \bm{X}_{t'}) \Big] =0,$$
unless $\la r,t \ra = \la r', t' \ra$. This is due to the assumption  $\mb{E} h_n (\bm{X}_1^{\top} \bm{X}_2)=0 $ and Lemma \ref{degenerate-Ustat}. Consequently,
\begin{align*}
    \mbox{Var} \lb Q_n^{*} \rb &= \mbox{Var} \Big( \sum_{(r,t) \in A} \Big[ n - \max \la r,t \ra -1 \Big] \cdot H_n(\bm{X}_r, \bm{X}_t) \Big) \\
    & = \sum_{(r,t) \in A} \Big[ n - \max \la r,t \ra -1 \Big]^2 \cdot \mb{E}  H_n^2(\bm{X}_r, \bm{X}_t) \\
    & \leq O(n^4) \cdot \mb{E}    H_n^2(\bm{X}_1, \bm{X}_2),
\end{align*}
where we use the fact that $|A| \leq n^2$ in the last bound. By \eqref{s_n}, $s_n^4= \sigma^4 n^4(1+o(1))$, so it suffices to verify that $\mb{E} H_n^2(\bm{X}_1, \bm{X}_2) \to 0$, which is the content of Lemma \ref{Diaconis-Freeman} in Section \ref{aux}. This concludes the proof of (\ref{cond-var}), which together with (\ref{lindeberg-cond}) implies (\ref{reduced-fd}). $\hfill$ $\square$


\subsection{Proof of Lemma \ref{group action}} \label{appendix group action}

Suppose $\mathbb{Q}(A)=0$ for some measurable subset $A \subset \mathcal{X}$. 
Let $X \sim \mathbb{P}$ and $Y \sim \mathbb{Q}$. 
By the disintegration theorem (see, for example, Appendix~F of \cite{pollard2002user}), 
\[
\mathbb{Q}(A)
=
\int_{\mathcal{Y}}
  \mathbb{Q}\bigl( Y \in A \cap T^{-1}(t) \mid T(Y)=t \bigr)
  \,\mathbb{Q}_T(dt),
\]
and similarly
\[
\mathbb{P}(A)
=
\int_{\mathcal{Y}}
  \mathbb{P}\bigl( X \in A \cap T^{-1}(t) \mid T(X)=t \bigr)
  \,\mathbb{P}_T(dt).
\]
For each $t$, define the $G$-invariant conditional laws on $T^{-1}(t)$ (such sets are called fibers) by
\[
\mathbb{K}_P^t(B)
  := \int_G \mathbb{P}(gX \in B \mid T(X)=t)\, \Pi(dg),
\qquad
\mathbb{K}_Q^t(B)
  := \int_G \mathbb{Q}(gY \in B \mid T(Y)=t)\, \Pi(dg),
\]
for measurable $B \subset T^{-1}(t)$, where $\Pi$ is the normalized Haar probability measure on $G$. 

Roughly speaking, we integrate the disintegration kernels over the group $G$ so that the resulting kernels are $G$-invariant, up to some null sets in $\mc{Y}$. Because $\mathbb{P}$ and $\mathbb{Q}$ are $G$-invariant, we still have  
\[
\mathbb{Q}(A)
  = \int_{\mathcal{Y}} \mathbb{K}_Q^t(A \cap T^{-1}(t)) \,\mathbb{Q}_T(dt),
\qquad
\mathbb{P}(A)
  = \int_{\mathcal{Y}} \mathbb{K}_P^t(A \cap T^{-1}(t)) \,\mathbb{P}_T(dt).
\]

Each fiber $T^{-1}(t)$ is a single $G$-orbit, in the sense that it is generated by $\la g x_t, g \in G \ra$ for some $x_t \in \mc{X}$, because $T$ is a maximal invariant. 
Thus, $G$ acts transitively on every such fiber. 
By Theorem~4.5 of \cite{eaton1989group}, a transitive compact group action admits a unique $G$-invariant probability measure on each orbit. 
Hence,
\begin{align} \label{kernel equal}
\mathbb{K}_P^t \equiv \mathbb{K}_Q^t
\quad\text{for $\mathbb{Q}_T$-almost all $t$}.
\end{align}

Since $\mathbb{Q}(A)=0$, we have 
\[
\mathbb{K}_Q^t(A \cap T^{-1}(t)) = 0
\quad\text{for $\mathbb{Q}_T$-almost surely $t$},
\]
and therefore for $\mathbb{P}_T$-almost every $t$ as well, because $\mathbb{P}_T \ll \mathbb{Q}_T$. 
Using $\mathbb{K}_P^t = \mathbb{K}_Q^t$ on this set of $t$,
\[
\mathbb{P}(A)
= \int_{\mathcal{Y}} \mathbb{K}_P^t(A \cap T^{-1}(t)) \,\mathbb{P}_T(dt)
= 0.
\]
Thus $\mathbb{P} \ll \mathbb{Q}$.

Now, for the statement regarding the likelihood ratio, let 
\[
L(t) := \frac{d\mathbb{P}_T}{d\mathbb{Q}_T}(t).
\]
Then, by \eqref{kernel equal},
\[
\mathbb{P}(A)
= \int_{\mathcal{Y}}
    \mathbb{K}_P^t(A \cap T^{-1}(t)) 
    \,\mathbb{P}_T(dt)
= \int_{\mathcal{Y}}
    \mathbb{K}_Q^t(A \cap T^{-1}(t)) L(t)
    \,\mathbb{Q}_T(dt).
\]
On the other hand,
\[
\int_{\mathcal{X}} \mathbf{1}_A(x)\, L(T(x)) \,\mathbb{Q}(dx)
=
\int_{\mathcal{Y}}
    \mathbb{K}_Q^t(A \cap T^{-1}(t)) L(t)
    \,\mathbb{Q}_T(dt).
\]
Hence,
\[
\mathbb{P}(A)
= \int_A L(T(x)) \,\mathbb{Q}(dx),
\]
which implies
\[
\frac{d\mathbb{P}}{d\mathbb{Q}}(x)
= L(T(x))
= \frac{d\mathbb{P}_T}{d\mathbb{Q}_T}(T(x)),
\quad \mathbb{Q}\text{-a.s.}
\]
This completes the proof. 
\hfill $\square$


\subsection{Proof of Lemma \ref{moment asymptotic}} \label{moment compute}

Let us start with some basic moment computations for $T$. It is easy to see that $T$ has density proportional to
\[
e^{\kappa t^2} (1-t^2)^{(p-3)/2}\, \mathbf{1}_{[-1,1]}(t).
\]
Recall that $\delta_p:= p/2-\kappa$ and $r_p:= \kappa/\delta_p$. Note that \eqref{watson scaling limit} is equivalent to 
\[
\frac{n}{p} r_p^2 \to \tau \in (0,\infty). 
\]
We first claim that
\begin{align} \label{T moment}
    \mb{E} T^2 = \frac{1+r_p}{p} + O \lb \frac{(1+r_p)^3}{p^2} \rb; \qquad \mb{E}T^4 =   O \lb \frac{(1+r_p)^2}{p^2} \rb.
\end{align}
Let us first relate the second moment of $T$ to its fourth moment. A direct computation gives
\[
\frac{d}{dt} \left[ t\, e^{\kappa t^2} \, (1-t^2)^{(p-1)/2}  \right] = e^{\kappa t^2} \lb 1 - t^2 \rb^{(p-3)/2} \left[ 1 + (2\kappa-p)t^2 -2\kappa t^4 \right].
\]
Integrating the above from $-1$ to $1$ and noting that the left-hand side vanishes, we obtain
\[
1 - 2\delta_p \mb{E} T^2 -2\kappa \mb{E} T^4 = 0.
\]
Equivalently,
\[
\mb{E}T^2 = \frac{1}{2\delta_p} - r_p \mb{E} T^4 = \frac{1+r_p}{p} -  r_p \mb{E} T^4.
\]
Therefore, to prove \eqref{T moment}, we only need to prove the second statement regarding the fourth moment. To do this, write 
\begin{align*}
    \mb{E} T^4 = \frac{\int_{-1}^1 t^4\,e^{\kappa t^2} (1-t^2)^{(p-3)/2} dt}{\int_{-1}^1 e^{\kappa t^2} (1-t^2)^{(p-3)/2}dt} = \frac{N_{4,p}}{D_{4,p}}.
\end{align*}
Let us bound the numerator above. By using the bound $\log(1-x)\leq -x$ for all $x \in [0,1)$, we obtain
\begin{align*}
    \kappa t^2 + \frac{p-3}{2} \log(1-t^2) \leq \lb \kappa - \frac{p-3}{2} \rb t^2 = \lb - \delta_p + \frac{3}{2} \rb t^2. 
\end{align*}
Therefore, for all large $p$, 
\[
N_{4,p} \leq e^{-3/2} \int_{\mb{R}} t^4 e^{-\delta_p t^2}dt = O \lb \delta_p^{-5/2} \rb.
\]
Next, we bound the denominator $D_{4,p}$ from below. Observe that
\begin{align} \label{delta_p}
\delta_p = \frac{p}{2(1+r_p)}
\end{align}
and $r_p = O\lb \sqrt{p/n} \rb$. Therefore, $\delta_p/\sqrt{p} \to \infty$ as $p \to \infty$.

Now, for $x \in [0,1/100]$, we have the bound $\log(1-x) \geq -x- x^2$. Therefore, for sufficiently large $p$, the following bound is valid on the set $\la |t| \leq \delta_p^{-1/2} \ra$:
\begin{align*}
    \kappa t^2 + \frac{p-3}{2} \log \lb 1 -t^2 \rb \geq \kappa t^2 - \frac{p-3}{2} t^2 - \frac{p-3}{2}t^4 &= -\lb \delta_p - \frac{3}{2} \rb t^2 - \frac{p-3}{2} t^4 \\
    &\geq -1 -\frac{2p}{\delta_p^2} \geq -2. 
\end{align*}
Thus,
\[
D_{4,p} \geq \int_{-\delta_p^{-1/2}}^{\delta_p^{-1/2}} e^{-2} dt = \Omega \lb \delta_p^{-1/2} \rb.
\]
Consequently, we get \eqref{T moment} from the bounds on $N_{4,p}$ and $D_{4,p}$, together with \eqref{delta_p}.

We are now ready to estimate the moments of $A,A_0,B,B_0$. Let us start with $A$ and $A_0$. Since both have symmetric densities, their odd moments vanish. A direct computation gives
\begin{align*}
    \mb{E} A^2 = p (\mb{E} T^2)^2 &= p \lb  \frac{1+r_p}{p} + O \lb \frac{(1+r_p)^3}{p^2} \rb \rb^2 \\
   & = \frac{(1+r_p)^2}{p} + O \lb  \frac{(1+r_p)^4}{p^2} +  \frac{(1+r_p)^6}{p^3} \rb \\
   &=  \frac{(1+r_p)^2}{p} + O \lb \frac{1}{p^2} + \frac{r_p^4}{p^2} + \frac{1}{p^3} + \frac{r_p^6}{p^3} \rb.
\end{align*}
Since $r_p = O(\sqrt{p/n})$, we conclude that
\[
    \mb{E} A^2 = \frac{(1+r_p)^2}{p} + O \lb \frac{1}{p^2} + \frac{1}{n^2} \rb.
\]
Let us estimate $\mb{E} |A|^3$. Write
\begin{align*}
    \mb{E} |A|^3 = p^{3/2} \lb  \mb{E}|T|^3 \rb^2  &\leq p^{3/2} \cdot \mb{E} T^2 \cdot \mb{E} T^4 \\
    &= p^{3/2} \cdot \left[ \frac{1+r_p}{p} + O \lb \frac{(1+r_p)^3}{p^2} \rb \right] \cdot    O \lb \frac{(1+r_p)^2}{p^2} \rb \\
    &= O \lb \frac{(1+r_p)^3}{p^{3/2}} \rb + O \lb \frac{(1+r_p)^5}{p^{5/2}} \rb \\
    &= O \lb \frac{1}{p^{3/2}} + \frac{1}{n^{3/2}} + \frac{1}{p^{5/2}} + \frac{1}{n^{5/2}} \rb = O \lb \frac{1}{p^{3/2}} + \frac{1}{n^{3/2}} \rb.
\end{align*}
For $A_0$, it is easy to check that $\mb{E} A_0^2=1/p$ and $\mb{E} |A_0|^3= O(p^{-3/2})$. 

Let us now prove the moment estimates for $B$ and $B_0$. It is obvious that $B$ and $B_0$ take values in $[-1,0]$. Put $S=T^2+ T'^2 - T^2T'^2$. By using the expansion $\sqrt{1-s}=1-s/2+r(s)$ with $|r(s)|\leq s^2/2$ for all $s \in [0,1]$, we obtain 
\begin{align*}
    \mb{E} B = -\frac{\mb{E} S}{2} + \mb{E}[r(S)]  &= - \mb{E}[T^2] + O \lb \mb{E} [S^2] \rb \\
    &=  - \mb{E}[T^2] +\frac{1}{2} \underbrace{\lb \mb{E}[T^2] \rb^2}_{=O \lb \mb{E}[T^4] \rb}  + O \lb \mb{E} [S^2] \rb  \\
    &=  - \mb{E}[T^2] + O \lb \mb{E}[T^4] + \mb{E} \left[ T'^4 \right] + \mb{E} \left[ T^4T'^4  \right] \rb \\
    &=   - \mb{E}[T^2] + O \lb \mb{E}[T^4] \rb \\
    &= -\frac{1+r_p}{p} + O \lb \frac{(1+r_p)^3}{p^2} \rb +  O \lb \frac{(1+r_p)^2}{p^2} \rb \\
    &=  -\frac{1+r_p}{p} + O \lb \frac{(1+r_p)^3}{p^2} \rb =   -\frac{1+r_p}{p} + O \lb \frac{1}{p^2} + \frac{1}{n\sqrt{pn}} \rb.
\end{align*}
The estimate for $\mb{E}B_0$ can be proved similarly (and is easier, since it is under the uniform distribution). 

We now estimate $\mb{E} B^2$. By using the inequality $1-\sqrt{(1-x)(1-y)} \leq x+y$ for all $0 \leq x,y \leq 1$, we have 
\[
\mb{E} B^2 \leq \mb{E} \lb T^2 + T'^2 \rb^2 \leq 2 \mb{E}T^4 = O \lb \frac{(1+r_p)^2}{p^2} \rb =  O \lb \frac{1}{p^2} + \frac{1}{pn} \rb
\]
due to \eqref{T moment}. The estimate for $\mb{E} B_0^2$ can be proved similarly and is in fact easier. This finishes the proof. $\hfill$ $\square$


\subsection{Proof of Lemma \ref{Taylor}}
Fix \(u\in\mathbb R\), and define
\[
h(t):=F_u^g(ta,tb)=g\!\left(\frac{u-ta}{1+tb}\right),
\qquad t\in[0,1].
\]
Write
\[
z_t:=\frac{u-ta}{1+tb}.
\]
Then \(h(t)=g(z_t)\), and
\[
z_t'=-\frac{a+bu}{(1+tb)^2},
\qquad
z_t''=\frac{2b(a+bu)}{(1+tb)^3},
\qquad
z_t'''=-\frac{6b^2(a+bu)}{(1+tb)^4}.
\]
Since
\[
u-ta=(1+tb)z_t,
\]
we have
\[
a+bu=a+b(u-ta)+tab=a+b(1+tb)z_t+tab,
\]
hence, on \(|b|\le 1/2\),
\[
|a+bu|\le \frac{3}{2}\bigl(|a|+|b|(1+|z_t|)\bigr).
\]

\noindent \underline{\textit{Step 1: Bounding the third derivative of $h$.}}
We will show that
\begin{align} \label{3rd derivative}
\max_{t\in [0,1]}|h'''(t)|\le C_g (|a|+|b|)^3
\end{align}
for some $C_g>0$ depending only on $M_g$.

A straightforward computation yields
\begin{align*}
h'''(t) &=   g'''(z_t) \lb z'_t \rb^3 + 3g''(z_t)z'_tz''_t + g'(z_t)z'''_t \\
&= -\frac{1}{(1+tb)^3}
\Bigl[
g'''(z_t)(a+bz_t)^3
+
6b\,g''(z_t)(a+bz_t)^2
+
6b^2\,g'(z_t)(a+bz_t)
\Bigr].
\end{align*}
For the first term,
\begin{align*}
|g'''(z_t)(a+bz_t)^3|
&\le \frac{M_g}{(1+|z_t|)^3}\bigl(|a|+|b|\,|z_t|\bigr)^3 \\
&\le  M_g (|a|+|b|)^3,
\end{align*}
because \(|z_t|/(1+|z_t|)\le 1\).

For the second term, by the same reasoning,
\begin{align*}
|b\,g''(z_t)(a+bz_t)^2|
&\le |b|\,\frac{M_g}{(1+|z_t|)^2}\bigl(|a|+|b|\,|z_t|\bigr)^2 \\
&\le  M_g\,|b|(|a|+|b|)^2 \\
&\le  M_g (|a|+|b|)^3.
\end{align*}

Similarly, for the third term,
\begin{align*}
|b^2 g'(z_t)(a+bz_t)|
&\le |b|^2\,\frac{M_g}{1+|z_t|}\bigl(|a|+|b|\,|z_t|\bigr) \\
&\le  M_g\,|b|^2(|a|+|b|) \\
&\le  M_g (|a|+|b|)^3.
\end{align*}
Putting everything together yields \eqref{3rd derivative}.



\noindent \underline{\textit{Step 2: Expanding $h$.}} Taylor expansion at \(t=0\) gives
\[
h(1)=h(0)+h'(0)+\frac12 h''(0)+\frac16 h'''(\xi)
\]
for some \(\xi\in(0,1)\). Therefore, it remains to compute \(h(0)\), \(h'(0)\), and \(h''(0)\).

Since \(z_0=u\), \(z_0'=-(a+bu)\), and \(z_0''=2b(a+bu)\), we obtain
\[
h(0)=g(u),
\]
\[
h'(0)=g'(u)z_0'=-a g'(u)-u b\,g'(u),
\]
and
\[
h''(0)=g''(u)(z_0')^2+g'(u)z_0''.
\]
Also,
\[
(z_0')^2=(a+ub)^2=a^2+2aub+u^2b^2,
\]
and
\[
z_0''=2ab+2u b^2.
\]
Therefore,
\begin{align*}
h''(0)
&=
a^2 g''(u)
+ \underbrace{2ab\bigl(g'(u)+u g''(u)\bigr)}_{=O(|ab|)}
+b^2\bigl(2u g'(u)+u^2 g''(u)\bigr).
\end{align*}
Substituting into the Taylor expansion of $h$ yields the desired expansion, with
\[
R_u^g(a,b)= ab(g'(u)+ug''(u))+\frac16 h'''(\xi).
\]
The proof is completed.
$\hfill$ $\square$


\subsection{Likelihood ratio analysis} \label{appendix LLR}
This section contains the analysis of the second moment of the likelihood ratio used in Sections \ref{sec fvml} and \ref{sec low-rank}. 

We first analyze the second moment of the likelihood ratio between the FvML distributions with randomized locations and the uniform distributions. We first parametrize the FvML distributions as
\begin{align*}
    d \mb{P}_{\rm Fvml} := C_{p}(\kappa) \cdot \exp \lb \kappa \langle \bm{\mu}, \bm{x} \rangle  \rb d\mb{P}_0.
\end{align*}
with $\kappa \in (0,\infty)$ and $\bm{\mu} \in \mb{S}^{p-1}$. From this, we find that
\begin{align} \label{Cp(kappa)}
C_{p}(\kappa) =  \left[ \mb{E}_{\mb{P}_0}  \exp \lb  \kappa \langle \bm{\mu}, \bm{x} \rangle \rb   \right]^{-1} =  \left[   \frac{1}{\sqrt{\pi}} \frac{\Gamma\left(\frac{p}{2}\right)}{\Gamma\left(\frac{p-1}{2}\right)} \cdot \int_{-1}^{1} e^{\kappa t} \left(1-t^2\right)^{\frac{p-3}{2}} dt \right]^{-1}.
\end{align}
Some basic properties of the above normalizing constant are collected in Lemma \ref{Bessel}. Define the likelihood ratio 
\begin{align} \label{fvml LLR}
   L_n:= \mb{E}_{\bm{\mu} \sim \mb{P}_0} \left[ \frac{d\mb{P}_{\rm Fvml}^{\otimes n}}{d \mb{P}_0^{\otimes n}} \right] = \lb C_{p}(\kappa) \rb^n \cdot  \mb{E}_{\bm{\mu} \sim \mb{P}_0} \left[ \prod_{i=1}^n \exp \Big( \kappa \langle \bm{\mu}, \bm{X}_i \rangle \Big)  \right].
\end{align}
Our first result is an asymptotic formula for the second moment of the likelihood ratio. 
\begin{prop} \label{2nd moment LLR}
    Let $\kappa=\tau p^{3/4}/\sqrt{n}$ for some $\tau>0$. Suppose $\min \la p,n \ra \to \infty$, then
    \[
    \mb{E} \lb L_n^2 \rb = \exp \lb \tau^4/2 + o(1) \rb
    \]
    where $L_n$  is defined as in \eqref{fvml LLR}.
\end{prop}

\noindent \textbf{Proof of Proposition \ref{2nd moment LLR}.} Recall the form of $L_n$ in \eqref{fvml LLR}. To compute the second moment, take two independent copies $\bm{\mu}, \bm{\mu}_1$ of the uniform distribution and write
\begin{align*}
    \mb{E} L_n^2 &=  \lb C_{p}(\kappa) \rb^{2n} \cdot \mb{E}_{\bm{X}}\left[ \mb{E}_{(\bm{\mu},\bm{\mu_1})}  \left[ \prod_{i=1}^n \exp \Big( \kappa \langle \bm{\mu}+\bm{\mu}_1, \bm{X}_i \rangle \Big) \right]  \right] \\
    &=    C_{p}(\kappa)^{2n} \cdot  \mb{E}_{(\bm{\mu},\bm{\mu_1})} \left[ \bm{E}_{\bm{X}} \left[ \prod_{i=1}^n \exp \lb \kappa \| \bm{\mu}+\bm{\mu}_1 \| \cdot \langle \frac{\bm{\mu}+\bm{\mu}_1}{\| \bm{\mu}+\bm{\mu}_1 \|}, \bm{X}_i \rangle \rb   \right]  \right] \\
    &= \mb{E} \left[  \frac{    C_{p}(\kappa)^{2n} }{C_{p}(\kappa\| \bm{\mu}+\bm{\mu}_1 \|)^{n} } \right].
\end{align*}
Note that $\| \bm{\mu}+\bm{\mu}_1 \| \stackrel{d}{=} \sqrt{2(1+U)}$, where $U$ has the law as in \eqref{density-inner}. Thus, we have 
\[
    \mb{E} L_n^2 =  \mb{E} \left[  \frac{    C_{p}(\kappa)^{2n} }{C_{p}(\kappa \sqrt{2(1+U)})^{n} } \right] = \mb{E} \left[ \exp \lb n \lb  2 \log C_p(\kappa) - \log C_p \lb \kappa \sqrt{2(1+U)}  \rb  \rb \rb   \right]
\]
where $U$ has a symmetric Beta-type distribution as in \eqref{density-inner}. 

Put 
\[
L_{n1}:=   2\log C_p(\kappa) - \log C_p \lb \kappa \sqrt{2(1+U)} \rb.
\]
By the third property in Lemma \ref{Bessel}, we have 
\begin{align*}
    \left| L_{n1} + \frac{\kappa^2}{p} - \frac{\kappa^2(1+U)}{p}  \right| &\leq 2 \left| \log C_p(\kappa) +  \frac{\kappa^2}{p} \right| + \left| \log C_p \lb \kappa \sqrt{2(1+U)}  \rb + \frac{\kappa^2(1+U)}{p} \right| \\
    &=O \lb \frac{\kappa^4}{p^3} \rb = O \lb \frac{\tau^4}{n^2} \rb.
\end{align*}
Consequently,
\begin{align*}
        \mb{E} L_n^2 =   \mb{E} \left[ \exp \lb n L_{n1} \rb \right] &= \mb{E} \left[ \exp \lb \frac{\kappa^2n}{p} \cdot U + O \lb n^{-1} \rb \rb \right] \\
        &= \mb{E} \left[ \exp \lb  \tau^2 \cdot \sqrt{p}U +  O \lb n^{-1} \rb \rb  \right].
\end{align*}
The proof is completed by noting that the sequence $\la \sqrt{p}U \ra$ converges to a standard normal distribution and is exponentially tight. $\hfill$ $\square$

The next asymptotic result was used in the proof of Theorem \ref{low-rank information bound}. Recall that the multivariate Gamma function $\Gamma_n(z)$ is defined as
\begin{align} \label{multi Gamma}
    \Gamma_n(z):= \pi^{n(n-1)/4} \prod_{k=1}^n \Gamma \lb z - \frac{k-1}{2} \rb
\end{align}
for all complex numbers $z$ such that $\mbox{Re}(z)> (n-1)/2$.

The multivariate Gamma function reduces to the usual Gamma function when $n=1$. The following lemma is taken from Lemma 5.1 and Proposition 5.1 in \cite{jiang2015likelihood}.

\begin{lemma} \label{gamma asymptotic}
Let $\Gamma(x)$ be the standard Gamma function and $\Gamma_n(x)$ be the multivariate Gamma function as in \eqref{multi Gamma}. We have
\begin{itemize}
    \item Uniformly for all $b \in [-x/2,x/2]$,
    \[
    \log \left[ \frac{\Gamma(x+b)}{\Gamma(x)} \right] = (x+b)\log(x+b) - x\log x - b-\frac{b}{2x} + O \lb \frac{b^2+1}{x^2} \rb
    \]
    as $x \to \infty$.
    \item Uniformly for all $t \in [-p/n,\, p/n]$, 
\[
    \log \left[ \frac{\Gamma_n\!\left( \frac{p}{2} + t \right)}{\Gamma_n\!\left( \frac{p}{2} \right)} \right]
    =  
    \alpha_{n,p}\, t 
    + \beta_{n,p}\, t^2 
    + \gamma_{n,p}(t)
    + o(1)
\]
as $n,p \to \infty$ with $p/n \to \infty$, where
\begin{align*}
    \alpha_{n,p}
    &:= -\left[
        2n 
        + \left(p-n-\frac{1}{2}\right)
          \log\!\left(1 - \frac{n}{p}\right)
    \right], \\
    \beta_{n,p}
    &:= -\left[
        \frac{n}{p}
        + \log\!\left(1 - \frac{n}{p}\right)
    \right], \\
    \gamma_{n,p}(t)
    &:= n\left[
        \left(\frac{p}{2} + t\right)
        \log\!\left(\frac{p}{2} + t\right)
        - \frac{p}{2}\log\!\left(\frac{p}{2}\right)
    \right].
\end{align*}
\end{itemize}
\end{lemma}

\noindent \textbf{Proof of Lemma \ref{gamma asymptotic}.} The first item follows from Lemma 5.1 in \cite{jiang2015likelihood}, and the second item follows from Proposition 5.1 in \cite{jiang2015likelihood}. $\hfill$ $\square$

\begin{prop} \label{Fn asymptotic}
    Recall $F_n$ in \eqref{Fn}. Suppose $p/n \to \infty$, $n(1-k/p) \to 0$, and $n+1 \leq k \leq p$. Then,
    \[
F_n(\Delta) - 2F_n \lb \frac{\Delta}{2} \rb \to 0
\]
where $\Delta:=k-p$.
\end{prop}

\noindent \textbf{Proof of Proposition \ref{Fn asymptotic}.} Write 
\begin{align*}
    F_n(\Delta) - 2F_n \lb \frac{\Delta}{2} \rb &= n \log \left[ \frac{\Gamma \lb \frac{p}{2} \rb}{\Gamma \lb \frac{p}{2} +\Delta \rb}  \right]  +\log \left[  \frac{\Gamma_n \lb \frac{p}{2}+\Delta \rb}{\Gamma_n \lb \frac{p}{2} \rb} \right] \\
    & - 2n \log \left[ \frac{\Gamma \lb \frac{p}{2} \rb}{\Gamma \lb \frac{p+\Delta}{2} \rb}  \right]  - 2\log \left[  \frac{\Gamma_n \lb \frac{p+\Delta}{2} \rb}{\Gamma_n \lb \frac{p}{2} \rb} \right].
\end{align*}
Obviously, $|\Delta|=O(p/n)=o(p)$, so Lemma \ref{gamma asymptotic} applies and yields
\begin{align*}
    & n \log \left[ \frac{\Gamma \lb \frac{p}{2} \rb}{\Gamma \lb \frac{p}{2} +\Delta \rb}  \right] - 2n \log \left[ \frac{\Gamma \lb \frac{p}{2} \rb}{\Gamma \lb \frac{p+\Delta}{2} \rb}  \right] \\
   = & n \left[ - \lb \frac{p}{2} + \Delta \rb \log \lb  \frac{p}{2} + \Delta \rb + \frac{p}{2} \log \lb \frac{p}{2} \rb + \Delta + \frac{\Delta}{p} + O \lb  \frac{\Delta^2+1}{p^2}  \rb      \right] \\
   - &2n \left[ - \lb \frac{p+ \Delta}{2}  \rb \log \lb  \frac{p+ \Delta}{2} \rb + \frac{p}{2} \log \lb \frac{p}{2} \rb + \frac{\Delta}{2} + \frac{\Delta}{2p} + O \lb  \frac{\Delta^2+1}{p^2}  \rb \right] \\
   = & -n \left[   \lb \frac{p}{2} + \Delta \rb \log \lb  \frac{p}{2} + \Delta \rb - 2 \lb \frac{p+ \Delta}{2}  \rb \log \lb  \frac{p+ \Delta}{2} \rb \right] +  O \lb  \frac{n(\Delta^2+1)}{p^2}  \rb \\
   - & \frac{np}{2} \log \lb \frac{p}{2} \rb \\
   = & -n \left[   \lb \frac{p}{2} + \Delta \rb \log \lb  \frac{p}{2} + \Delta \rb - 2 \lb \frac{p+ \Delta}{2}  \rb \log \lb  \frac{p+ \Delta}{2} \rb \right] -\frac{np}{2} \log \lb \frac{p}{2} \rb + o(1),
\end{align*}
where the last line follows from the fact that $n \Delta^2/p = O(1/n)=o(1)$.

Similarly, with $\alpha_{n,p}, \beta_{n,p}, \gamma_{n,p}(t)$ as in Lemma \ref{gamma asymptotic}, we have 
\begin{align*}
  &  \log \left[  \frac{\Gamma_n \lb \frac{p}{2}+\Delta \rb}{\Gamma_n \lb \frac{p}{2} \rb} \right] - 2\log \left[  \frac{\Gamma_n \lb \frac{p+\Delta}{2} \rb}{\Gamma_n \lb \frac{p}{2} \rb} \right] \\
  = & \alpha_{n,p}\, \Delta + \beta_{n,p}\, \Delta^2 + \gamma_{n,p}(\Delta)  -  \lb \alpha_{n,p}\, \Delta + \beta_{n,p}\, \frac{\Delta^2}{2} + 2\gamma_{n,p}(\Delta/2) \rb + o(1) \\
  = &  \beta_{n,p}\, \frac{\Delta^2}{2} +  \gamma_{n,p}(\Delta) -  2\gamma_{n,p}(\Delta/2) + o(1) \\
  = & \beta_{n,p}\, \frac{\Delta^2}{2} +  n\left[
        \left(\frac{p}{2} + \Delta\right)
        \log\!\left(\frac{p}{2} + \Delta\right)
        - \frac{p}{2}\log\!\left(\frac{p}{2}\right)
    \right] \\
    - &2n \left[
        \left(\frac{p+\Delta}{2} \right)
        \log\!\left(\frac{p+\Delta}{2} \right)
        - \frac{p}{2}\log\!\left(\frac{p}{2}\right)
    \right] + o(1) \\
    = & \beta_{n,p}\, \frac{\Delta^2}{2} + n \left[ \left(\frac{p}{2} + \Delta\right)
        \log\!\left(\frac{p}{2} + \Delta\right) - 2  \left(\frac{p+\Delta}{2} \right)
        \log\!\left(\frac{p+\Delta}{2} \right) \right] + \frac{np}{2} \log \lb \frac{p}{2} \rb.
\end{align*}
Thus,
\begin{align*}
 F_n(\Delta) - 2F_n \lb \frac{\Delta}{2} \rb =  \beta_{n,p}\, \frac{\Delta^2}{2} + o(1) &= -\left[ 
        \frac{n}{p}
        + \log\!\left(1 - \frac{n}{p}\right)
     \right]  \frac{\Delta^2}{2} + o(1) \\
     &=  \frac{n^2 \Delta^2}{2p^2}(1+o(1))+o(1)
\end{align*}
which tends to $0$ since $n\Delta/p=-n(1-k/p) \to 0$. The proof is completed. $\hfill$ $\square$

\subsection{Kolmogorov distance asymptotic results} \label{sec Kol asymptotic}

Let us start with a simple observation. 

\begin{lemma} \label{contiguity mean}
    Suppose $\la \mb{P}_n \ra$ and $\la \mb{Q}_n \ra$ are two sequences of probability measures such that the likelihood ratio $L_n:=d \mb{Q}_n/d\mb{P}_n$ exists. Let $\la X_n; n \geq 1 \ra$ be a sequence of random variables. If
    \begin{align*}
    \sup_{n \geq 1} \la \mb{E}_{\mb{P}_n} \lb L_n^2 \rb + \mb{E}_{\mb{P}_n} \lb X_n^4 \rb \ra &< \infty,
    \end{align*}
   then 
    \[
\sup_{n\geq 1} \la  \mb{E}_{\mb{Q}_n} \lb  X_n^2 \rb \ra < \infty \ \text{and} \ \left| \mb{E}_{\mb{P}_n} \lb  X_n \rb - \mb{E}_{\mb{Q}_n} \lb  X_n \rb \right| \leq \sqrt{\mb{E}_{\mb{P}_n} \lb  X_n^2 \rb \cdot \mbox{Var}_{\mb{P}_n} \lb L_n \rb}.
    \]
\end{lemma}

\noindent \textbf{Proof of Lemma \ref{contiguity mean}.} Observe that
\begin{align*}
      \mb{E}_{\mb{Q}_n} X_n^2 &=   \mb{E}_{\mb{P}_n} \lb X_n^2 L_n \rb \leq \sqrt{   \mb{E}_{\mb{P}_n} \lb X_n^4 \rb \cdot    \mb{E}_{\mb{P}_n} \lb L_n^2 \rb} \leq \frac{\mb{E}_{\mb{P}_n} \lb X_n^4 \rb +   \mb{E}_{\mb{P}_n} \lb L_n^2 \rb}{2}.
\end{align*}

The second inequality can be proved similarly. The proof is completed. $\hfill$ $\square$



\noindent \textbf{Proof of Proposition \ref{Kolmogorov expansion fvml}.}  The proof is based on Lemma \ref{contiguity mean} and Proposition \ref{2nd moment LLR} and Le Cam's third lemma. The interesting feature of this approach is that no growth condition on $n$ and $p$ is assumed. We do not know whether direct analyses based on the Edgeworth expansion can yield the same result. The reason that Proposition \ref{Kolmogorov expansion fvml} holds without any growth condition on $n$ and $p$ is due to some special properties of the Bessel functions of the first kind, which were exploited in \cite{Cutting-P-V} and Proposition \ref{2nd moment LLR}.

Fix $u \in \mb{R}$ and consider the sequence of random variables
\[
A_n(u)=A_n:= \sqrt{\frac{2}{n(n-1)}} \sum_{1\leq i<j \leq n} \left[ \mathbf{1}_{\la \sqrt{p}\bm{X}_i^\top \bm{X}_j \leq u \ra} - \mb{P}_{\mu_0} \lb \sqrt{p}\bm{X}_1^\top \bm{X}_2 \leq u \rb \right].
\]
Recall $L_n$ in \eqref{fvml LLR}. It was shown in \cite{Cutting-P-V} that, under uniformity, we have the LAN expansion
\begin{align} \label{fvml LAN}
\log \lb L_n \rb = \frac{\tau^2}{\sqrt{2}}R_n - \frac{\tau^4}{4} + o_{\mb{P}}(1)
\end{align}
where $R_n$ is the Rayleigh test in \eqref{Rayleigh}. 

Recall that $\Phi$ is the CDF of a standard normal. Using Proposition \ref{clt for Ustat} and the Cr\'ammer--Wold device, we have
\[
\lb A_n,  R_n  \rb  \stackrel{d}{\to} N \lb \lb 0, 0  \rb^\top, \lb \begin{matrix}
    \Phi(u) \lb 1 - \Phi(u) \rb &  g(u) \\
    g(u) & 1
\end{matrix} \rb \rb
\]
under uniformity, where
\begin{align*}
\quad g(u):&= \lim_{n \to \infty} \mb{E}_{\mu_0} \left[ \sqrt{p}\bm{X}_1^\top \bm{X}_2 \cdot \mathbf{1}_{\la  \sqrt{p}\bm{X}_1^\top \bm{X}_2 \leq u\ra} \right] = \mb{E}_{\mu_0} \lb Z\cdot \mathbf{1}_{\la Z \leq u \ra} \rb =  \frac{-\exp \lb -u^2/2 \rb}{\sqrt{2\pi}}
\end{align*}
with $Z$ denoting a standard normal random variable in the expression above. The function $g$ is simply the negative of the standard Gaussian density.

The convergence in expectation in the display above follows from the fact that the normalized inner products are asymptotically normal and have uniformly bounded fourth moments. Thus,
\[
\lb A_n,  \log(L_n)  \rb  \stackrel{d}{\to} N \lb \lb 0, \frac{-\tau^4}{4}  \rb^\top, \lb \begin{matrix}
    \Phi(u) \lb 1 - \Phi(u) \rb &  \frac{\tau^2}{\sqrt{2}}g(u) \\
    \frac{\tau^2}{\sqrt{2}}g(u) & \frac{\tau^4}{4}
\end{matrix} \rb \rb
\]
under uniformity. By Le Cam's third lemma and \eqref{fvml LAN},
we obtain that under $\mu_n$,
\[
A_n \stackrel{d}{\to} N \lb \frac{\tau^2g(u)}{\sqrt{2}},     \Phi(u) \lb 1 - \Phi(u) \rb  \rb.
\]
Next, we claim that
\begin{align} \label{4th moment bounded}
    \sup_{n \geq 1} \mb{E} A_n^4 < \infty.
\end{align}
In fact, one can prove \eqref{4th moment bounded} quickly by using the following combinatorial argument. Put $\xi_{ij}=  \mathbf{1}_{\la \sqrt{p}\bm{X}_i^\top \bm{X}_j \leq u \ra} - \mb{P}_{\mu_0} \lb \sqrt{p}\bm{X}_1^\top \bm{X}_2 \leq u \rb$. Write
\[
\mb{E} A_n^4 = \frac{4}{n^2(n-1)^2} \sum_{e_1, e_2, e_3, e_4} \mb{E} \left[ \prod_{r=1}^4  \xi_{e_r} \right]
\]
where $e_r=(i_r,j_r)$ for some $1 \leq i_r<j_r \leq n$.

Associate $\la  e_1, e_2, e_3, e_4 \ra$ with a multigraph on the involved indices. If some vertex appears only once among the eight endpoints, then conditioning on all other $\bm{X}$'s and using the degeneracy of the kernel (Proposition \ref{degenerate-Ustat}), we get
\[
\mb{E} \left[ \prod_{r=1}^4  \xi_{e_r} \right] =0
\]
for such a configuration  $\la  e_1, e_2, e_3, e_4 \ra$. 

Therefore, the only non-zero configurations are those in which each vertex has degree at least two. This means that there are at most four indices in a non-zero configuration, and the number of such configurations is at most $O(n^4)$. Consequently, we get \eqref{4th moment bounded}.

Combining \eqref{4th moment bounded}, the expansion \eqref{fvml LAN}, Lemma \ref{contiguity mean}, and Proposition \ref{2nd moment LLR}, we obtain the mean convergence
\begin{align*}
\frac{\tau^2g(u)}{\sqrt{2}} &=\lim_{n \to \infty} \mb{E}_{\mu_n} \lb A_n \rb \\
&= \lim_{n \to \infty} \la  \frac{n(1+o(1))}{\sqrt{2}} \cdot \left[  \mb{P}_{\mu_n} \lb \sqrt{p}\bm{X}^\top \bm{Y} \leq u  \rb - \mb{P}_0 \lb \sqrt{p}\bm{X}^\top \bm{Y} \leq u  \rb   \right]  \ra.
\end{align*}
This completes the proof of the first claim. 

We now prove the uniform convergence. Observe that for a fixed $M>0$
\begin{align*}
       &  \sup_{u \in \mb{R}} \left|     \sqrt{\frac{n(n-1)}{2}} \Big[ \mb{P}_{\mu_n} \lb \sqrt{p}\bm{X}^\top \bm{Y} \leq u  \rb - \mb{P}_{\mu_0} \lb \sqrt{p}\bm{X}^\top \bm{Y} \leq u  \rb \Big] + \frac{\tau^2}{\sqrt{2\pi}} \exp \lb -u^2/2 \rb \right| \\
      = & \sup_{u \in \mb{R}} \left|  \mb{E}_{\mu_n} \lb A_n(u) \rb +  \frac{\tau^2}{\sqrt{2\pi}} \exp \lb -u^2/2 \rb  \right| \\
      \leq & \sup_{u \in [-M,M] } \left|  \mb{E}_{\mu_n} \lb A_n(u) \rb +  \frac{\tau^2}{\sqrt{2\pi}} \exp \lb -u^2/2 \rb \right| 
      + \sup_{ |u|>M } \left|  \mb{E}_{\mu_n} \lb A_n(u) \rb +  \frac{\tau^2}{\sqrt{2\pi}} \exp \lb -u^2/2 \rb \right|.
\end{align*}
It is a standard fact that if a sequence of equicontinuous functions converges pointwise on a compact set, then the convergence is uniform (this is a corollary of the Arzel\'a--Ascoli theorem; see Theorem 4.43 in \cite{folland1999real}). To use this fact, let us check that the functions
\begin{align} \label{equicontinuity}
u \mapsto   \mb{E}_{\mu_n} \lb A_n(u) \rb +   \frac{\tau^2}{\sqrt{2\pi}} \exp \lb -u^2/2 \rb
\end{align}
are equicontinuous. 

 The second term is obviously smooth and has bounded derivatives, so we only have to treat the first term. By the second estimate in Lemma \ref{contiguity mean}, for all $u<v \in [-M,M]$, we have 
\begin{align*}
    & \left|   \Big[ \mb{E}_{\mu_n} \lb A_n(u) \rb -    \mb{E}_{\mu_n} \lb A_n(v) \rb \Big]   \right| \\
  = &  \left|   \Big[ \mb{E}_{\mu_n} \lb A_n(u) \rb -    \mb{E}_{\mu_n} \lb A_n(v) \rb \Big] - \underbrace{ \Big[ \mb{E}_{\mu_0} \lb A_n(u) \rb -    \mb{E}_{\mu_0} \lb A_n(v) \rb \Big]}_{=0}  \right| \\
  \leq & \sqrt{ \mb{E}_{\mu_0} \left[ A_n(u) - A_n(v) \right]^2 \cdot \mbox{Var}_{\mb{P}_{\mu_0}} \lb L_n \rb } \\
  \lesssim & \sqrt{\mb{E}_{\mu_0} \left[ A_n(u) - A_n(v) \right]^2} \leq \sqrt{ \mb{P}_{\mu_0} \lb u \leq  \sqrt{p} \bm{X}^\top \bm{Y} \leq v \rb }.
\end{align*}
Since the densities of $\sqrt{p}  \bm{X}^\top \bm{Y} $ under $\mu_0$ are uniformly bounded for large $p$ (the upper bound can be taken as $1/2\pi$ for $p \geq 3$), equicontinuity follows. In fact, the argument above also gives H\"older continuity of the sequence in \eqref{equicontinuity} with exponent $1/2$ and with the same H\"older constant. 

Now, since the sequence of functions in \eqref{equicontinuity} is equicontinuous and converges pointwise to $0$ on $[-M,M]$, the convergence is uniform. Thus, we have 
\begin{align*}
 & \limsup_{n \to \infty} \la \sup_{u \in \mb{R}} \left|    \mb{E}_{\mu_n} \lb A_n(u) \rb  + \frac{\tau^2}{\sqrt{2\pi}} \exp \lb -u^2/2 \rb \right|  \ra  \\
 \leq &   \limsup_{n \to \infty} \la \sup_{ |u|>M } \left|  \mb{E}_{\mu_n} \lb A_n(u) \rb +  \frac{\tau^2}{\sqrt{2\pi}} \exp \lb -u^2/2 \rb \right|  \ra  \\
\leq &  \limsup_{n \to \infty} \la \sup_{ |u|>M } \Big|  \mb{E}_{\mu_n} \lb A_n(u) \rb  \Big|  \ra +   \frac{\tau^2}{\sqrt{2\pi}} \exp \lb -M^2/2 \rb
\end{align*} 
for all $M>0$. 

It suffices to show that the first term in the last display goes to $0$ as $M \to \infty$. To see this, use the second estimate in Lemma \ref{contiguity mean} again to get
\begin{align*}
    \left| \mb{E}_{\mb{P}_{\mu_n}} \lb  A_n(u) \rb - \underbrace{\mb{E}_{\mb{P}_{\mu_0}} \lb  A_n(u) \rb}_{=0} \right| \leq \sqrt{ \mb{E}_{\mb{P}_{\mu_0}} \lb  A_n^2(u) \rb \cdot \mbox{Var}_{\mb{P}_{\mu_0}} \lb L_n \rb}.
\end{align*}
Consequently,
\begin{align*}
    \sup_{ |u|>M } \Big|  \mb{E}_{\mu_n} \lb A_n(u) \rb  \Big| &\lesssim  \sup_{|u|>M} \la \sqrt{\mb{E}_{\mb{P}_{\mu_0}} \lb  A_n^2(u) \rb } \ra \\
    &\lesssim    \sup_{|u|>M} \la  \sqrt{a_{n,u} \lb 1- a_{n,u} \rb }\ra
\end{align*}
where $a_{n,u}:= \mb{P}_{\mu_0} \lb \sqrt{p} \bm{X}^\top \bm{Y} \leq u \rb$. By \eqref{edgeworth} below, we have 
\[
    \sup_{|u|>M} \la  a_{n,u} \lb 1- a_{n,u} \rb \ra = \Phi(M) \lb 1 - \Phi(M) \rb + O(1/p).
\]
Therefore,
\[
  \limsup_{n \to \infty} \la \sup_{ |u|>M } \Big|  \mb{E}_{\mu_n} \lb A_n(u) \rb  \Big|  \ra \lesssim  \sqrt{\Phi(M) \lb 1 - \Phi(M) \rb}.
\]
The last term goes to $0$ as $M \to \infty$. The proof is completed.
$\hfill$ $\square$



\noindent \textbf{Proof of Proposition \ref{Kolmogorov low-rank}.}
Recall from \eqref{edgeworth1} that 
\begin{align} \label{edgeworth}
    \sup_{u \in \mb{R}} \left| \mb{P}_{\mu_0} \lb \sqrt{p} \bm{X}^\top \bm{Y} \leq u \rb - \Phi(u)  \right| = O \lb p^{-1} \rb
\end{align}
as $p \to \infty$.

To finish the proof, note that 
\begin{align*}
n \cdot \left[ \mb{P}_{\mu_k} \lb \sqrt{p} \bm{X}^\top \bm{Y} \leq u \rb - \mb{P}_{\mu_0} \lb \sqrt{p} \bm{X}^\top \bm{Y} \leq u \rb  \right] &= n \left[ \Phi \lb u\sqrt{\frac{k}{p}} \rb - \Phi \lb u \rb \right] +O \lb n/p \rb \\
&= -nu \lb 1 - \sqrt{k/p}  \rb \phi(u)(1+o(1))  \\
&+ O(n/p) \\
&\to -\frac{\tau}{2}u\phi(u)
\end{align*}
for all $u \in \mb{R}$. The convergence above is uniform whenever $p/n \to \infty$. The proof is complete. $\hfill$ $\square$


\subsection{Other technical results}
Let us start with a concentration inequality for degenerate U-processes, taken from \cite{cattaneo2024uniform}. Let $(S,\mathcal{S})$ be a measurable space, and let $X_1,\dots,X_n$ be i.i.d.\ $S$-valued random
variables with common law $P$. Let $\mathcal{F}$ be a pointwise measurable class of measurable
functions $f : S \times S \to \mathbb{R}$.

Define the (degenerate) degenerate $U$-process of order two by
\[
U_n(f)
  := \frac{2}{n(n-1)} \sum_{i < j}
  \Big\{
    f(X_i,X_j)
    - \mathbb{E}\big[ f(X_i,X_j) \mid X_i \big]
    - \mathbb{E}\big[ f(X_i,X_j) \mid X_j \big]
    + \mathbb{E}\big[ f(X_i,X_j) \big]
  \Big\}
\]
for $f \in \mathcal{F}$. Assume that
\begin{enumerate}
  \item Each $f \in \mathcal{F}$ is symmetric, i.e.\ $f(s_1,s_2) = f(s_2,s_1)$ for all $s_1,s_2 \in S$.
  \item There exists a measurable envelope $F : S \times S \to \mathbb{R}$ such that
        $|f(s_1,s_2)| \le F(s_1,s_2)$ for all $f \in \mathcal{F}$ and all $s_1,s_2 \in S$.
  \item For any probability measure $Q$ on $(S \times S,\mathcal{S}\otimes\mathcal{S})$ and $q \ge 1$,
        let
        \[
          \|f\|_{Q,q} := \big( \mathbb{E}_Q[|f|^q] \big)^{1/q}.
        \]
        Suppose that the envelope $F$ is VC-type in the sense that there exist constants
        $C_1 \ge e$ and $C_2 \ge 1$ such that, for all $\varepsilon \in (0,1]$,
        \[
          \sup_{Q} N\Big(
            \mathcal{F}, \|\cdot\|_{Q,2}, \varepsilon \|F\|_{Q,2}
          \Big)
          \;\le\;
          \Big(\frac{C_1}{\varepsilon}\Big)^{C_2}
        \]
        where the supremum is taken over all finite discrete probability measures $Q$ on $S \times S$.
\end{enumerate}
Under the three conditions above, we have

 \begin{lemma}[Lemma SA37 from \cite{cattaneo2024uniform}] \label{U-process concentration}
   Let $\sigma > 0$ be any deterministic quantity satisfying
\[
\sup_{f \in \mathcal{F}} \|f\|_{P,2} \;\le\; \sigma \;\le\; \|F\|_{P,2},
\]
and define the random variable
\[
M := \max_{1 \le i,j \le n} |F(X_i,X_j)|.
\]

Then there exists a universal constant $C_3 > 0$ such that
\[
n \, \mathbb{E}\Bigg[
  \sup_{f \in \mathcal{F}} |U_n(f)|
\Bigg]
\;\le\;
C_3 \sigma \Bigg(
  C_2 \log \frac{C_1 \|F\|_{P,2}}{\sigma}
\Bigg)
\;+\;
\frac{C_3 \|M\|_{P,2}}{\sqrt{n}}
\left[
  C_2 \log \lb \frac{C_1 \|F\|_{P,2}}{\sigma} \rb^{\!2}.
\right]
\]
 \end{lemma}

\noindent \textbf{Proof of Lemma \ref{U-process concentration}.} See \cite{cattaneo2024uniform}. $\hfill$ $\square$


\begin{lemma} \label{maximal invariant}
    Suppose $\bm{X}, \bm{Y}$ are matrices of size $p \times n$ such that $\bm{X}^\top \bm{X} =  \bm{Y}^\top \bm{Y}$. Then, there exists an orthogonal matrix $\bm{Q}$ of size $p \times p$ such that $\bm{X}= \bm{Q} \bm{Y}$.   
\end{lemma}

\noindent \textbf{Proof of Lemma \ref{maximal invariant}.} Notice that the assumption  $\bm{X}^\top \bm{X} =  \bm{Y}^\top \bm{Y}$ implies that $\mbox{Ker} \lb \bm{X} \rb =  \mbox{Ker} \lb \bm{Y} \rb$ since $\bm{X}$ and $\bm{X}^\top \bm{X}$ have the same kernel. We will use the notation $\mbox{Col} \lb \bm{A} \rb$ to indicate the column space of a matrix $\bm{A}$. Therefore, the map
\begin{align*}
    T: \mbox{Col} \lb \bm{X} \rb \subset \mb{R}^p &\to \mbox{Col} \lb \bm{Y} \rb  \subset \mb{R}^p  \\
 \bm{X}\bm{u} &\mapsto  \bm{Y} \bm{u} 
\end{align*}
is well-defined for all $\bm{u} \in \mb{R}^n$. Note that $ \mbox{dim} \lb \mbox{Col} \lb \bm{X} \rb  \rb = \mbox{dim} \lb \mbox{Col} \lb \bm{Y} \rb  \rb $ since their kernels are identical.

Moreover, $T$ is an isometry because
\[
\langle T\lb      \bm{X}\bm{u} \rb , T \lb      \bm{X}\bm{v} \rb \rangle =  \langle      \bm{Y}\bm{u},      \bm{Y}\bm{v} \rangle = \bm{u}^\top \bm{Y}^\top \bm{Y} \bm{v} = \bm{u}^\top \bm{X}^\top \bm{X} \bm{v} = \langle       \bm{X}\bm{u},      \bm{X}\bm{v} \rangle.
\]
Since $ \mbox{dim} \lb \mbox{Col} \lb \bm{X} \rb  \rb = \mbox{dim} \lb \mbox{Col} \lb \bm{Y} \rb  \rb $, $T$ admits a linear isometric extension to $\mb{R}^p$. Since linear isometries are represented by orthogonal matrices, we have
\[
\bm{Q} \bm{X} \bm{u} = \bm{Y} \bm{u}
\]
for some orthogonal matrix $\bm{Q}$ of size $p \times p$ and all $\bm{u} \in \mb{R}^n$. 
The proof is complete. $\hfill$ $\square$

 \begin{lemma} \label{degenerate-Ustat}
Let $\bm{X}$, $\bm{Y}$, and $\bm{Z}$ be i.i.d.\ realizations of the uniform distribution on $\mb{S}^{p-1}$, and let $f: \mb{R} \mapsto \mb{R}$ be a bounded measurable function. Then, we have 
   \begin{itemize}
       \item 
    $\mb{E} \lb f(\bm{X}^T \bm{Y}) \Big|  \bm{Y} \rb = \mb{E} f(\bm{X}^T \bm{Y}) $
    almost surely. 
    \item $\bm{X}^T \bm{Y}$ and $\bm{X}^T \bm{Z}$ are independent.
   \end{itemize}
    
\end{lemma}

\noindent\textbf{Proof of Lemma \ref{degenerate-Ustat}}. The first claim is a consequence of rotational invariance. Conditioning on $\bm{Y}$, there exists an orthogonal matrix $O$ such that $O^T \bm{Y} = \bm{e}_1= (1,0,\dots,0)$. Thus, with probability one, we have
\begin{align*}
    \mb{E} \lb f(\bm{X}^T \bm{Y}) \Big|  \bm{Y} \rb &= \mb{E} \lb f \lb \bm{X}^T O^T \bm{Y} \rb \Big|  \bm{Y} \rb \\
    &= \mb{E} \lb f(\bm{X}^T \bm{e}_1) \Big|  \bm{Y} \rb\\
    &= \mb{E} \lb f(\bm{X}^T \bm{e}_1) \rb ,
\end{align*}
where we use the fact that $\bm{X}^T \bm{e}_1$ is independent of $\bm{Y}$ in the last equality. Similarly, one can also show that $\mb{E} f(\bm{X}^T \bm{Y}) = \mb{E} \lb f(\bm{X}^T \bm{e}_1) \rb$. This concludes the proof of the first claim.

For the second claim, take any bounded measurable functions $f$ and $g$. By conditioning on $\bm{X}$, we get
\begin{align*}
    \mb{E} \lb f(\bm{X}^T \bm{Y}) g(\bm{X}^T \bm{Z}) \rb &= \mb{E} \Big[ \mb{E} \lb f(\bm{X}^T \bm{Y}) \Big| \bm{X} \rb \cdot \mb{E} \lb g(\bm{X}^T \bm{Z}) \Big| \bm{X} \rb  \Big] \\
    &= \mb{E} f(\bm{X}^T \bm{Y}) \cdot \mb{E} g(\bm{X}^T \bm{Z}),
\end{align*}
where we use the conclusion of the first statement in the last equality. This concludes the proof.
$\hfill\square$

\medskip

A direct consequence of Lemma \ref{degenerate-Ustat} is that, for any bounded measurable function $h$, the term $Z_n$ defined in (\ref{def-Zn}) is a degenerate U-statistic of order $1$ (see, for example, Section 5.3 of the monograph \cite{Dehling} for a comprehensive introduction to U-statistics and their limit theory). This fact was used frequently throughout the proof of Theorem 1. The next lemma is used in checking the conditions of the martingale CLT, which was used in the proof of Theorem 1. It gives a simpler form for the distribution of joint angles. 

\begin{lemma}\label{qi_hard} 
Let $p\geq 2$ be an integer. Assume $\bm{a}\in \mathbb{S}^{p-1}$ and $\bm{b} \in \mathbb{S}^{p-1}$ are fixed vectors. Let $\bm{X}$ be a random vector uniformly distributed over $\mathbb{S}^{p-1}$ and $\xi_1, \xi_2,\dots,\xi_p$ be i.i.d. standard normal. Set $\bm{\xi}=\lb \xi_1, \xi_2,\dots,\xi_p \rb^{\top}$. Then, for any bounded measurable function $f(x,y): \mathbb{R}^2 \to \mathbb{R}$, we have 
\[
\mb{E}f(\bm{a}^T \bm{X}, \bm{b}^T \bm{X})= \mb{E}f\Big(\frac{\xi_1}{\|\bm{\xi}\|}, (\bm{a}^T \bm{b})\frac{\xi_1}{\|\bm{\xi}\|}+ \sqrt{1-(\bm{a}^T \bm{b})^2}\,\frac{\xi_2}{\|\bm{\xi}\|}\Big).
\]
\end{lemma}

\noindent\textbf{Proof of Lemma \ref{qi_hard}}. 
We may assume
$\bm{X}=(\xi_1, \dots, \xi_p)^{\top}/\|\bm{\xi}\|$ without loss of generality. Suppose $ \bm{c}\in \mathbb{S}^{p-1}$ and $\bm{d}\in \mathbb{S}^{p-1}$ are constant vectors with $\bm{c}^{\top} \bm{d}=0$. Let $A$ be an orthogonal matrix whose first two rows are $\bm{c}^{\top}$ and $\bm{d}^{\top}$, respectively. Then, by the Haar invariance of the uniform distribution on the sphere, $A\bm{X}$ and $\bm{X}$ are identically distributed. In particular, the top two entries of $A\bm{X}$, that is, $(\bm{c}^{\top} \bm{X}, \bm{d}^{\top} \bm{X})$, have the same distribution as that of $(\xi_1, \xi_2)^{\top}/\|\bm{\xi}\|.$ Write
\[
\bm{b}=(\bm{a}^{\top}\bm{b})\bm{a}+ \sqrt{1-(\bm{a}^{\top}\bm{b})^2}\cdot \frac{\bm{b}-(\bm{a}^{\top}\bm{b})\bm{a}}{\sqrt{1-(\bm{a}^{\top}\bm{b})^2}}.
\]
The advantage of doing so is the trivial observation that $\bm{a}$ and $\frac{\bm{b}-(\bm{a}^{\top}\bm{b})\bm{a}}{\sqrt{1-(\bm{a}^{\top}\bm{b})^2}}$ are orthogonal unit vectors. Hence $\bm{a}^T \bm{X}$ and $\frac{\bm{b}-(\bm{a}^{\top}\bm{b})\bm{a}}{\sqrt{1-(\bm{a}^{\top}\bm{b})^2}} \bm{X}$ have the same law as that of $(\xi_1, \xi_2)'/\|\bm{\xi}\|.$ This implies that $(\bm{a}^T \bm{X}, \bm{b}^T \bm{X})$ has the same law as
\[
\Big(\frac{\xi_1}{\|\bm{\xi}\|}, (\bm{a}^{\top}\bm{b})\frac{\xi_1}{\|\bm{\xi}\|}+ \sqrt{1-(\bm{a}^{\top}\bm{b})^2}\,\frac{\xi_2}{\|\bm{\xi}\|}\Big).
\]
As a result,
\[
\mb{E}f(\bm{a}^T \bm{X}, \bm{b}^T \bm{X})= \mb{E}f\Big(\frac{\xi_1}{\|\bm{\xi}\|}, (\bm{a}^T \bm{b})\frac{\xi_1}{\|\bm{\xi}\|}+ \sqrt{1-(\bm{a}^T \bm{b})^2}\,\frac{\xi_2}{\|\bm{\xi}\|}\Big),
\]
where the last expectation is taken over $\xi_1, \dots, \xi_p$, and hence it is a function of $\bm{a}^T \bm{b}$. The proof is completed. $\hfill\square$
\medskip

\begin{lemma} \label{Diaconis-Freeman}
    Let $H_n$ be defined in (\ref{def-H}) for a sequence of measurable functions $h_n: \mb{R} \mapsto \mb{R}$. Assume additionally that 
    \begin{align*}
          \mb{E} h_n(\bm{X}^{\top}_1 \bm{X}_2) &=0; \\
          \mbox{Var} \lb h_n(\bm{X}^{\top}_1 \bm{X}_2) \rb &\leq C_1
    \end{align*}
for some constant $C_1$ independent of $n$. Then, we have
    $$\mb{E} H_n^2(\bm{X}_1, \bm{X}_2) \to 0$$
    as $n \to \infty$.
\end{lemma}

\noindent\textbf{Proof of Lemma \ref{Diaconis-Freeman}}. It suffices to prove Lemma \ref{Diaconis-Freeman} when $h_n$ is bounded. Indeed, suppose we have proved $\mb{E} H_n^2(\bm{X_1}, \bm{X}_2) \to 0$ for all bounded $h_n$. Then write 
\[
h_n = \underbrace{h_n \cdot \mathbf{1}_{\la |h_n| \leq L \ra} - \mb{E} \lb h_n \cdot \mathbf{1}_{\la |h_n| \leq L \ra} \rb}_{f_{n,L}} + \underbrace{h_n \cdot \mathbf{1}_{\la |h_n| > L \ra} - \mb{E} \lb h_n \cdot \mathbf{1}_{\la |h_n| > L \ra} \rb}_{g_{n,L}}
\]
where the expectation is taken with respect to the law of $\bm{X}_1^\top \bm{X}_2$. Then,
\begin{align*}
        \mb{E} H_n^2(\bm{X_1}, \bm{X}_2) &= \mb{E} \mb{E}^2 \Big[  h_n \lb \bm{X}_1^T \bm{Y} \rb \cdot h_n \lb \bm{X}_2^T \bm{Y}  \rb \Big| \bm{X}_1, \bm{X}_2 \Big] \\
        &= \mb{E} \left[ \mb{E}^2 \left[  f_{n,L} \lb \bm{X}_1^T \bm{Y} \rb \cdot f_{n,L} \lb \bm{X}_2^T \bm{Y}  \rb \Big| \bm{X}_1, \bm{X}_2 \right]  \right] + O \lb \mbox{Var} \lb g_{n,L} \lb (\bm{X_1}, \bm{X}_2 \rb \rb \rb.
\end{align*}
For every fixed $L$, the first term tends to zero as $n \to \infty$. We then deduce the result by letting $L \to \infty$ and noting that
\[
\sup_{n \geq 1} \mbox{Var} \lb g_{n,L} \lb (\bm{X_1}, \bm{X}_2 \rb \rb \leq \frac{C_1}{L^2}.
\]
Now assume that $h$ is bounded. Let $\bm{Y}$ be drawn from the uniform distribution on $\mb{S}^{p-1}$ independently of $\bm{X}_1$ and $\bm{X}_2$. Thanks to Lemma \ref{qi_hard}, we can write
\begin{align*}
    \mb{E} H_n^2(\bm{X_1}, \bm{X}_2) &= \mb{E} \mb{E}^2 \Big[  h_n \lb \bm{X}_1^T \bm{Y} \rb \cdot h_n \lb \bm{X}_2^T \bm{Y}  \rb \Big| \bm{X}_1, \bm{X}_2 \Big] \\
    &= \mb{E} \mb{E}^2 \bigg[  h_n \lb \frac{\xi_1}{\| \bm{\xi} \|} \rb \cdot h_n \lb \bm{X}_1^T \bm{X}_2 \cdot \frac{\xi_1}{\| \bm{\xi} \|} + \sqrt{1 - (\bm{X}_1^T \bm{X}_2)^2} \cdot \frac{\xi_2}{\| \bm{\xi} \|} \rb \Big| \bm{X}_1, \bm{X}_2 \bigg],
\end{align*}
where  $\bm{\xi}=(\xi_1, \xi_2,\dots,\xi_p)^{\top}$ is a vector with i.i.d. standard normal components. Set $U= \bm{X}_1^T \bm{X}_2$ and let $f(u)$ be the density of $U$. We can write
\begin{align*}
    \mb{E} H_n^2(\bm{X_1}, \bm{X}_2) &= \int_{-1}^{1} f(u) \cdot \mb{E}^2 \bigg[   h_n \lb \frac{\xi_1}{\| \bm{\xi} \|} \rb \cdot h_n \lb u \cdot \frac{\xi_1}{\| \bm{\xi} \|} + \sqrt{1 - u^2} \cdot \frac{\xi_2}{\| \bm{\xi} \|} \rb  \bigg] du \\
    &= \int_{|u| \leq \ve}  f(u) \cdot \mb{E}^2 \bigg[   h_n \lb \frac{\xi_1}{\| \bm{\xi} \|} \rb \cdot h_n \lb u \cdot \frac{\xi_1}{\| \bm{\xi} \|} + \sqrt{1 - u^2} \cdot \frac{\xi_2}{\| \bm{\xi} \|} \rb  \bigg] du \\
    &+ \int_{|u| > \ve}  f(u) \cdot \mb{E}^2 \bigg[   h_n \lb \frac{\xi_1}{\| \bm{\xi} \|} \rb \cdot h_n \lb u \cdot \frac{\xi_1}{\| \bm{\xi} \|} + \sqrt{1 - u^2} \cdot \frac{\xi_2}{\| \bm{\xi} \|} \rb  \bigg] du \\
    & \leq \int_{|u| \leq \ve}  f(u) \cdot \mb{E}^2 \bigg[   h_n \lb \frac{\xi_1}{\| \bm{\xi} \|} \rb \cdot h_n \lb u \cdot \frac{\xi_1}{\| \bm{\xi} \|} + \sqrt{1 - u^2} \cdot \frac{\xi_2}{\| \bm{\xi} \|} \rb  \bigg] du \\
    &+ \| h_n \|_{\infty}^2 \cdot \mb{P} \lb |\bm{X}_1^T \bm{X}_2| > \ve \rb
\end{align*}
for any fixed $\ve>0$. By Proposition 5 in \cite{Jiang13}, we get $\mb{P} \lb |\bm{X}_1^T \bm{X}_2| > \ve \rb \to 0$ as $p \to \infty$, and hence it suffices to bound the first integrand over $(-\ve, \ve)$. Thanks to Lemma \ref{Diaconis-Freeman2}, for $|u| \leq \ve$ we get the bound
\begin{align*}
   & \bigg| \mb{E} \Big[   h_n \lb \frac{\xi_1}{\| \bm{\xi} \|} \rb \cdot h_n \lb u \cdot \frac{\xi_1}{\| \bm{\xi} \|} + \sqrt{1 - u^2} \cdot \frac{\xi_2}{\| \bm{\xi} \|} \rb  \Big] - \mb{E} \Big[   h_n \lb \frac{\xi_1}{\sqrt{p}} \rb \cdot h_n \lb   \frac{\xi_2}{\sqrt{p}} \rb \Big] \bigg| \\
   \leq & (\mbox{const}) \cdot \| h_n \|_{\infty}^2 \cdot \lb p^{-1} + \ve \rb,
\end{align*}
by choosing $g_n(x,y) =   h_n(x) h_n(y)$ in Lemma \ref{Diaconis-Freeman2}. Note that $\mb{E} h_n \lb \xi_1 / \| \bm{\xi} \| \rb =0$ since $\mb{E} h_n(\bm{X}^{\top}_1 \bm{X}_2)=0$, which gives us
\begin{align*}
    \bigg| \mb{E} \Big[   h_n \lb \frac{\xi_1}{\sqrt{p}} \rb \cdot h_n \lb   \frac{\xi_2}{\sqrt{p}} \rb \Big]  \bigg| &= \bigg|  \mb{E} \Big[   h_n \lb \frac{\xi_1}{\sqrt{p}} \rb \Big] \bigg|^2 \\
    &= \bigg| \mb{E} \Big[   h_n \lb \frac{\xi_1}{\sqrt{p}} \rb \Big] - \mb{E} \Big[   h_n \lb \frac{\xi_1}{\| \bm{\xi} \|} \rb \Big]  \bigg|^2 \\
    & \leq (\mbox{const}) \cdot \frac{ \| h_n \|_{\infty}^2}{p},
\end{align*}
where we use the bound (\ref{tv1}) in the last inequality. This in turn yields
\begin{align*}
    \bigg| \mb{E} \Big[   h_n \lb \frac{\xi_1}{\| \bm{\xi} \|} \rb \cdot h_n \lb u \cdot \frac{\xi_1}{\| \bm{\xi} \|} + \sqrt{1 - u^2} \cdot \frac{\xi_2}{\| \bm{\xi} \|} \rb  \Big] \bigg| \leq (\mbox{const}) \cdot \lb p^{-1} + \ve \rb,
\end{align*}
since the $L_{\infty}$ norm of $h_n$ is uniformly bounded. Consequently,
\begin{align*}
   & \int_{|u| \leq \ve}  f(u) \cdot \mb{E}^2 \bigg[   h_n \lb \frac{\xi_1}{\| \bm{\xi} \|} \rb \cdot h_n \lb u \cdot \frac{\xi_1}{\| \bm{\xi} \|} + \sqrt{1 - u^2} \cdot \frac{\xi_2}{\| \bm{\xi} \|} \rb  \bigg] du\\
    \leq & (\mbox{const}) \cdot \int_{|u| \leq \ve} f(u) \cdot (p^{-1} + \ve) du.
\end{align*}
The proof is completed by taking $p \to \infty$ and then letting $\ve \to 0$. $\hfill\square$
\medskip

\begin{lemma} \label{Diaconis-Freeman2}
   Let $\bm{\xi}=(\xi_1, \xi_2,\dots,\xi_p)^{\top}$ be a vector consisting of i.i.d. standard normal variables. Then, for any bounded, measurable function $g: \mb{R}^2 \mapsto \mb{R}$ and $u \in [-1/2,1/2]$, we have
   $$ \Big| \mb{E} g \lb \frac{ \xi_1}{\| \bm{\xi} \|}, u \cdot \frac{ \xi_1}{\| \bm{\xi} \|} + \sqrt{1-u^2} \cdot \frac{ \xi_2}{\| \bm{\xi} \|}  \rb - \mb{E} g \lb \frac{\xi_1}{\sqrt{p}}, \frac{\xi_2}{\sqrt{p}}  \rb \Big| \leq C \| g \|_{\infty} \lb \frac{1}{p} + |u| \rb $$
   for some universal constant $C>0$.
\end{lemma}

\noindent\textbf{Proof of Lemma \ref{Diaconis-Freeman2}}. The conclusion follows from the following two total variation distance bounds:
\begin{align} \label{tv1}
    d_{TV} \lb  \lb \frac{\sqrt{p} \cdot \xi_1}{\| \bm{\xi} \|}, \frac{\sqrt{p} \cdot \xi_2}{\| \bm{\xi} \|}  \rb, \lb \xi_1, \xi_2 \rb \rb \leq  \frac{C}{p} 
\end{align}
and 
\begin{align} \label{tv2}
    d_{TV} \Big( N\left(\begin{pmatrix}
0 \\
0
\end{pmatrix},\begin{pmatrix}
1 & u \\
u & 1
\end{pmatrix}\right) , N( \bm{0}, \bm{I}_2)
\Big) \leq C|u|
\end{align}
for some universal constant $C$ and for all $p$ large enough in the first estimate. We next explain \eqref{tv1} and \eqref{tv2}.

The first estimate (\ref{tv1}) is a consequence of the Diaconis--Freedman theorem (see Theorem 2.8 in the monograph \cite{Meckes} and also the paper \cite{D-F}). It provides a sharp bound in total variation between the joint distribution of the first few entries of  $\mbox{Unif} \lb \mb{S}^{p-1} \rb$ and the standard multivariate normal random vector of the same length. The second estimate (\ref{tv2}) is elementary and can be proved directly by estimating the difference between the two corresponding densities. To see this, write
\begin{align*}
      &  d_{TV} \Big( N\left(\begin{pmatrix}
0 \\
0
\end{pmatrix},\begin{pmatrix}
1 & u \\
u & 1
\end{pmatrix}\right) , N( \bm{0}, \bm{I}_2)
\Big) \\
\leq &  \int_{\mb{R}} \int_{\mb{R}} \bigg| \frac{1}{2 \pi} \exp \la -\frac{x^2+y^2}{2} \ra - \frac{1}{2 \pi \sqrt{1-u^2} } \exp \la - \frac{x^2 + y^2 - 2uxy}{2(1-u^2)} \ra  \bigg| dxdy \\
\leq & \int_{\mb{R}} \int_{\mb{R}} \frac{1}{2 \pi} \exp \la -\frac{x^2+y^2}{2} \ra \cdot  \bigg|  1 - \frac{1}{\sqrt{1-u^2}} \exp \la - \frac{(x^2+y^2)u^2}{2(1-u^2)} +\frac{uxy}{1-u^2} \ra  \bigg| dxdy \\
 \leq & I_1 + I_2,
\end{align*}
where 
\begin{align*}
    I_1 &=   \bigg| 1- \frac{1}{\sqrt{1-u^2}} \bigg| \cdot \int_{\mb{R}} \int_{\mb{R}} \frac{1}{2 \pi} \exp \la -\frac{x^2+y^2}{2} \ra dxdy, \\
    I_2 &=  \frac{1}{2 \pi \sqrt{1-u^2}} \cdot  \int_{\mb{R}} \int_{\mb{R}} \exp \la -\frac{x^2+y^2}{2} \ra \cdot \bigg| 1 -  \exp \la - \frac{(x^2+y^2)u^2}{2(1-u^2)} +\frac{uxy}{1-u^2} \ra  \bigg|. 
\end{align*}
The term $I_1$ is obviously of order $O(|u|)$ as $u \to 0$, and thus we only need to bound $I_2$. In polar coordinates, $I_2$ can be rewritten as
\begin{align*}
    I_2 &=   \frac{1}{2 \pi \sqrt{1-u^2}} \cdot \int_{0}^{\infty} r \cdot \exp \la - \frac{r^2}{2}  \ra \cdot \int_{0}^{2 \pi} \bigg| 1- \exp \la - \frac{r^2u^2 - u \sin 2 \theta}{2(1-u^2)} \ra \bigg| d \theta dr.
\end{align*}
Moreover, we have
\begin{align*}
     \bigg| 1- \exp \la - \frac{r^2u^2 - u \sin 2 \theta}{2(1-u^2)} \ra \bigg| & \leq \bigg| 1 - \exp \la \frac{u \sin 2 \theta}{2(1-u^2)} \ra  \bigg| \\
     &+ \bigg| \exp \la \frac{u \sin 2 \theta}{2(1-u^2)} \ra  \bigg| \cdot \bigg| 1 - \exp \la - \frac{r^2 u^2}{2(1-u^2)} \ra \bigg| \\
     & \leq (\mbox{const}) \cdot \bigg[ \Big| \frac{u \sin 2 \theta}{2(1-u^2)} \Big| + \frac{r^2 u^2}{2(1-u^2)} \bigg],
\end{align*}
where we use the elementary inequalities $1-e^{-x} \leq x$ for all $x>0$ and $|1-e^{x}| \leq C_1 |x|$ for all $|x| \leq C$, where $C_1$ depends only on $C$. Thus, we have
\begin{align*}
    I_2 & \leq (\mbox{const}) \cdot \frac{|u|}{4 \pi (1-u^2)^{3/2}} \cdot  \int_{0}^{\infty} r \cdot \exp \la - \frac{r^2}{2}  \ra \cdot \bigg( \int_{0}^{2 \pi} |\sin 2 \theta| + r^2u d\theta \bigg) dr \\
    & \leq (\mbox{const}) |u|.
\end{align*}
This concludes the proof of (\ref{tv2}). 

Now we are ready to prove the main estimate in Lemma \ref{Diaconis-Freeman2}. Write 
\begin{align*}
   & \Big| \mb{E} g \lb \frac{ \xi_1}{\| \bm{\xi} \|}, u \cdot \frac{ \xi_1}{\| \bm{\xi} \|} + \sqrt{1-u^2} \cdot \frac{ \xi_2}{\| \bm{\xi} \|}  \rb - \mb{E} g \lb \frac{\xi_1}{\sqrt{p}}, \frac{\xi_2}{\sqrt{p}}  \rb \Big| \\
    \leq & \Big| \mb{E} g \lb \frac{ \xi_1}{\| \bm{\xi} \|}, u \cdot \frac{ \xi_1}{\| \bm{\xi} \|} + \sqrt{1-u^2} \cdot \frac{ \xi_2}{\| \bm{\xi} \|}  \rb - \mb{E} g \lb \frac{\xi_1}{\sqrt{p}}, u \cdot \frac{\xi_1}{\sqrt{p}} + \sqrt{1-u^2} \cdot \frac{\xi_2}{\sqrt{p}}  \rb \Big| \\
    + & \Big| \mb{E} g \lb \frac{ \xi_1}{\sqrt{p}}, u \cdot \frac{ \xi_1}{\sqrt{p}} + \sqrt{1-u^2} \cdot \frac{ \xi_2}{\sqrt{p}}  \rb - \mb{E} g \lb \frac{\xi_1}{\sqrt{p}}, \frac{\xi_2}{\sqrt{p}}  \rb \Big| \\
    \leq & \frac{(\mbox{const}) \cdot \|g\|_{\infty}}{p} + (\mbox{const}) \cdot  \| g \|_{\infty} |u|.
\end{align*}
The proof is completed. $\hfill\square$
\medskip

The next lemma collects some elementary properties and bounds for the normalizing constant $C_{p}(\kappa)$ of the FvML distributions.
\begin{lemma} \label{Bessel}
    Recall $C_{p}(\kappa)$ in \eqref{Cp(kappa)}. The following statements hold:
    \begin{itemize}
        \item If $I_{\nu}(x)$ is the modified Bessel function of the first kind (see \cite{Ley-Verdebout} for details), then
        \[ \frac{C_{p}(\kappa)^{'}}{C_p(\kappa)} = - \frac{I_{p/2}(\kappa)}{I_{p/2-1}(\kappa)}. \]
        \item For all $\nu>0$, we have 
        \[
         G_{\nu+ \frac{1}{2},\nu+ \frac{3}{2}} (\kappa)  \leq \frac{I_{\nu+1}(\kappa)}{I_{\nu}(\kappa)}   \leq  G_{\nu,\nu+2} (\kappa)
        \]
        where
        \[
        G_{\alpha,\beta}(t):= \frac{t}{\alpha + \sqrt{\beta^2+t^2}}.
        \]
        \item For all $\kappa>0$, we have 
        \[
        \left| \log C_p(\kappa) + \frac{\kappa^2}{2p}  \right| \leq \frac{\kappa^4}{2p^4}
        \]
        for all $p \geq 3$ and $\kappa >0$.
    \end{itemize}
\end{lemma}

\noindent \textbf{Proof of Lemma \ref{Bessel}.}   
The first property can be found in \cite{M-Jupp}, pages $169$ and $170$. The second property can be found in Section 3 of \cite{hornik2013amos}. Let us use the first two properties to prove the last one.

We first show that
\begin{align} \label{Bessel bound}
\Big|  \frac{I_{p/2}(\kappa)}{I_{p/2-1}(\kappa)} - \frac{\kappa}{p}  \Big| \leq \frac{2 \kappa^3}{p^3}.
\end{align}

To see this, we apply the second property with $\nu= (p-2)/2$ to obtain 
\[
\frac{\kappa}{\nu+1/2 + \sqrt{\kappa^2 + \lb \nu + 3/2 \rb^2}}  \leq \frac{I_{p/2}(\kappa)}{I_{p/2-1}(\kappa)} \leq \frac{\kappa}{\nu + \sqrt{\kappa^2 + \lb \nu +2 \rb^2}}
\]
From the upper bound, it is clear that the Bessel ratio is always less than $1$. Consider the following two cases.

{\it  \underline{Case 1: $\kappa \geq p$}.} In this case, the result is trivial since 
\[\Big|  \frac{I_{p/2}(\kappa)}{I_{p/2-1}(\kappa)} - \frac{\kappa}{p}  \Big|  =  \frac{\kappa}{p} - \frac{I_{p/2}(\kappa)}{I_{p/2-1}(\kappa)} \leq \frac{\kappa}{p} \leq \frac{\kappa^3}{p^3}.
\]

{\it \underline{Case 2: $\kappa \leq p$}.} In this case, put $x:= \kappa/p \leq 1$. Use the lower bound to get 
\begin{align*}
      x - \frac{I_{p/2}(\kappa)}{I_{p/2-1}(\kappa)} &\leq x - \frac{\kappa}{(p-1)/2 + \sqrt{\kappa^2 + \lb (p+1)/2 \rb^2}} \\
      &= x \left[ 1 - \frac{1}{1+ \sqrt{x^2 + \lb 1/2 + 1/(2p) \rb^2} - 1/2-1/(2p)}  \right] \\
      &= x \cdot \frac{  \sqrt{x^2 + \lb 1/2 + 1/(2p) \rb^2} - [1/2+1/(2p)] }{ 1+ \sqrt{x^2 + \lb 1/2 + 1/(2p) \rb^2} - 1/2-1/(2p) } \\
      &\leq \frac{2x^3}{ 1+ \sqrt{x^2 + \lb 1/2 + 1/(2p) \rb^2} - 1/2-1/(2p)} \leq 2x^3
\end{align*}
where the last line follows from the fact that 
\[
  \sqrt{x^2 + \lb 1/2 + 1/(2p) \rb^2} - [1/2+1/(2p)] = \frac{x^2}{  \sqrt{x^2 + \lb 1/2 + 1/(2p) \rb^2} + [1/2+1/(2p)] } \leq 2x^2.
\]
Finally, to deduce the third property, we integrate \eqref{Bessel bound} to get
\begin{align*}
    \left| \log C_p(\kappa) + \frac{\kappa^2}{2p}  \right| = \left| \int_0^{\kappa} \lb  - \frac{I_{p/2}(t)}{I_{p/2-1}(t)} + \frac{t}{p} \rb dt\right| &\leq \int_0^{\kappa} \left|   - \frac{I_{p/2}(t)}{I_{p/2-1}(t)} + \frac{t}{p}  \right| dt \\
    &\leq \frac{\kappa^4}{2p^4}.
\end{align*}
This completes the proof. $\hfill$ $\square$

\end{appendix}

\bibliographystyle{plainnat}
\bibliography{paper-ref}

 \end{document}